\documentclass[oneside,reqno,11pt]{article}
\usepackage[margin=0.9in]{geometry}
\usepackage[nottoc]{tocbibind}

\usepackage{amssymb}
\usepackage{amsmath}
\usepackage{amsthm}      
\usepackage{graphicx} 
\usepackage{tikz}
\usepackage{tikz-cd}
\usetikzlibrary{matrix}
\usetikzlibrary{arrows.meta}
\usetikzlibrary{calc}
\usepackage[capitalise]{cleveref}
\usepackage{setspace}
\usepackage{verbatim}
\usepackage{xfrac}
\usepackage{pgffor}
\usepackage{calc}
\usepackage{subcaption}

\usepackage{enumerate}

\usepackage{array}   
\newcolumntype{C}{>{$}c<{$}}

\newtheorem{thm}{Theorem}[subsection]
\newtheorem*{thm*}{Theorem}
\newtheorem{lem}[thm]{Lemma}
\newtheorem{prop}[thm]{Proposition}

\theoremstyle{definition}
\newtheorem{defn}[thm]{Definition}
\newtheorem*{defn*}{Definition}
\theoremstyle{definition}

\theoremstyle{remark}

\newtheorem{rem}[thm]{Remark}

\crefname{defn}{Definition}{Definitions}
\crefname{figure}{Figure}{Figures}

\newcommand{\Emb}{\mathsf{Emb}}
\newcommand{\Embp}{\mathsf{Emb}_\partial}

\newcommand{\EmbIM}{\mathsf{Emb}_\partial(I, M)}
\newcommand{\ImmIM}{\mathsf{Imm}_\partial(I, M)}
\newcommand{\EmbIS}{\mathsf{Emb}_\partial(I, S^1\times B^3)}
\newcommand{\EmbIB}{\mathsf{Emb}_\partial(I, B^4)}
\newcommand{\Gpqr}{G(p, q, r)}
\renewcommand{\dim}{\mathsf{dim}}

\newcommand{\C}{\mathcal{C}}
\newcommand{\Top}{\mathsf{Top}}
\renewcommand{\S}{\mathcal{S}}

\newcommand{\prl}{\mathbin{||}}

\newcommand{\CSB}[1]{C_{#1}\langle S^1\times B^3 \rangle}
\newcommand{\twij}[4]{t_{#1}^{#2}\cdot w_{#3 #4}}

\newcommand{\sans}{\mathbin{\backslash}}
\newcommand{\lr}[1]{\langle #1 \rangle}

\newcommand{\To}{\xrightarrow}

\newcommand{\coloneq}{\colon\!\!\!=}
\newcommand{\Tto}[3][20pt]{\begin{tikzcd}[sep=#1, cramped, ampersand replacement=\&, text height=1ex, text depth=.3ex]\ar[r, shift left=0pt, "#2"]\ar[r, shift right=0pt, "#3"']\&{}\end{tikzcd}}
\newcommand{\Q}{\mathbb{Q}}
\renewcommand{\sb}[2]{\left\{#1\;\middle|\;#2\right\}}

\newcommand{\sand}{\; \textrm{and} \;}
\newcommand{\qand}{\quad \textrm{and} \quad}
\newcommand{\qqand}{\qquad \textrm{and} \qquad}

\newcommand{\cubeop}[1]{\mathbf{#1}}
\newcommand{\id}{\cubeop{id}}
\newcommand{\rev}{\cubeop{rev}}
\newcommand{\fold}{\cubeop{fold}}
\newcommand{\rot}{\cubeop{rot}}
\newcommand{\face}{\cubeop{face}}
\newcommand{\twist}{\cubeop{twist}}

\newcommand{\compass}[2]{\begin{tikzcd}[ampersand replacement=\&,row sep=huge, column sep=large] {} \\ \color{white}{\bullet} \uar{#2} \rar[swap]{#1} \& {} \end{tikzcd}}
\newcommand{\compassss}[2]{\begin{tikzcd}[ampersand replacement=\&,row sep=huge, column sep=large] {} \\ \color{white}{\bullet} \uar{#2} \rar[swap]{#1} \& {} \\ \\ \\ \\ \\ \end{tikzcd}}
\newcommand{\compasss}[2]{\begin{tikzcd}[ampersand replacement=\&,row sep=huge, column sep=large] \\ \\ \\ \\ \\ \\ \\ {} \\ \color{white}{\bullet} \uar{#2} \rar[swap]{#1} \& {} \end{tikzcd}}

\newcommand{\outcompass}[3]{\begin{tikzcd}[ampersand replacement=\&,row sep=huge, column sep=large] \&[-20] {} \\ \& \color{white}{\bullet} \ar[shorten <=-5]{dl}[swap]{#3} \uar{#2} \rar[swap]{#1} \& {} \\[-25] {} \end{tikzcd}}

\newcommand{\revolve}[2]{{\begin{array}{|c c c|}
     \hline  &  & \\
     \;\;\;\; & \;\boxed{\;#1\;}\; & \;#2\;\; \\
     & & \\ 
     \hline
    \end{array}}}

\newcommand{\flip}[1]{{\begin{array}{|cc|}
        \hline &
        {#1}\!\!\!\!
        \\[-\arraystretch ex] & \\ \hline
    \end{array}}}

\title{Families of Arcs in 4-Manifolds and Maps of Configuration Spaces}

\author{Shruthi Sridhar-Shapiro}

\date{}

\begin{document}

\maketitle

\begin{abstract}
In this thesis we construct 3-parameter families $\Gpqr$ of embedded arcs with fixed boundary in a 4-manifold. We then analyze these elements of $\pi_3\EmbIM$ using embedding calculus by studying the induced map from the embedding space to ``Taylor approximations" $T_k\EmbIM$. We develop a diagrammatic framework inspired by cubical $\omega$-groupoids to depict $\Gpqr$ and related homotopies. We use this framework extensively in Chapter 4 to show explicitly that $\Gpqr$ is trivial in $\pi_3T_3\EmbIM$ (however, we conjecture that it is non-trivial in $\pi_3T_4\EmbIM$). In Chapter 5 we use the Bousfield-Kan spectral sequence for homotopy groups of cosimplicial spaces to show that the rational homotopy group $\pi^\Q_3\EmbIS$ is $\Q$. This thesis extends work by Budney and Gabai in \cite{budney} which proves analogous results for $\pi_2\EmbIM$.
\end{abstract}

\tableofcontents

\section{Introduction}\label{subsec: intro}

$M$ will denote a smooth, connected, compact 4-manifold with boundary with two specified points $\ast_0, \ast_1 \in \partial M$. We specify an outgoing (from $\partial M$) unit vector $v_0$ at $\ast_0$ and an incoming unit vector $v_1$ at $\ast_1$. $I$ will denote the unit interval $[0,1]$. 

In this thesis, we study the embedding space $\EmbIM$ defined below. 
\begin{defn}
    $\EmbIM$ is the space of $C^1$ embeddings of $I$ into $M$ with constant speed such that $0\mapsto \ast_0$ and $1\mapsto \ast_1$, and unit tangent vectors at those points are $v_0, v_1$ respectively.  
\end{defn}

\begin{defn}
$\gamma$ is a chosen interval in $M$ which has endpoints $\ast_0, \ast_1$ which will serve as the base point in $\EmbIM$. 
\end{defn}

In \cite{budney}, Budney and Gabai construct non trivial elements of $\pi_k(\EmbIM)$  based at $\gamma$ for $k=1, 2$ and $M=S^1\times B^3$. This thesis builds on their work and in \cref{sec.gpqr} we construct 3-parameter families of embeddings $\Gpqr\colon I^3 \to \EmbIM$ that map to the base loop $\gamma$ on the boundary. 
 
Work of Goodwillie, Klein, Weiss in \cite{goodwillie} describe highly connected approximations --- $T_k\mathsf{Emb}_\partial (I, M)$ for the embedding space we are studying --- and show that $\pi_3\EmbIM$ is isomorphic to $\pi_3 T_4\EmbIM$ and surjective onto $\pi_3T_3\EmbIM$. Sinha in \cite{sinha} proves that that $T_n\EmbIM$ is homotopy equivalent to the space of strata preserving, \emph{aligned} maps between compactified $n$-point configuration spaces of $I$ and $M$. We write this as $$T_n\EmbIM \simeq Map^{sp}(C_n'\lr{I}, C_n'\lr{M})$$
From work in \cite{budney}, the element in $\pi_3T_3\EmbIM$ induced by $\Gpqr$ is torsion. In \cref{sec.nullhomotopy} we show (with an explicit homotopy) that this element is trivial in $\pi_3T_3\EmbIM$ using the mapping space model for $T_k\EmbIM$ that Sinha defines. \begin{thm}\label{thm: null homotopy}
    The map $T_3\Gpqr\colon I^3 \to \text{ Map}^{sp}(C_3\lr{I}, C_3 \lr{M})$ is homotopic to the map induced by the constant map $\id_\gamma\colon I^3 \to \EmbIM$.
\end{thm}
We conjecture that these elements are non-trivial in $\pi_3T_4\EmbIM$. 

Compared to the constructions in \cite{budney}, the maps we define have up to two additional parameters with a maximum of 7, making purely geometric definitions generally infeasible. To alleviate some of this dimensional burden, we develop notation and diagrams for concatenations and other operations on maps into a space from arbitrarily-high dimensional cubes. Inspired by the theory of cubical $\omega$-groupoids,\footnote{See for instance \cite{browncubes}.} these operations allow us to construct maps from higher dimensional cubes in an algebraic manner from geometrically defined building blocks, and our diagrams permit us to represent concatenations of high dimensional cubes in up to 4 directions at once using only 2-dimensional pictures. It is our hope that the use of this formalism makes our constructions more easily replicable for the reader, and that the new operations we define are of independent interest.

One strategy to show an element like $\Gpqr$ is non trivial is to look at the image of the induced map between configuration spaces into $\pi_7C_4\lr{M}$. To be precise, we look at the map $$\Gpqr^\ast\colon \text{ Map}^{sp}(I^3\times C_4\lr{I} , C_4\lr{M}) \to \pi_7C_4\lr{M}/R $$ 
The superscript $sp$ denotes strata preserving. The relations $R$ are necessary to define a well defined map to $\pi_7C_4\lr{M}$ when quotienting out the boundary of $I^3\times C_4\lr{I} \cong I^7$. 

In \cref{sec.piCk} we compute this group $\pi_7(\CSB{4})/R$ to be $\mathbb{Q}$ rationally where $R$ is subgroup generated by the relations obtained from 5 inclusions of $\pi_7(\CSB{3})$ into $\pi_7(\CSB{4})$ induced by the 5 face inclusions $\CSB{3} \hookrightarrow \CSB{4}$.

\begin{thm}\label{thm: pi7c4}
$\pi_7\CSB{4}/R\cong \mathbb{Q}$ rationally and is generated  by $[w_{12}, [w_{13}, w_{14}]]$
\end{thm}

This group being non trivial would allow us to potentially show that $\Gpqr$ is non trivial. \cref{sec.progress} we give some strategies to create invariants to prove the conjecture that $G(0,0,0)$ is the generator of $\pi_3 \EmbIS$. 

In \cite{sinha}, Sinha shows that $\EmbIM$ is homotopy equivalent to the totalization of a certain cosimplical space involving $C_i\lr{M}$. They use this to define a Bousefield-Kan spectral sequence that converges to the homotopy groups of $\EmbIM$. A related spectral sequence for homology of $\EmbIM$ has been shown to converge on the $E_2$ page when $M=B^4$ in \cite{lambrechts}. In \cite{sinha-scannell}, they use the above mentioned spectral sequence to compute $\pi_3\EmbIB$. We compute $\pi_3\EmbIS$ in \cref{subsec: relationsspectral} and show that that the map $\pi_3\EmbIB \to \pi_3\EmbIS$ is an isomorphism rationally giving the following theorem in \cref{subsec: relationsspectral}. 

\begin{thm}\label{thm: pi3EmbS1B3}
For rational homotopy groups, $$\pi_3\EmbIS \cong \pi_3\EmbIB \cong \Q.$$
\end{thm}

\section{Background}\label{sec.background}

$I$ will denote the unit interval $[0,1]$. $\gamma$ denotes the chosen base interval in $\EmbIM$. When $M=B^4$ and $M=S^1\times B^3$, $\gamma$ will be along the $x-axis$ of $B^4$ and a $B^3$ slice respectively. 

Let $\gamma_1 \in \EmbIM$. As described in \cite{budney} the \emph{domain support} of $\gamma_1$ is the closure of the subset of the embedded $I$ on which $\gamma_1$ does not agree with $\gamma$. The \emph{support range} of $\gamma_1$ is the image of the domain support of $\gamma_1$. We say that two embeddings $\gamma_1$ and $\gamma_2$ have disjoint supports if they have disjoint domain supports and disjoint range supports.

\begin{defn}\label{defn: concatenation}
    Let $\gamma_1, \gamma_2 \in \EmbIM$ have disjoint supports. We use $\gamma_1 \prl \gamma_2$ to be the embedding agreeing with $\gamma_1$ and $\gamma_2$ on their respective supports and the base loop everywhere else. This operation extends to maps $X\to \EmbIM$.
\end{defn}

\subsection{Loops in embedding spaces}\label{subsec: lassos}

\begin{figure}
    \centering
    \begin{subfigure}{\textwidth}
        \includegraphics[width=\textwidth, trim = {0 24cm 0 1cm}, clip]{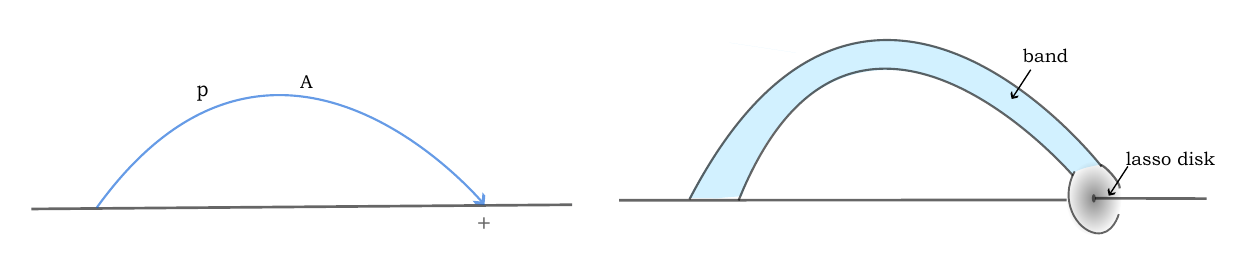}
        \caption{Positive Lasso $L_A$}
    \label{fig: pos lasso}
    \end{subfigure}
    
    \begin{subfigure}{\textwidth}
        \includegraphics[width=\textwidth, trim = {0 12cm 0 13cm}, clip]{lassoLA.pdf}
        \caption{Arc traversing `past' disk (left) and `future' disk (right)}
    \label{fig: pastfuturedisk}
    \end{subfigure}
    
    \begin{subfigure}{\textwidth}
        \includegraphics[width=\textwidth, trim = {0 18cm 0 6.5cm}, clip]{lassoLA.pdf}
        \caption{Negative Lasso $L_A$}
    \label{fig: neg lasso}
    \end{subfigure}
    \caption{Lasso $L_A$}
\end{figure}

We depict loops in embedding space (which we call lassos) via chord diagrams where all chords have the same color. Chords are labeled with an uppercase letter (like $A$) and decorated with an element of $\pi_1 M$ ($p$ in left figure in \cref{fig: pos lasso}). $p$ is the homotopy class described by the loop based at the base of the chord $A$, travels along $A$ until a specific point on $I$ (decorated either by $+$ or $-$) and then returns along $I$ until the base of the chord. (If $n$ colors of chords are in a chord diagram it will be used to depict a map $I^n \to \EmbIM$ like in \cref{subsec: gpq}.)

\begin{defn}
We denote a lasso given by a chord $A$ by $L_A\colon I \to \EmbIM$.

We define a lasso around a loop $p\in \pi_1(M)$ (as described in Figure 14 from \cite{budney} and \cref{fig: pos lasso}) by concatenating the following stages.
\begin{enumerate}
    \item The arc traverses upwards along a band in a neighborhood of the chord A.
    \item The lasso sphere normal to the lasso point (the end of the chord A) can be split into two hemispherical disks. The first disk is traversed in the past and the second is traversed in the future. 
    \item These two disks intersect in a boundary circles that lies in the present and is the unit normal bundle in the present of the arc containing the lasso point at that point . 
    \item We call the disk normal bundle at the same point in the present the "lasso disk". Hence the past and future hemispherical disks project to the lasso disk in the present.
    \item The lasso arc traverses the "past disk" and then the "future disk" and at this stage is at the end of the band closest to the lasso point.
    \item The arc then returns to the base along the band.
\end{enumerate}

$L_A$ is shown in \cref{fig: pos lasso}. In \cref{fig: pastfuturedisk}, any arcs in green are in the present. In this figure, the arc starts at the top of the `past' lasso disc/hemsiphere (shown in red), and as the arc traverses the past disk it gradually changes from red to green. The arc returns along the `future' disk/hemisphere (this is shown as the sequence of arc changing from green to purple). 

\end{defn}

\begin{defn}\label{def: pos lasso}
    A positive lasso has the right boundary of the band pass ``over" the arc. 
\end{defn}

\begin{rem}
    The positive lasso is defined identically to \cite{budney}, and thus by Lemma 4.4 of \cite{budney} a ``negative lasso" (an inverse in $\pi_1\EmbIM$ to the corresponding positive lasso) has the right boundary of the band go under the arc (see \cref{fig: neg lasso})
\end{rem}

\begin{defn}
Let $A_1, \cdots A_n$ denote non intersecting chords, we write $L_{A_1 \cdots A_n}$ to mean the loop of embeddings $L_{A_1}\prl L_{A_2} \cdots \prl L_{A_n}$. 

The chord diagram for $L_{A_1B_1}$ is shown in \cref{fig: undo chords} on the left. 
\end{defn}

\begin{defn} \label{defn: undo}
If $A_1$ and $B_1$ are parallel chords of opposite sign (see \cref{fig: undo chords}), $L_{A_1B_1}$ is null homotopic in $\EmbIM$ via the \emph{undo} null homotopy $U_{A_1B_1}\colon I^2 \to \EmbIM$ defined by the following stages (as described in \cite[Figure 63]{budney} and shown in \cref{fig: undo stages}).

\begin{figure}

    \includegraphics[width=\textwidth, trim = {0 23cm 0 2cm}, clip]{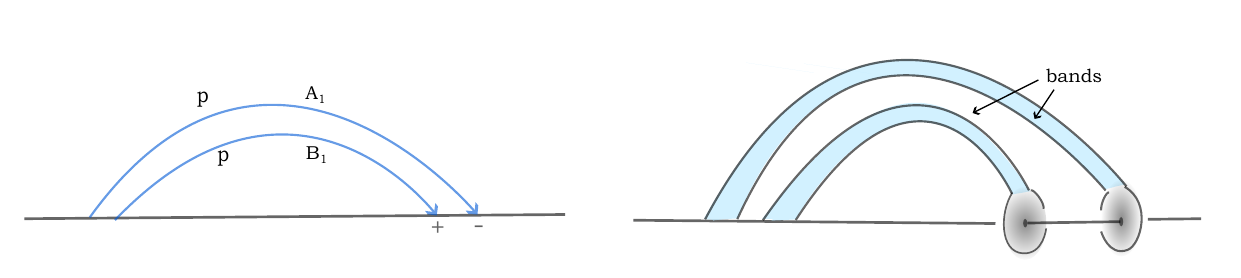}
    \caption{Undo Homotopy: Chord diagram}
    \label{fig: undo chords}
\end{figure}

\begin{figure}
    \centering
    
\begin{subfigure}{\textwidth}
    \includegraphics[width=\textwidth, trim = {0 16cm 0 9cm}, clip]{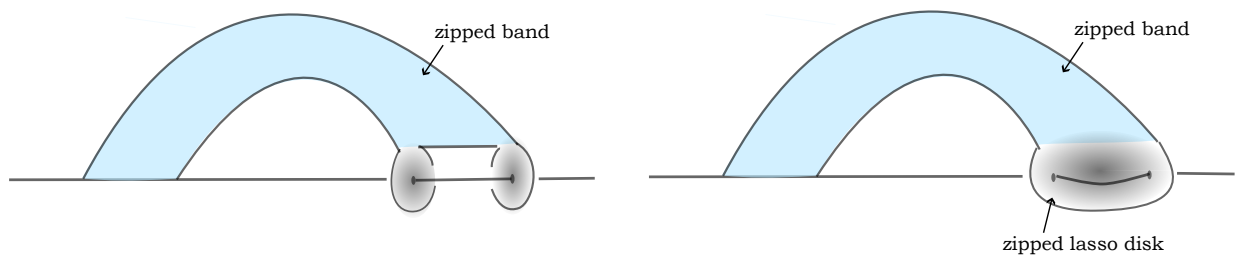}
    \caption{Undo Homotopy: Zipping bands and lassos}
    \label{fig: undo zip}
\end{subfigure}
\hfill
\begin{subfigure}{\textwidth}
    \includegraphics[width=\textwidth, trim = {0 10.5cm 0 14cm}, clip]{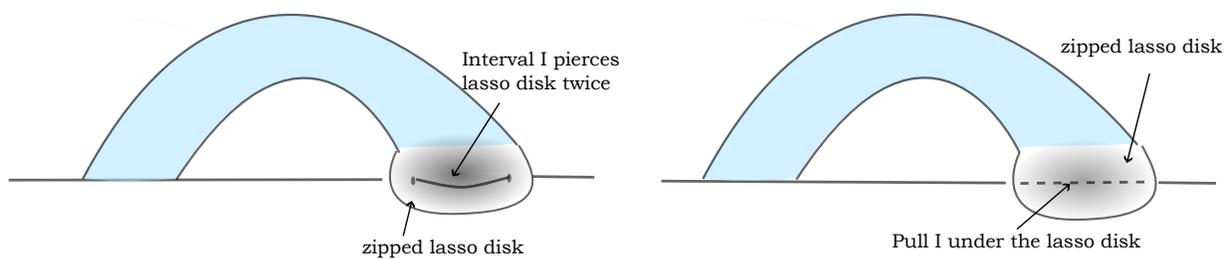}
    \caption{Undo Homotopy: Pulling interval out of the lasso disk}
    \label{fig: undo pull}
\end{subfigure}

\begin{subfigure}{\textwidth}
    \centering
    \includegraphics[width=0.5\textwidth, trim = {0 4cm 8cm 21cm}, clip]{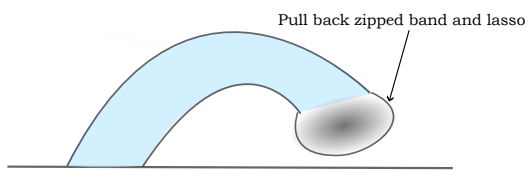}
    \caption{Undo Homotopy: Retracting the zipped chords and zipped lasso disks}
    \label{fig: undo retract}
\end{subfigure}

\caption{Stages of the undo homotopy}
\label{fig: undo stages}
\end{figure}

\begin{enumerate}
    \item Zip up the band to one whose base joins the leftmost point of the left band's base and the rightmost point of the right band's base (\cref{fig: undo zip}). At this stage, the loop in $\EmbIM$ has the arc travel up the zipped band, and then two portions of that arc travel down the \textit{past} disk and then travel back to the zipped band along the \textit{future} disk before returning back along the zipped band. 
    \item Zip the lasso disk to one that contains both lasso disks (\cref{fig: undo zip}). This can be done because the positive and negative lasso disks can be zipped together without passing through the lasso'd portion of $I$. In this stage, one can see the arc being lasso'd around starts behind the zipped lasso disk, pokes out of the lasso disk and pokes it again to leave. (left of \cref{fig: undo pull}) 
    \item Pull out the arc from the lasso disk. This can be done because the lasso sphere exists in either the past or the future except for the boundary of the lasso disk, which exists in the present. 
    (right of \cref{fig: undo pull})
    \item Now that the lasso disk doesn't intersect $I$, we can retract the zipped lasso disks and band back to the base of the chords. (\cref{fig: undo retract}) 
\end{enumerate}    
\end{defn}


\begin{defn} \label{defn: backtrack}
    The \emph{backtrack} null homotopy $B_{A_1 \cdots A_n}$ of $L_{A_1 \cdots A_n}$ in $\ImmIM$ is given by gradually retracting the lasso bands and disks back to the base of the lasso. If we are doing the backtrack homotopy on all the chords (or if it is clear from context which chords get the backtrack homotopy) we may simply denote it as $B$.
\end{defn}

\begin{defn} \label{defn: full}
    Given the map $L_{ABCD} \colon I \to \EmbIM$ for positive (or negative) chords $A, C$ and negative (or positive respectively) chords $B, D$ nested in the order $A, B, C, D$ from innermost to outermost as shown in \cref{fig: full chords}, the \emph{full} null homotopy $F_{ABCD}$ of $L_{ABCD}$ is given by $U_{BC}$ followed by $U_{AD}$. When the chord labels are clear from context, we may simply denote this by $F$.
\begin{figure}
    \centering
    \includegraphics[width=7cm, trim = {1cm 21cm 12cm 5cm}, clip]{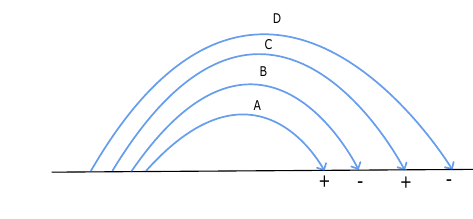}
    \caption{Chords for full null homotopy}
    \label{fig: full chords}
\end{figure}

\end{defn}

\begin{rem}
    In \cref{subsec: cubical} we define notation that makes the full null homotopy $F= (U_{BC} \; \prl \; \id_2 L_{AD}) \star_2 U_{AD}$.
\end{rem}

\subsection{Homotopy limits and stratified spaces}\label{subsec: homlim}

The limit of a diagram in the category of sets or spaces can be defined as the set whose elements consist of a point in each space of the diagram which is equal to the image of each other such point under the maps in the diagram. The homotopy limit of a diagram of spaces relaxes the requirement of equality to merely paths with coherence homotopies, so that for a sequence of spaces $A_0 \to \cdots A_k$ in the diagram with respective elements $a_0,...,a_k$, the data of an element of the homotopy limit includes a $k$-simplex in $A_k$ between the $k+1$ images of those points in $A_k$. Among many equivalent definitions of homotopy limits, the following succinctly packages the data described above of an element of the homotopy limit of a diagram. 

\begin{defn}[{\cite[Definition 1.2]{sinha}}]
    The homotopy limit of a diagram $F \colon \C \to \Top$ is the space of natural transformations $|\C/-| \to F$, where for an object $c$ in $\C$ the space $|\C/c|$ is the geometric realization of the nerve of the category of morphisms into $c$ and commuting triangles between them.
\end{defn}

We will sometimes say ``a homotopy limit'' for any space with a cone over $F$ which is weakly equivalent to ``the'' homotopy limit of $F$ as defined above (just as a limit of a diagram is defined only up to isomorphism, a homotopy limit is defined only up to weak equivalence). A motivating property of homotopy limits is that they preserve weak equivalences, in the sense that the homotopy limits of two naturally weakly equivalent diagrams of spaces will themselves be weakly equivalent.

When the category $\C$ is a poset with a terminal object $e$ and the diagram $F$ consists only of inclusion maps, as is the case in our diagram, its homotopy limit admits a simplified description. The images of those inclusions can be considered as ``strata,'' or subspaces of $F(e)$ which are nested according to the morphisms in $\C$. There are several competing definitions of stratified spaces in the literature, but we define them here in the simplest possible way for how they are used in the relevant homotopy limits.

\begin{defn}
    For a poset $P$, a $P$-stratification of a space $X$ is a functor from $P$ to subsets of $X$, and a $P$-stratified map between such $P$-stratified spaces is a map $f \colon X \to Y$ such that for each $p \in P$ and $x$ in the subset (``stratum'') $X_p$ corresponding to $p$, $f(x) \in Y_p$.
\end{defn}

In particular, for each object $c$ of $\C$, the stratum $F(e)_c$ is the image of $F(c)$ in $F(e)$. In this setting, an element of the homotopy limit of $F$ can be reduced to the data coming from the space $F(e)$, which is described using the stratification structure. This data is based on the $\C$-stratified space $|\C|$, the geometric realization of the nerve of $\C$, with strata given by the images of the inclusions $|\C/c| \to |\C/e| \cong |\C|$.

\begin{prop}[{\cite[Proposition 1.3]{sinha}}]\label{stratapreserving}
    The homotopy limit of a diagram $F \colon \C \to \Top$, where $\C$ is a poset with a terminal object $e$ and the maps in $F$ are all suitably nice inclusions, is given by the space of stratified maps from $|\C| \to F(e)$.
\end{prop}

\subsection{Configuration spaces}\label{subsec: confspace}

\begin{defn}

We denote the set $\{ 1, 2, \cdots k\}$ as $[k]$
\end{defn}

\begin{defn}
    The $k$ point configuration space of a manifold $M$ is denoted by $C_k(M)$ and is defined as $$C_k(M) \coloneq \sb{(p_1, p_2, \cdots p_k) \in M^k}{p_i \neq p_j \text{ when } i \neq j}$$
\end{defn}

We will need a variant of configuration spaces with some extra data -  unit tangent vectors associated to each point - that we define below.

\begin{defn}
    $$C'_{k}(M) \coloneq \sb{\left((p_i, v_i)\right)_{i\in [k]} \in (STM)^k}{p_i\neq p_j \text{ when } i\neq j}$$
\end{defn}

We define a compactification of $C_k(M)$ as in \cite[Definition 4.1]{sinha}
\begin{defn}
    Suppose $f\colon M \to \mathbb{R}^n$ is an embedding of $M$ into Euclidean space, and $\S\subseteq[k]$ we define the following.
    
    \begin{enumerate}
        \item Maps $\pi_{i,j}\colon C_k(M) \to S^{n-1}$ given by $\frac{f(p_i) - f(p_j)}{\Vert f(p_i)- f(p_j)\Vert}$. 
        \item $C_2(k)\coloneq \sb{(i,j)}{1\leq i<j\leq k}$
        \item $C_2(\S)\coloneq \sb{(i,j)}{i< j, i,j \in \S}$
        \item $A_k\lr{M}\coloneq M^k \times \left(S^{(n-1)}\right)^{C_2(k)}$
    \end{enumerate}

\end{defn}

\begin{defn}
The compactified $k$ point configuration space of $M$ is denoted by $C_k\lr{M}$ and is the closure of $C_k(M)$ in $A_k\lr{M}$ via the map $\left( \iota , (\pi_{i,j})_{(i,j) \in C_2(k)}\right)$.
\\\\We can similarly define $C'_k\lr{M}$ as the closure in $(STM)^k\times\left(S^{(n-1)}\right)^{C_2(k)}$.

Points in $C_k\lr{M}$ consist of tuples $(p_1, \cdots, p_k)$ with pairwise disjoint points along with boundary points where we could have $p_i=p_j$, in which case we add the data of a unit tangent vector $v_{ij}$ in M for every pair of colliding points $p_i, p_j$ that specifies the \emph{direction} that those two points collide in. 

\end{defn}

These `colliding' faces along with $C_k\lr{M}$ make $C_k\lr{M}$ a stratified space. We can describe strata $C^{\S}_k\lr{M}$ for each subset $\S \in [k]$. 
$$C^S_k\lr{M} \coloneq \sb{\left((p_i)_{i\in[k]}, (v_{ij})_{(i,j)\in C_2(\S)}\right)}{p_i = p_j \text{ when } i, j\in \S}$$

When $\S=\phi$, $C^{\S}_k\lr{M} = C_k\lr{M}$. We define $C'^{\;\S}_k\lr{M}$ analogously. When $\S_1 \subseteq \S_2$, we have $C'^{\;\S_2}_k\lr{M} \hookrightarrow C'^{\;\S_1}_k\lr{M}$. 
This allows us to define maps $\partial^i\colon C'_{k-1}\lr{M} \to C'_{k}\lr{M}$ for $1\leq i \leq k-1$ which shifts up by 1 the indices of all points $p_j$ for $j>i$, sets $p_{i+1}=p_i$ and sets $v_{i\; {i+1}} = v_{i+1} = v_i$. 

We will use a more specialized subspace called the aligned stratum.
$$C^{\textrm{align}}_k\lr{M} \subset C'_k\lr{M} \text{ such that } v_{ij} = v_i = v_j \text{ when } p_i = p_j$$

We see that the connected component where $0=p_0\leq p_1\leq \cdots p_k \leq p_{k+1}=1$ of $C^{align}_{\partial,k}\lr{I}$ is homeomorphic as a stratified space to the standard $k$-simplex $\Delta^k$.

We also describe here some special elements of $\pi_{\dim M -1} C_k\lr{M}$. We will define them here for dim$(M) = 4$, but they generalize accordingly.

\begin{defn}\label{defn: wij}
    The element $w_{ij} \in \pi_3C_k\lr{M}$ is the point $p_i$ traversing the sphere normal bundle of $p_j$ in M. 
\end{defn}

\begin{defn}\label{defn: twij}
    Suppose $\alpha \in \pi_1(M)$, we define $\twij{i}{\alpha}{i}{j} \in  \pi_3C_k\lr{M}$ as point $i$ traversing the loop $\alpha$ before traversing the sphere normal bundle of $p_j$ in $M$
\end{defn}

\begin{rem}
When $M = S^1\times B^3$, $\pi_1(M) \cong \mathbb{Z}$, so for $p\in \mathbb{Z}$, we will write $\twij{i}{p}{i}{j}$ to mean the element obtained by point $i$ circling the $S^1$ direction $p$ times before traversing the sphere normal bundle of $p_2$ in $M$. 
\end{rem}

To suit the spaces that we will use to approximate $\EmbIM$ in \cref{subsec: embeddingcalc}, we define a variant of configuration spaces where the first and last point are fixed on $\partial M$

\begin{defn}\label{defn: confbound}
    For $k\geq 0$, we define
    $$C_{\partial, k}(M) \coloneq \sb{(p_0, p_1 \cdots p_{k+1}) \in M^{k+2}}{p_0 = \ast_0, p_{k+1} = \ast_1, p_i\neq p_j \text{ when } i\neq j}$$
\end{defn}

We can anaogously define $C'_{\partial, k} \lr{M}$ and $C^{align}_{\partial, k}\lr{M}$.

\subsection{Embedding calculus}\label{subsec: embeddingcalc}


Functors such as $\Embp(-,M)$ to spaces from the opposite category of open subsets of $I$ containing the endpoints, which have relatively few convenient properties beyond preserving weak equivalences, are often studied using a sequence of increasingly accurate approximations in analogy with the Taylor approximation of a smooth function. These approximations come equipped with connectivity results that show the homotopy groups of embedding spaces such as $\EmbIM$ in sufficiently low dimensions to agree with those of its approximations. We give here the basic definitions of this ``embedding calculus'' and describe how it is used to simplify the study of $\pi_n\EmbIM$.

\begin{defn}
    For $\square_k$ the poset $(0 < 1)^k$ and $P_k$ the ``punctured $n$-cube category'' given by the poset $\square_k \sans (1,...,1)$, a diagram $D \colon \square_k \to \C$ is cocartesian if $D(1,...,1)$ is a colimit of the restricted diagram $P_k \to \square_k \to \C$. 
    
    A diagram $D \colon \square_k^{op} \to \Top$ is homotopy cartesian if $D(1,...,1)$ is a homotopy limit of the restricted diagram $P_k^{op} \to \square_k^{op} \to \Top$.
    
    A functor $F \colon \C^{op} \to \Top$ to spaces is $k$-polynomial if for every cocartesian diagram $D \colon \square_{k+1} \to \C$, the composite diagram $\square_{k+1}^{op} \To{D} \C^{op} \To{F} \Top$ is homotopy cartesian.
\end{defn}

Goodwillie, Klein, and Weiss showed in \cite{goodwillie} that for any functor $F \colon \C^{op} \to \Top$ which preserves weak equivalences, where $\C$ is a poset of open subsets of some space, there is a $k$-polynomial functor $T_k F \colon \C^{op} \to \Top$ with a natural transformation $F \to T_k F$. There are also natural fibrations $T_k F \to T_{k-1} F$ commuting under $F$.

We are particularly interested in the functor $\Embp(-,M)$ for $M$ a manifold and $\C$ the poset of open subsets of $I$ containing the endpoints. In this case, the maps $\Embp(X,M) \to T_k\Embp(X,M)$ are $(k-1)(\dim M - 3)$-connected, and $\Embp(X,M)$ is the homotopy limit of the sequence $$\cdots \to T_1\Embp(X,M) \to T_0\Embp(X,M)$$

We use the same model for $T_k\EmbIM$ as used in \cite{sinha}. When $I = I'_0 \cup I_1 \cup I'_1 \cup I_2 \cdots I'_{k+1}$, a concatenation of intervals, $T_k\EmbIM$ is given by the homotopy limit of the punctured cubical diagram that sends a subset $S\subset\{ 1, \cdots k+1\}$ to $\Embp\left(I\sans \left(\cup_{i\in S} I_{i}\right), M\right)$.  

We work out the example for $T_2$. Let $I = I'_0 \cup I_1 \cup I'_1 \cup I_2 \cup I'_2 \cup I_3 \cup I'_3$. As the functor $T_2\Embp(-,M)$ is 2-polynomial, the space $T_2\EmbIM$ will be the homotopy limit of the diagram $P_3^{op} \to \Top$ pictured below.

\[
\begin{tikzcd}[column sep=-40,row sep=large]
    & & \Embp(I \sans I_2,M) \ar{dl} \ar{dr} \\[10pt]
    & \Embp(I \sans (I_1 \cup I_2),M) \ar{dr} & & \Embp(I \sans (I_2 \cup I_3),M) \ar{dl} \\[-10pt]
    & & \Embp(I \sans (I_1 \cup I_2 \cup I_3),M) \\
    \Embp(I \sans I_1,M) \ar{uur} \ar{rr} & & \Embp(I \sans (I_1 \cup I_3),M) \uar & & \Embp(I \sans I_3,M) \ar{ll} \ar{uul}
\end{tikzcd}
\]

We know that that $\Emb(I, M) \simeq STM$ (note here that we don't require fixed endpoints). Suppose $I = I'_0 \cup I_1 \cup I'_1 \cdots I'_{k+1}$. A similar argument shows that $\Emb(I_1 \cup I_2 \cup \cdots \cup I_k, M) \simeq C'_{k}(M)$. This, along with the homotopy invariance and condition that endpoints of embeddings in $\Embp(U,M)$ are fixed, shows that 
$$\Embp\left(I \setminus {\left(\cup_{i\in [k+1]} I_i\right)} , M\right) \simeq \text{Emb}_\partial\left(I'_0 \cup I'_1 \cup \cdots \cup I'_{k+1} , M\right) \simeq C'_{\partial, k}(M) \simeq C'_{\partial, k}\lr{M}.$$ 
This allows us to replace our punctured cubical diagram above with the following while preserving the homotopy type of its homotopy limit. 
\[
\begin{tikzcd}[row sep=large]
    & & C'_{\partial,0}\lr{M} \ar{dl} \ar{dr} \\[5pt]
    & C'_{\partial,1}\lr{M} \ar{dr} & & C'_{\partial,1}\lr{M} \ar{dl} \\[-20pt]
    & & C'_{\partial,2}\lr{M} \\[-10pt]
    C'_{\partial,0}\lr{M} \ar{uur} \ar{rr} & & C'_{\partial,1}\lr{M} \uar & & C'_{\partial,0}\lr{M} \ar{ll} \ar{uul}
\end{tikzcd}
\]

In this diagram of suitably nice inclusions, the space $C'_{\partial,2}\lr{M}$ has strata given by the images of the three copies of $C'_{\partial,1}\lr{M}$ and their pairwise intersections which are the images of $C'_{\partial,0}\lr{M}$. 

The stratified space $\left|P_3^{op}\right|$ is precisely the 2-simplex $\Delta^2$, based on the shape of the diagram above, and its strata are given by the edges and vertices of $\Delta^2$. This agrees with the aligned component of $C'_{\partial,2}\lr{I}$.

The homotopy limit of this diagram then (\cref{stratapreserving}) is the space of strata preserving maps from the 2-simplex to $C'_{\partial,2}\lr{M}$. In particular, this means that the three edges are sent to the strata arising from $C'_{\partial,1}\lr{M}$ and the vertices are sent to the strata arising from $C'_{\partial,0}\lr{M}$. 

In a similar manner, $T_k\EmbIM$ can be shown to be the space of strata preserving maps from the $k$-simplex to $C'_{\partial,k}\lr{M}$. This leads to a theorem nearly identical to \cite[Theorem 5.4]{sinha}.

\begin{thm}
Let $Map^{sp}(C'_{\partial, k}\lr{I}, C'_{\partial, k}\lr{M})$ denote the space of strata preserving maps that send the aligned stratum of $C'_{\partial, k}\lr{I}$ to the aligned stratum of $C'_{\partial, k}\lr{M}$, then 
$$T_k\EmbIM \simeq Map^{sp}(C'_{\partial, k}\lr{I}, C'_{\partial, k}\lr{M})$$
and by the connectivity result this space agrees with $\EmbIM$ on $\pi_i$ for $i=0,...,(k-1)(\dim M - 3)$. 
\end{thm}

So for $M$ a 4-manifold where we are interested in $\pi_3\EmbIM$, it suffices to consider $T_4\EmbIM$ which is the space of strata-preserving aligned maps $C'_4\lr{I} \to C'_4\lr{M}$ which the above theorem shows will agree on $\pi_3$. We will often drop the $'$ in $C'_k\lr{M}$ and restrict our attention to strata preserving maps from $C_k\lr{I}$ to $\C_k\lr{M}$.
\begin{defn}
Given a map $F\colon X\to \EmbIM$, we define $T_k F\colon X \times C_k\lr{I} \to C_k\lr{M}$ to be the induced map on compactified configuration spaces.    
\end{defn}
Sometimes we will use the same notation $T_k F$ when studying the induced map $X \to Map(C_k\lr{I}, C_k\lr{M})$.

\subsection{Cosimplicial model for $\EmbIM$}\label{subsec: cosimplicial}

In \cref{subsec: embeddingcalc}, we discussed how $\EmbIM$ is the homotopy limit of the tower of fibrations $T_0\EmbIM\leftarrow T_1\EmbIM \leftarrow \cdots $ where each level is given by a homotopy limit of a punctured cubical diagram of configuration spaces. Sinha \cite[Theorem 7.1]{sinha} shows that this is equivalent to $\EmbIM$ being the totalization of the cosimplicial space that sends $[n] \to C'_{\partial, n}\lr{M}$. The $i^{th}$ codegeneracy map $s_i\colon C'_{\partial, n}\lr{M} \to C'_{\partial, n-1}\lr{M}$ is the map that `drops' the $i^{th}$ point for $1\leq i \leq n$. The $i^{th}$ coface map, $\partial^i\colon C'_{\partial, n}\lr{M} \to C'_{\partial, n+1}\lr{M}$ `doubles' the $i^{th}$ point when $0\leq i \leq n+1$. (Note that when $i=0$ or $i=n+1$, the doubled point $p_i$ is one of the fixed endpoints from $\partial M$.)

Sinha then shows that this gives rise to a second quadrant (Bousfield-Kan) spectral sequence such that 
$$E_1^{-p, q} = \bigcap_i \text{ker } s_i \subseteq \pi_q(C'_{\partial, p}\lr{M}) \cong \pi_qC'_{p}(M)$$ for $p, q \geq 0$, where the $d_1\colon E^{-p,q}_1 \to E^{-p-1, q}_1$ differential is the restriction of the map $$\sum\limits_i (-1)^i \partial^i\colon \pi_q(C'_p\lr{M}) \to \pi_q(C'_{p+1}\lr{M}).$$
In general, the $d_r$ differential goes from $E^{-p,q}_r$ to $E^{-p-r, q+r-1}_r$.

In \cref{subsec: relationsspectral} we compute some $d_1$ differentials when $M=S^1\times B^3$.

\subsection{Whitehead products}\label{subsec: whitehead}

\begin{defn}
Given maps $f\colon (D^n, \partial D^n) \to (X, x_0)$ and $g\colon (D^m, \partial D^m) \to (X,x_0)$ we can define its \textit{Whitehead product} as a map $$[f, g]\colon (D^{n+m-1}, \partial D^{n+m-1}) \to (X, x_0)$$ as follows. 
\begin{itemize}
    \item Inside $D^{n+m-1}$, we can find a generalized Hopf link of $ S^{m-1}$ and $ S^{n-1}$. The disk normal bundles of these $S^{m-1}$ and $S^{n-1}$ are $N_a: S^{m-1}\times D^n$ and $N_b: S^{n-1}\times D^m$ respectively. 
    \item $[f,g]$ maps $D^{n+m-1} \sans (N_a \cup N_b)$ to the base point $x_0$. 
    \item $[f,g]$ maps points $(p_a, q_a) \in S^{m-1} \times D^n = N_a$ to $f(q_a)$ and maps points $(p_b, q_b) \in  S^{n-1} \times D^m = N_b$ to $g(q_b)$. 
\end{itemize}

This induces a well defined map on the product of homotopy classes $$[\cdot,\cdot]: \pi_n(X,x_0) \times \pi_m(X, x_0) \to \pi_{n+m-1}(X, x_0)$$ 

The Whitehead product is bilinear, graded symmetric ($[f, g] = (-1)^{kl}[g, f]$) and satisfies a Jacobi relation: 

$$(-1)^{km}[[ f , g], h] + (-1)^{lm}[[g, h], f ] + (-1)^{mk}[[h, f ], g] = 0$$
where $f \in \pi_{k}X, g \in \pi_{l}X, h \in \pi_{m}X$ and $k,l,m \geq 2$.
\end{defn}

Milnor and Moore \cite{milnor-moore} first described the rational homotopy groups $\Q\otimes \pi_\ast C_k(B^n)$ as generated by the classes $w_{ij}$ defined in \cref{subsec: confspace} subject to the following relations: 

\begin{itemize}
    \item $w_{ii}=0$
    \item $w_{ij} = (-1)^{n}w_{ji} $
    \item $[w_{ij}, w_{kl}] = 0$ if $\{i,j\} \cap \{k,l\} = \phi$
    \item $[w_{ij}, w_{jk}] = [w_{jk}, w_{ki}] = [w_{ki}, w_{ij}]$
\end{itemize}

Budney and Gabai extend this approach to describe rational homotopy groups of $S^1\times B^n$ which are generated by $\twij{i}{p}{i}{j}$ subject to some additional relations that we describe in \cref{subsec: proppinck}.

\subsection{Algebraic constructions of cubical maps}\label{subsec: cubical}

The constructions and results of \cref{sec.gpqr,sec.nullhomotopy} involve increasingly complicated maps from cubes $I^n$ into various spaces associated to embeddings. To more easily describe these maps and how they are build up in piecewise fashion, we use operations for concatenated, constant, folded, rotated, and reversed maps from cubes inspired by the theory of cubical $\omega$-groupoids (see \cite{browncubes}), as suggested by Brandon Shapiro. While we define from scratch all of these operations, the idea is that a space $X$ has an associated algebraic structure consisting of all maps $I^n \to X$ regarded as ``$n$-dimensional cubical cells'' which are equipped with operations including composition (concatenation), units (constant maps, folds), and inverses (reverse maps) in all $n$ directions, and that this ``cubical $\omega$-groupoid'' contains all of the homotopical information of $X$.

\begin{defn}
    Given $f \colon I^n \to X$, we define $\face^\varepsilon_{i}\colon I^{n-1} \to X$ for $\varepsilon = 0,1$ and $i = 1, \dots n$ as the map 
    $$I^{n-1} \cong I^{i-1} \times I^{n-i} \To{\varepsilon} I^{i-1} \times I \times I^{n-i} \cong I^n \To{f} X$$ 
    This denotes the restriction of $f$ to the front face in the $i^{th}$ direction when $\varepsilon = 0$ (and back face when $\varepsilon = 1$).
\end{defn}

If $h \colon I^n \to X$ with $\face_i^0 h = f$ and $\face_i^1 h = g$, we will sometimes write $f \Tto{h}{i} g$. If furthermore $\face_j^0 h = k$ and $\face_j^1 h = \ell$, we will often depict $h$ as below. This style of picture will also occasionally be used with three dimensions displayed.
\[
\begin{tikzcd}
    \cdot \rar{g} & \cdot \\ \cdot \uar{k} \rar[swap]{f} \ar[phantom]{ur}{h} & \cdot \uar[swap]{\ell}
\end{tikzcd}\qquad\compass{i}{j}
\]

\begin{defn}
    Given $f \colon I^n \to X$, we define the map $\id_i \colon I^{n+1} \to X$ as a projection along the $i^{th}$ coordinate followed by $f$. So, $$\id_i \colon I^{n+1} \cong I^{i-1} \times I \times I^{n+1-i} \to I^{i-1} \times I^{n+1-i} \cong I^n \To{f} X$$
\end{defn}

$\id_i f$ will be depicted as $f \begin{tikzcd}[sep=13pt, cramped, ampersand replacement=\&, text height=1ex, text depth=.3ex]\ar[equals,r,"i"']\&{}\end{tikzcd} f$, and for $f \Tto{h}{i} g$, $\id_j h$ will be depicted as below left for any $j > i$ and as below right for any $j < i$.
\[
\begin{tikzcd}
    f \rar{h} & g \\ f \uar[equals]{\id_j f} \rar[swap]{h} \ar[phantom]{ur}{\id_j h} & g \uar[equals,swap]{\id_j g}
\end{tikzcd}\qquad\compass{i}{j}
\qquad\qquad\qquad
\begin{tikzcd}
    f \rar{h} & g \\ f \uar[equals]{\id_j f} \rar[swap]{h} \ar[phantom]{ur}{\id_j h} & g \uar[equals,swap]{\id_j g}
\end{tikzcd}\qquad\compass{i+1}{j}
\]

We will often consider maps which are constant in not just one but multiple directions.

\begin{defn}
    Given $f \colon I^n \to X$ and $0 < i_1 < \cdots < i_k < n+k$, we write $\id_{i_1,...,i_k}$ for the map $\id_{i_k}\cdots\id_{i_1} f \colon I^{n+k} \to X$ which is constant in the $i_1,...,i_k$ directions. For $f \colon I^n \to X$, when $k$ is clear from context we write $\id f \colon I^{n+k} \to X$ to denote $\id_{n+1,...,n+k}f$ (which is constant in all directions above $n$). 
\end{defn}

\begin{defn}
    Given $f \colon I^n \to X$, we define the map $\rev_i\colon I^{n} \to X$ as the map that reverses $f$ along the $i^{th}$ coordinate. So $\rev_i$ is the map 
    $$I^n \cong I^{i-1} \times I \times I^{n-i} \To{a \; \mapsto \; 1-a} I^{i-1} \times I \times I^{n-i} \cong I^n \To{f} X$$
\end{defn}

For $f \Tto{h}{i} g$, $\rev_i h$ has the form $g \Tto[35pt]{\rev_i h}{i} f$, and for $h$ as below left, $\rev_i h$ has the form below right

\[
\begin{tikzcd}
    c \rar{g} & d \\ a \uar{k} \rar[swap]{f} \ar[phantom]{ur}{h} & b \uar[swap]{\ell}
\end{tikzcd}\qquad\compass{i}{j}
\qquad\qquad\qquad
\begin{tikzcd}
    d \rar{\rev_i g} & c \\ b \uar{\ell} \rar[swap]{\rev_i f} \ar[phantom]{ur}{\rev_i h} & a \uar[swap]{k}
\end{tikzcd}\qquad\compass{i}{j}
\]

\begin{defn}
    Given $f \colon I^n \to X$, we define the map $\fold^{0,0}_{i,j}\colon I^{n+1} \to X$ as the map 
    $$I^{n+1} \cong I^{i-1} \times I \times I^{j-i-1} \times I \times I^{n+1-j} \To{(a,b) \; \mapsto \; 1 - (1-a)(1-b)} I^{i-1} \times I \times I^{n-i} \cong I^n \To{f} X$$
    Given $f \colon I^n \to X$, we define the map $\fold^{1,1}_{i,j}\colon I^{n+1} \to X$ as the map 
    $$I^{n+1} \cong I^{i-1} \times I \times I^{j-i-1} \times I \times I^{n+1-j} \To{(a,b) \; \mapsto \; ab} I^{i-1} \times I \times I^{n-i} \cong I^n \To{f} X$$
\end{defn}

For $f \Tto{h}{i} g$, $\fold^{0,0}_{i,j} h$ has the form below left and $\fold^{1,1}_{i,j}$ has the form below right.

\[
\begin{tikzcd}[column sep=large]
    g \rar[equals] & g \\ f \uar{h} \rar[swap]{h} \ar[phantom]{ur}{\fold^{0,0}_{i,j} h} & g \uar[equals]
\end{tikzcd}\qquad\qquad\qquad\compass{i}{j}
\qquad\qquad\qquad
\begin{tikzcd}[column sep=large]
    f \rar{h} & g \\ f \uar[equals] \rar[equals] \ar[phantom]{ur}{\fold^{1,1}_{i,j} h} & f \uar[swap]{h}
\end{tikzcd}
\]
We will also frequently use the additional folded maps
\[
\fold^{1,0}_{i,j} h \coloneq \rev_j\fold^{1,1}_{i,j}h \qqand
\fold^{0,1}_{i,j} h \coloneq \rev_j\fold^{0,0}_{i,j}h,
\]
which respectively have the forms below left and right.
\[
\begin{tikzcd}[column sep=large]
    f \rar[equals] & f \\ f \uar[equals] \rar[swap]{h} \ar[phantom]{ur}{\fold^{1,0}_{i,j} h} & g \uar[swap]{\rev_i h}
\end{tikzcd}\qquad\qquad\qquad\compass{i}{j}
\qquad\qquad\qquad
\begin{tikzcd}[column sep=large]
    f \rar{h} & g \\ g \uar{\rev_i h} \rar[equals] \ar[phantom]{ur}{\fold^{0,1}_{i,j} h} & g \uar[equals]
\end{tikzcd}
\]

\begin{rem}
We can see that $\fold^{1-\varepsilon}_{i,j} f = \rev_i \rev_{j} \fold^\varepsilon_{i,j} \rev_{i} f$, but we define $\fold^{0,0}_{i,j}$ and $\fold^{1,1}_{i,j}$ separately for convenience. 
\end{rem}

\begin{defn}
    Given $f \colon I^n \to X$, we define the map $\rot_{i,j} \colon I^n \to X$ for $1 \leq i < j \leq n$ as the map that interchanges the $i^{th}$ and $j^{th}$ coordinates. So $\rot_{i,j}$ is the map 
    $$I^{n} \cong I^{i-1} \times I \times I^{j-i-1} \times I \times I^{n-j} \To{(a,b) \; \mapsto \; (b,a)} I^{i-1} \times I \times I^{j-i-1} \times I \times I^{n-j} \cong I^n \To{f} X$$
\end{defn}

For $h$ of the form below left, $\rot_{i,j} h$ has the form below right.

\[
\begin{tikzcd}
    c \rar{g} & d \\ a \uar{k} \rar[swap]{f} \ar[phantom]{ur}{h} & b \uar[swap]{\ell}
\end{tikzcd}\qquad\compass{i}{j}
\qquad\qquad\qquad
\begin{tikzcd}
    b \rar{\ell} & d \\ a \uar{f} \rar[swap]{k} \ar[phantom]{ur}{\rot_{i,j} h} & c \uar[swap]{g}
\end{tikzcd}\qquad\compass{i}{j}
\]

\begin{defn}
    Given $f,g \colon I^n \to X$ such that $\face_i^1 f = \face_i^0 g$, we define the map $f \star_i g \colon I^n \to X$ as the concatenation of $f$ and $g$ in the $i^{th}$ direction along their shared face. So $f \star_i g$ is the map 
    $$I^n \cong I^{i-1} \times I \times I^{n-i} \cong I^{i-1} \times (I \cup_\ast I) \times I^{n-i} \cong I^n \cup_{I^{n-1}} I^n \To{f \cup g} X$$
\end{defn}

For $f \Tto{k}{i} g$ and $g \Tto{\ell}{i} h$, we have $f \Tto[25pt]{k \star_i \ell}{i} h$, and for $k,\ell$ as below left, $k \star_i \ell$ has the form below right.

\[
\begin{tikzcd}
    \cdot \rar{c} & \cdot \rar{d} & \cdot \\ 
    \cdot \uar{f} \rar[swap]{a} \ar[phantom]{ur}{k} & \cdot \uar[swap]{g} \rar[swap]{b} \ar[phantom]{ur}{\ell} & \cdot \uar[swap]{h}
\end{tikzcd}\qquad\compass{i}{j}
\qquad\qquad\qquad
\begin{tikzcd}
    \cdot \rar{c \star_i d} & \cdot \\ 
    \cdot \uar{f} \rar[swap]{a \star_i b} \ar[phantom]{ur}{k \star_i \ell} & \cdot \uar[swap]{h}
\end{tikzcd}\qquad\compass{i}{j}
\]

Note that $\star_i$ is associative up to homotopy, and we may sometimes write 
\[
f_1 \star_i \cdots \star_i f_k \colon I^n \cong I^n \cup_{I^{n-1}} \cdots \cup_{I^{n-1}} I^n \To{f_1 \cup \cdots \cup f_k} X
\]
for the $k$-fold concatenation in the $i^{th}$ direction (without addressing associativity homotopies). Associativity also applies (strictly in fact) to concatenations in multiple directions at once, so that we can at once compose grids as below where adjacent squares are presumed to agree on their appropriate faces.

\[
\begin{array}{|c|c|c|}
    \hline \;\;f\;\; &  \;\;g\;\; & \;\;h\;\; \\
    \hline k & \ell & m \\
    \hline
\end{array}
\qquad\compass{i}{j}
\]
\[
(k \star_i \ell \star_i m) \star_j (f \star_i g \star_i h) = (k \star_j f) \star_i (\ell \star_j g) \star_i (m \star_j h)
\]

We will often denote such a bidirectional concatenation simply by the grid of its factors as above, rather than as a convoluted expression of nested $\star_i$'s and $\star_j$'s. This notation, which we call \emph{concatenation diagrams}, also conveniently allows us to depict bidirectional concatenations of higher dimensional cubical maps without over-complicating the figures with extraneous dimensions.

\begin{rem}\label{rem:unital}
The operation $\star_i$ is also unital up to homotopy with respect to $\id_i$. This means that for an $f \Tto{h}{i} g$, where $h \colon I^n \to X$, there are maps $I^{n+1} \to X$ of the form
\[
\id_i f \star_i h \; \Tto{}{i} \;h\; \begin{tikzcd}[sep=20pt, cramped, ampersand replacement=\&, text height=1ex, text depth=.3ex]{}\&\ar[l,"i"]\end{tikzcd} \; h \star_i \id_i g.
\]
These maps are called \emph{unitors}, and generalize the standard homotopies witnessing unitality of constant maps in homotopy groups. 
\end{rem}

Finally, we describe several particular combinations of the above operations that arises repeatedly in our constructions. The first corresponds to ``revolving'' a map $I^{n-1} \to X$ around a suitable map $I^n \to X$. 

\begin{defn}\label{defn:revolve}
    Given $f \colon I^n \to X$ of the form below 
    \[
    \begin{tikzcd}
        g \rar[equals]{\id_i g} & g \\ g \uar[equals]{\id_j g} \rar[equals,swap]{\id_i g} \ar[phantom]{ur}{f} & g \uar[equals,swap]{\id_j g}
    \end{tikzcd}\qquad\compass{i}{j}
    \]
    and $g \Tto{\ell}{k} h$, we define the composite
    \[
    \revolve{f}{\ell}
    \quad \coloneq \quad 
    \begin{array}{|c|c|c|}
    \hline \rev_i \fold_{i,j} \rot_{k,i} \ell & \id_i \rot_{k,j} \ell & \fold^{0,0}_{i,j} \rot_{k,i} \ell \\
    \hline \rev_i \id_j \rot_{k,i} \ell & f & \id_j \rot_{k,i} \ell \\
    \hline \rev_j\rev_i\fold^{1,1}_{i,j} \rot_{k,i} \ell & \rev_j \id_i \rot_{k,j} \ell & \rev_j \fold_{i,j} \rot_{k,i} \ell \\
    \hline
    \end{array}
    \]\[\compass{i}{j}
    \]
    which has boundary as below.
    \[
    \begin{tikzcd}[column sep=large,row sep=large]
        \cdot \rar[equals]{\id_i h} & \cdot \\ \cdot \uar[equals]{\id_j h} \rar[equals,swap]{\id_i h} \ar[phantom]{ur}{\;\;\boxed{f}\;\ell} & \cdot \uar[equals,swap]{\id_j h}
    \end{tikzcd}\qquad\compass{i}{j}
    \]
\end{defn}

More generally, we will often consider composites of grids with reflectional symmetry and use similarly simplified notation to only specify their upper right corner.

\begin{defn}\label{defn:flip}
    Given adjacent cubes of the form below, 
    \[
    \begin{tikzcd}
        \cdot \rar & \cdot \rar & \cdot \\ \cdot \uar{o} \rar{\ell} \ar[phantom]{ur}{h} & \cdot \uar[swap]{o} \rar{n} \ar[phantom]{ur}{k} & \cdot \uar \\ \cdot \uar{m} \rar[swap]{\ell} \ar[phantom]{ur}{f} & \cdot \uar[swap]{m} \rar[swap]{n} \ar[phantom]{ur}{g} & \cdot \uar
    \end{tikzcd}\qquad\compass{i}{j}
    \]
    we define the composite
    \[
    \flip{\begin{array}{|c|c|}
    \hline h & k \\
    \hline f & g \\
    \hline
    \end{array}}
    \quad \coloneq \quad 
    \begin{array}{|c|c|c|}
    \hline \rev_i k & h & k \\
    \hline \rev_i g & f & g \\
    \hline \rev_j\rev_i k & \rev_j h & \rev_j k \\
    \hline
    \end{array}
    \qquad\compass{i}{j}
    \]
\end{defn}

%
%

We will occasionally need ``twisting'' homotopies from a map constant in one parameter to a map constant in a different parameter.

\begin{lem}\label{lem:twist}
    For any map $h \colon I^n \to X$ with $f \Tto{h}{i} g$, there is a map $$\twist_i h \colon I^{n+2} \to X$$ with $\id_{i+1} h \Tto[45pt]{\twist_i h}{n+2} \id_i h$ of the form
    \[
    \begin{tikzcd}[column sep=huge,row sep=huge]
        h \rar{\fold^{0,0}_{i,n+1}h} & \id_n g \\ 
        \id_n f \uar{\fold^{1,1}_{i,n+1}h} \rar[swap]{\fold^{1,0}_{i,n+1}h} \ar[phantom]{ur}{\twist_i h} & \rev_i h \uar[swap]{\fold^{0,1}_{i,n+1}h}
    \end{tikzcd}\qquad\compass{i}{i+1}
    \]
    \[
    \begin{tikzcd}
        & f \ar{rr}{h} \ar{dl}[swap]{h} \ar[equals,from=dd] & & g \ar[equals]{dl} \\ 
        g \ar[equals,crossing over]{rr} & & g \\ 
        & f \ar{rr}[swap,pos=.25]{h} \ar[equals]{dl} & & g \ar[equals]{uu} \ar{dl}{\rev_ih} \\
        f \ar{uu}{h} \ar[equals]{rr} & & f \ar[crossing over]{uu}[swap,pos=.75]{h}
    \end{tikzcd}\qquad\outcompass{i}{i+1}{n+2}
    \]
\end{lem}

Topologically, this map could be defined by regarding $I^{n+2}$ as a cylinder with the round part in the $i,i+1$ directions and rotating as one progresses in the $n+2$ direction, but it can also be described using the ``algebraic'' operations we have defined.

\begin{proof}
     We first consider the concatenation of the pair
     \[
     \id_{i+1}h \Tto[100pt]{\fold^{1,1}_{i+1,n+2}\fold^{0,0}_{i,i+1}h}{n+2} \fold^{0,0}_{i,i+1}h \Tto[120pt]{\fold^{1,0}_{i,n+2}\fold^{0,0}_{i,i+1}h}{n+2} \id_i h,
     \]
     where the two component maps have the form below left and below right respectively.
     \[
    \begin{tikzcd}
        & f \ar{rr}{h} \ar{dl}[swap]{h} \ar[equals,from=dd] & & g \ar[equals]{dl} \\ 
        g \ar[equals,crossing over]{rr} & & g \\ 
        & f \ar{rr}[swap,pos=.25]{h} \ar[equals]{dl} & & g \ar[equals]{uu} \ar[equals]{dl} \\
        f \ar{uu}{h} \ar{rr}[swap]{h} & & g \ar[equals,crossing over]{uu}
    \end{tikzcd}\qquad\begin{tikzcd}
        & g \ar[equals]{rr} \ar[equals]{dl} \ar[from=dd]{}[pos=.25]{h} & & g \ar[equals]{dl} \\ 
        g \ar[equals,crossing over]{rr} & & g \\ 
        & f \ar{rr}[swap,pos=.25]{h} \ar[equals]{dl} & & g \ar[equals]{uu} \ar{dl}{\rev_ih} \\
        f \ar{uu}{h} \ar[equals]{rr} & & f \ar[crossing over]{uu}[swap,pos=.75]{h}
    \end{tikzcd}\]\[\outcompass{i}{i+1}{n+2}
    \]
    To get $\twist_i h$ then with the desired faces, we concatenate $$\left(\fold^{1,1}_{i+1,n+2}\fold^{0,0}_{i,i+1}h\right) \star_{n+2} \left(\fold^{1,0}_{i,n+2}\fold^{0,0}_{i,i+1}h\right)$$ with unitors (\cref{rem:unital}) on all four of the faces in the $i$- and $(i+1)$-directions.    
\end{proof}

%
%

\section{Construction of $\Gpqr$}\label{sec.gpqr}

Typically we will work with lassos along the $1$ direction and null-homotopies of lassos pointing in the $2$ direction (and transitions between those in the $3$ direction), as shown below.
\[
\begin{tikzcd}
    \gamma \rar[equals] & \gamma \\ \gamma \uar[equals] \rar[swap]{L_{ABCD}} \ar[phantom]{ur}{U} & \gamma \uar[equals]
\end{tikzcd}\qquad\begin{tikzcd}
    \gamma \rar[equals] & \gamma \\ \gamma \uar[equals] \rar[swap]{L_{ABCD}} \ar[phantom]{ur}{F} & \gamma \uar[equals]
\end{tikzcd}\qquad\compass{1}{2}
\]

\subsection{Defining $G(p,q)$}\label{subsec: gpq}

\begin{figure}
    \centering
    \includegraphics[width=10cm, trim = {3cm 11cm 6cm 13cm}, clip]{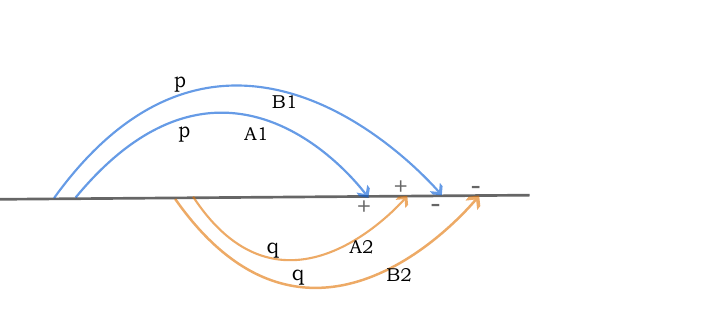}
    \caption{Chord diagram for $G(p,q)$}
    \label{fig: gpqchord}
\end{figure}

Given elements $p, q \in \pi_1(M)$, we depict the chord diagram of the map $G(p,q) \colon I^2 \to \EmbIM$ in  \cref{fig: gpqchord}

It is given by the concatenation shown below using the notation from \cref{defn:flip}. 
\[
\flip{
\begin{array}{|c|c|}
    \hline U_{A_1B_1} & \id\gamma \\
    \hline \id_2 L_{A_1B_1} \prl \id_1 L_{A_2B_2} & \rot_{1,2}U_{A_2B_2} \\
    \hline
\end{array}
}\qquad\qquad\qquad\compass{1}{2}
\]\[
=\begin{array}{|c|c|c|}
    \hline \id\gamma &  U_{A_1B_1} & \id\gamma \\
    \hline \rev_1 \rot_{1,2}U_{A_2B_2} & \id_2 L_{A_1B_1} \prl \id_1 L_{A_2B_2} & \rot_{1,2}U_{A_2B_2} \\
    \hline \id\gamma & \rev_2 U_{A_1B_1} & \id\gamma \\
    \hline
\end{array}\qquad\qquad\qquad{}
\]

\begin{figure}
    \centering
    \includegraphics[width=8cm, trim = {3cm 16.5cm 7cm 3cm}, clip]{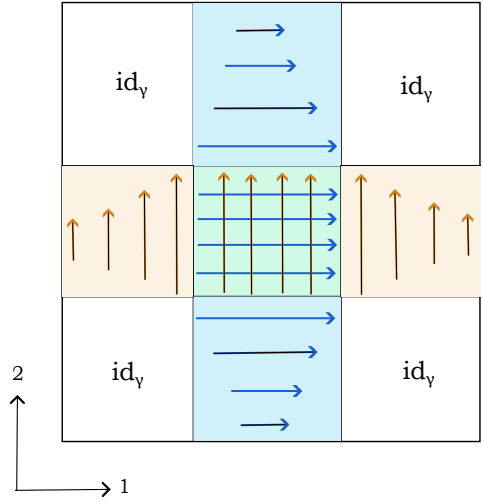}
    \caption{Geometric picture of $G(p,q)$}
    \label{fig: gpq}
\end{figure}

This concatenation diagram can be visualized in \cref{fig: gpq}.

The blue lassos ($L_{A_1, B_1}$) progress in the $1$ direction while the orange lassos ($L_{A_2B_2}$) progress perpendicular to it in the $2$ direction. This allows us to cap off the blue lassos with end homotopies because the orange chords are stationary at the base loop at those squares, and vice versa.

Note that $G(p,q)|_{\partial I^2}$ is constant ($\gamma$), and in \cite{budney} $G(p,q)$ was shown to be non-trivial in $\pi_2\EmbIM$. This was shown by inducing a non-trivial map to $\pi_2 T_3\EmbIM$. However, the induced map to $\pi_2T_2\EmbIM$ is shown to be trivial. 

In general, given a chord diagram with disjoint blue chords $Blue$ and orange chords $Oran$, and null homotopies (via embeddings) of those chords $U_{Blue}$, $U_{Oran}$ respectively, we can define a element of $\pi_2\EmbIM$ given by 

\[
\flip{
\begin{array}{| c | c |}
    \hline U_{Blue} & \id\gamma \\
    \hline \id_2 L_{Blue} \prl \id_1 L_{Oran} & \rot_{1,2}(U_{Oran})  \\
    \hline
\end{array}
}
\]

\subsection{Defining $\Gpqr$}

\begin{figure}
    \centering
    \includegraphics[width=10cm, trim = {3cm 19cm 3cm 2cm}, clip]{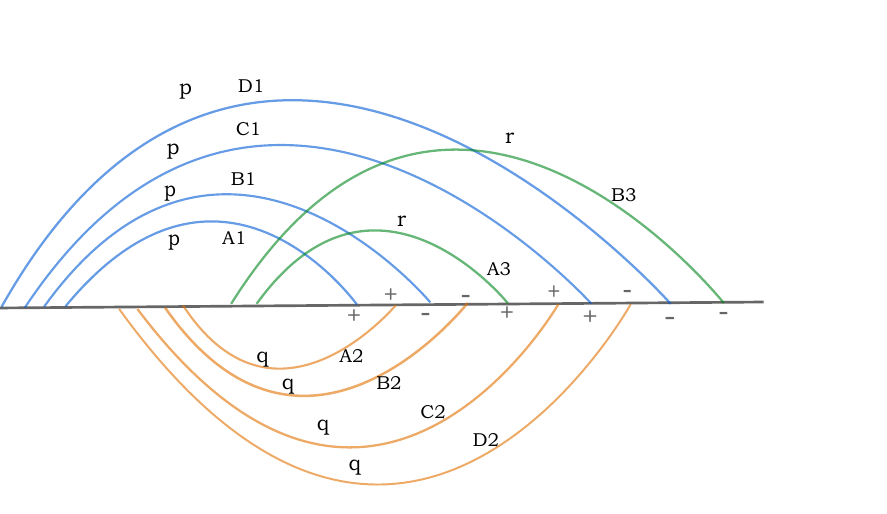}
    \caption{}
    \label{fig: gpqrchord}
\end{figure}

Given elements $p, q, r \in \pi_1(M)$, we now define the map $\Gpqr\colon I^3 \to \EmbIM$ such that $\Gpqr|_{\partial I^3} = \gamma$. We conjecture in \cref{sec.progress} that $\Gpqr$ is non-trivial in $\pi_3\EmbIM$. The chord diagram for $\Gpqr$ is given in \cref{fig: gpqrchord}. 

Let $F_0 \colon I^2 \to \EmbIM$ denote a representative of a cancelling pair of elements of $\pi_2\EmbIM$ that is represented by the chord diagram in \cref{fig: F0 chord diagram} which is obtained from \cref{fig: gpqrchord} with the green chords $A_3, B_3$ removed, and the concatenation diagram in \eqref{eqn.Fzero}.
\begin{equation}\label{eqn.Fzero}
F_0 \coloneq\flip{
\begin{array}{| c | c |}
    \hline U_{A_1B_1} \prl U_{C_1D_1} & \id\gamma \\
    \hline \id_2 L_{A_1B_1C_1D_1} \prl \id_1 L_{A_2B_2C_2D_2} & \rot_{1,2}(U_{A_2B_2}\prl U_{C_2D_2})  \\
    \hline
\end{array}
}
\qquad\compass{1}{2}
\end{equation}

$F_0$ is trivial in $\pi_2(\EmbIM)$ because it is a sum of cancelling elements from \\$\pi_2(\EmbIM)$. This can be shown with a sequence of chord moves from \cite{budney}, but we will describe a specific null homotopy we call the capping null homotopy $H_F$ of this in $\EmbIM$ in \cref{subsec:full basic}. 

The idea of $\Gpqr$ is given by the figure on the left of \cref{fig: gpqr}, while its formal description in terms of concatenation diagrams on the right.

\begin{figure}
    \centering
    \includegraphics[width=15cm, trim = {0cm 7.5cm 1cm 11cm}, clip]{Gpqr.pdf}
    \caption{$G(p,q,r)$}
    \label{fig: gpqr}
\end{figure}

In more detail, the green portion of \cref{fig: gpqr} (which in the center overlaps with the blue and orange) is given by the concatenation diagram in \eqref{eqn.gpqrcenter}
\begin{equation}\label{eqn.gpqrcenter}
\begin{array}{|c|c|c|}
    \hline  \rev_1\fold^{0,0}_{1,2} \rot_{1,2}U_{A_3B_3} & \id_1\rot_{1,2} U_{A_3B_3} &  \fold^{0,0}_{1,2} \rot_{1,2}U_{A_3B_3}\\
    \hline  \rev_1\id_2\rot_{1,2} U_{A_3B_3}   & \id_{2,1} L_{A_3B_3} & \id_2\rot_{1,2} U_{A_3B_3} \\
    \hline  \rev_1 \rev_2\fold^{0,0}_{1,2} \rot_{1,2}U_{A_3B_3}& \rev_2\id_1\rot_{1,2} U_{A_3B_3} &  \rev_2\fold^{0,0}_{1,2} \rot_{1,2} U_{A_3B_3}\\
    \hline
\end{array}
\end{equation}
\[\compass{1}{2}\]
which we have denoted in \cref{defn:revolve} as in \eqref{eqn.gpqrcentersimple}.

\begin{equation}\label{eqn.gpqrcentersimple}
\revolve{\id_{2,1}L_{A_3B_3}}{U_{A_3B_3}}
\qquad\compass{1}{2}
\end{equation}

\begin{figure}
    \centering
    \includegraphics[width = 6cm, trim = {7.5cm 15cm 2.5cm 5cm}, clip]{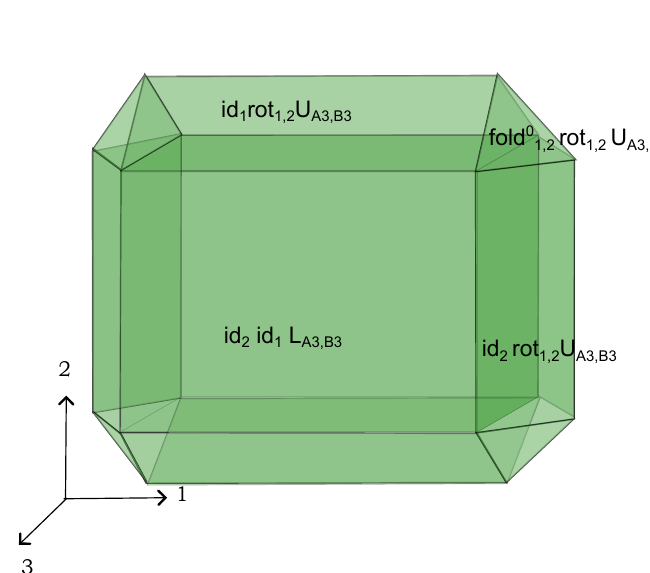}
    \caption{Green portion visualized}
    \label{fig: green}
\end{figure}

We can visualize the pieces in the green portion in \cref{fig: green}

Hence we can define the entirety of $\Gpqr$ following the figure in \cref{fig: gpqr} as a 3-term concatenation in the 3 direction as in \eqref{eqn.gpqrconcat}.
\begin{equation}\label{eqn.gpqrconcat}
\revolve{\rev_3H_F}{\id\gamma}
\;\star_3\;
\revolve{\substack{\id_3F_0 \;\; \prl \\ \id_{2,1} L_{A_3B_3}}}{U_{A_3B_3}}
\;\star_3\;
\revolve{H_F}{\id\gamma}
\end{equation}
\[\compass{1}{2}\]

\subsection{The ``undo-full" transition homotopy $T_{UF}$}\label{subsec:TUF}

Consider the map $L_{ABCD} \colon I \to \EmbIM$ for positive chords $A, C$ and negative chords $B, D$ nested in the order $A, B, C, D$ from outermost to innermost as shown in \cref{fig: full chords}. 

There are two possible null homotopies of this loop given by $U \coloneq U_{AB} \prl U_{CD}$ and $F \coloneq (U_{BC} \; \prl \; L_{AD}) \star_2 U_{AD}$, shown in our cubical diagrams as below. 
\[
\begin{tikzcd}
    \gamma \rar[equals] & \gamma \\ \gamma \uar[equals] \rar[swap]{L_{ABCD}} \ar[phantom]{ur}{U} & \gamma \uar[equals]
\end{tikzcd}\qquad\begin{tikzcd}
    \gamma \rar[equals] & \gamma \\ \gamma \uar[equals] \rar[swap]{L_{ABCD}} \ar[phantom]{ur}{F} & \gamma \uar[equals]
\end{tikzcd}\qquad\compass{1}{2}
\]

We omit the labels of the chords and call this the ``undo-full" transition homotopy $T_{UF}$ because it describes a homotopy from the undo null homotopy $U$ (which is supported in a neighborhood of pairwise zipped bands of $A, B$ and $C, D$) to the full null homotopy $F$ (which is supported in a neighbourhood of the fully zipped bands and lasso disks).

\begin{defn}
    $T_{UF}$ is the transition homotopy from the undo homotopy $ U_{AB} \prl U_{CD}$ to the full null homotopy $U_{BC} \star_2 U_{AD}$. We depict $T_{UF}$ as a concatenation diagram below.
\[
\begin{tikzcd}
    L_{ABCD} \rar{F} & \id\gamma \\
    L_{ABCD} \uar[equals] \rar[swap]{U} \ar[phantom]{ur}{T_{UF}} & \id\gamma \uar[equals]
\end{tikzcd}\qquad\compass{2}{3}
\]

We may denote $\rev_3T_{UF}$ as $T_{FU}$ because it is a homotopy from $F$ to $U$.
\end{defn}

\begin{figure}
    \centering
    \begin{subfigure}{0.38\textwidth}
        \includegraphics[width=\textwidth, trim={3cm 21cm 12cm 5.5cm}, clip]{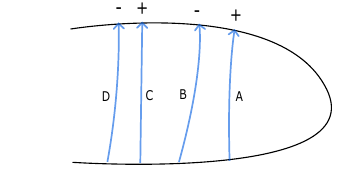}
        \caption{$L_{ABCD}$ deformed}
        \label{fig: LABCD deformed}
    \end{subfigure}
    \hfill
    \begin{subfigure}{0.45\textwidth}
        \includegraphics[width=0.9\textwidth, trim={11.5cm 21cm 2cm 5cm}, clip]{undofullmushrooms.pdf}
        \caption{Zipped bands and disks for $U_{AB} \prl U_{CD}$}
        \label{fig: deformed undo}
    \end{subfigure}
    \hfill
    \begin{subfigure}{\textwidth}
        \includegraphics[width=0.9\textwidth, trim={3cm 14cm 4cm 11cm},clip]{undofullmushrooms.pdf}
        \caption{Fully zipped picture (marked vertices on the right figure)}
        \label{fig: fully zipped LABCD}
    \end{subfigure}
    \caption{Undo Homotopy stages}
    \label{fig: TUF1}
\end{figure}

\begin{figure}
    \centering
    \begin{subfigure}{0.45\textwidth}
        \includegraphics[width=\textwidth, trim={3cm 8cm 12cm 17cm},clip]{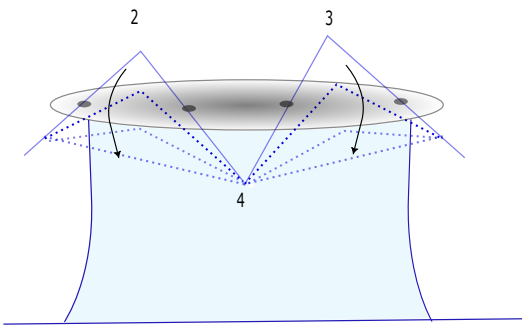}
        \caption{Fully zipped undo homotopy}
        \label{fig: undo vertices}
    \end{subfigure}
        \hfill
    \begin{subfigure}{0.45\textwidth}
        \includegraphics[width=\textwidth, trim={10.5cm 8cm 4cm 17cm},clip]{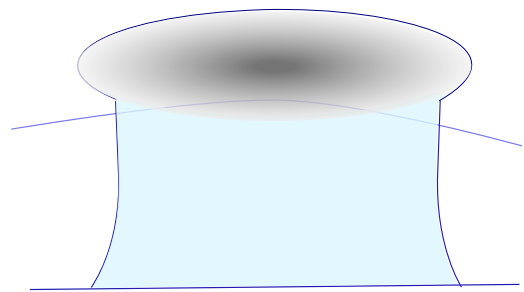}
        \caption{After the fully zipped undo}
        \label{fig: undo vertices end}
    \end{subfigure}
    \caption{Fully Zipped Undo Homotopy}
    \label{fig: TUF1b}
\end{figure}

\begin{figure}
    \centering
    \begin{subfigure}{0.45\textwidth}
        \includegraphics[width=\textwidth, trim={2.5cm 14cm 12cm 11cm}, clip]{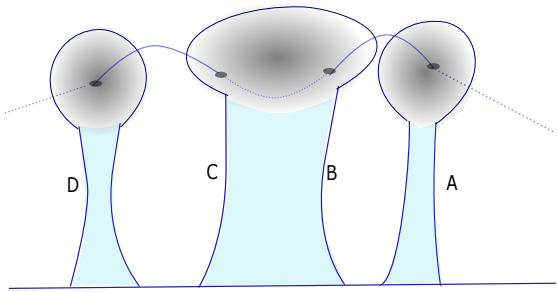}
        \caption{Zipped bands and disks for\\ $U_{BC}\prl \id_2 L_{AD}$}
        \label{fig: deformed full}
    \end{subfigure}
    \hfill
    \begin{subfigure}{0.4\textwidth}
        \includegraphics[width=\textwidth, trim={3.5cm 13.5cm 12cm 12cm},clip]{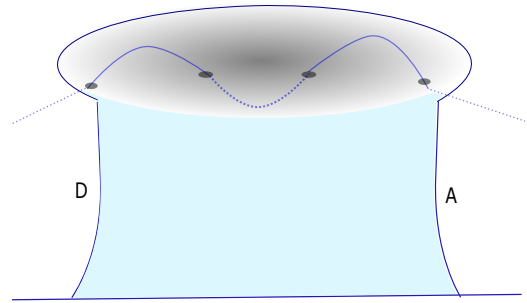}
        \caption{Fully zipped picture}
        \label{fig: fully zipped LABCD again}
    \end{subfigure}
    
    \begin{subfigure}{0.4\textwidth}
        \includegraphics[width=\textwidth, trim={2.5cm 7.5cm 12cm 16cm}, clip]{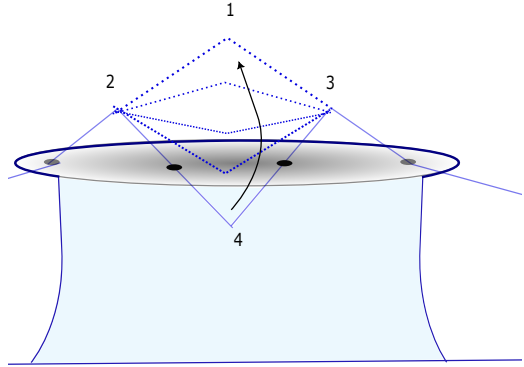}
        \caption{$U_{BC}\prl \id_2 L_{AD}$}
        \label{fig: middle undo}
    \end{subfigure}
    \caption{Fully Zipped Full Homotopy stages}
    \label{fig: TUF2a}
\end{figure}

\begin{figure}
\centering
    \begin{subfigure}{0.4\textwidth}
        \includegraphics[width=\textwidth, trim={11.5cm 14.5cm 4cm 11cm},clip]{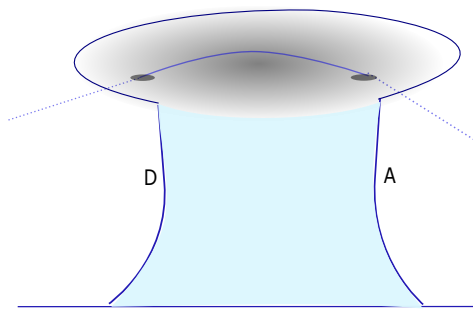}
        \caption{After $U_{BC}\prl \id_2 L_{AD}$}
        \label{fig: after UBC}
    \end{subfigure}
    \hfill
    \begin{subfigure}{0.4\textwidth}
        \includegraphics[width=\textwidth, trim={11.5cm 7cm 2cm 16cm},clip]{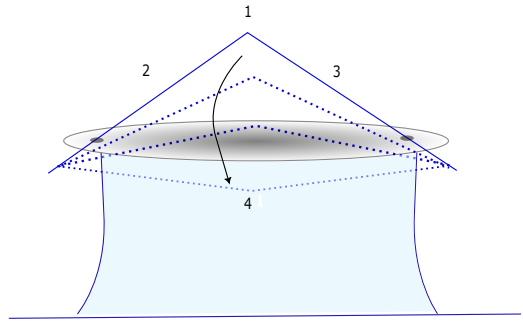}
        \caption{$U_{AD}$}
        \label{fig: undo outer}
    \end{subfigure}
        \hfill
    \begin{subfigure}{0.4\textwidth}
        \includegraphics[width=\textwidth, trim={10.5cm 8cm 4cm 17cm},clip]{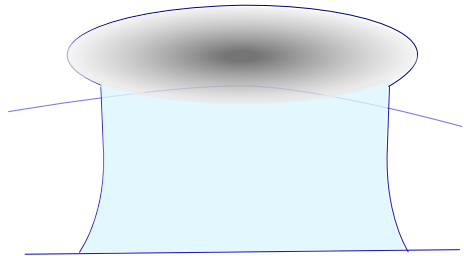}
        \caption{After full null homotopy}
        \label{fig: full vertices end}
    \end{subfigure}
    \caption{Fully Zipped Full Null Homotopy second half}
    \label{fig: TUF2b}
\end{figure}

\begin{figure}
    \centering
    \includegraphics[width=0.9\textwidth, trim={1.5cm 0.5cm 0.5cm 22.5cm},clip]{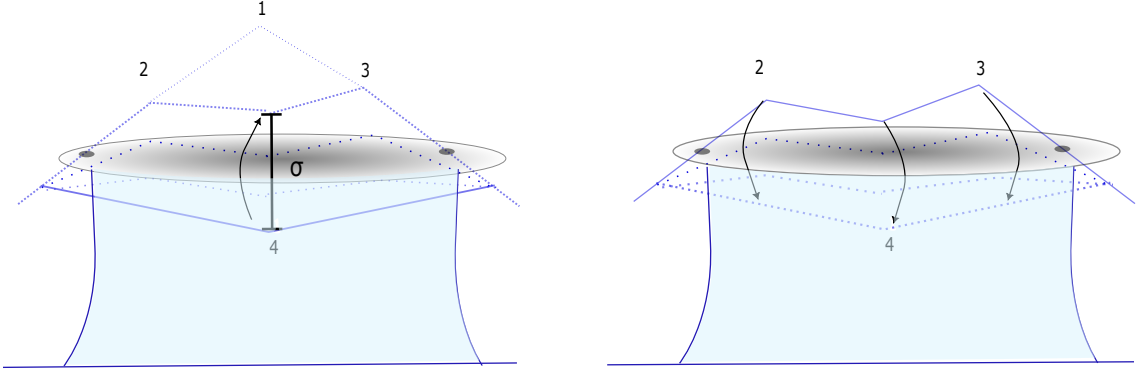}
    \caption{Transition between fully zipped undo and fully zipped full null homotopies}
    \label{fig: TUF}
\end{figure}

\newpage
\begin{itemize}
\item We deform the part of the arc that gets lasso'd around to be situated above the source of the bands as shown in \cref{fig: LABCD deformed}. The zipped pairwise zipped bands and lassos for $U_{AB} \prl U_{CD}$ are shown in \cref{fig: deformed undo}. The zipped bands for the first undo portion of $U_{BC} \;\star_2\; U_{AD}$ is shown in \cref{fig: deformed full}.
\item The first stage is to deform $U_{AB}\; \prl\; U_{CD}$ in the beginning and $U_{BC} \;\star_2\; U_{AD}$ in the end to similar null homotopies where the only difference is that all 4 bands are zipped together and all 4 lasso disks are zipped together. When we zip all the bands and the lasso disks, we see the arc being lasso'd around starts "below" the lasso disc, and pierces the lasso disk four times as shown in the left \cref{fig: fully zipped LABCD}.
\newpage
\item We label the peak of the pierced arc of chords $A, B$ as vertex 2, the peak of the pierced arc of chords $C,D$ as vertex 3, and the lowest point of the pierced arc between chords $B$, $C$ as vertex 4. We also label a point in space as vertex 1 which is the reflection of vertex 4 across the fully zipped lasso disk. See the figure on the right in \cref{fig: fully zipped LABCD}.
\item The ``fully zipped" version of $U_{AB} \;\prl \;U_{CD}$ null homotopy involves translating the vertices labelled 2 and 3 downwards (and the edges attaching to them as well) until the edge joining 2 from the left becomes parallel to the edge joining 2 to 4. and similarly for the edge joining 3 to 4. See \cref{fig: undo vertices}. At the end, we can pull the fully zipped band and lasso back because $I$ does not pierce the lasso disk at this point. See \cref{fig: undo vertices end}. 
\item The fully zipped version of $(U_{BC}\; \prl \;\id_2L_{AD}) \star_2 U_{AD}$ null homotopy first translates vertex 4 until it reaches vertex 1 (like $U_{BC}$) (see \cref{fig: middle undo,fig: after UBC}). We then translates vertex 1 and and all attached edges below the disk (like $U_{AD}$) (see \cref{fig: undo outer,fig: full vertices end}). 
\item Let the distance from vertex 4 to vertex 1 be 1 unit.
\item The homotopy (indexed by $\sigma$) from $U_{AB} \;\prl\; U_{CD}$ to $U_{BC}\;\star_2\; U_{AD}$ translates vertex 4 up by $\sigma$ units (see the left part of \cref{fig: TUF}) and then translates vertices 2, 4, and 3 simultaneously below the lasso disk (see the right part of \cref{fig: TUF}) until they reach the end position described earlier. When $\sigma=0$ this is the $U_{AB} \;\prl\; U_{CD}$ null homotopy and when $\sigma=1$ this is the $(U_{BC}\; \prl \;\id_2L_{AD}) \star_2 U_{AD}$ null homotopy.  
\end{itemize}

\pagebreak

\subsection{The capping null homotopy $H_F$}\label{subsec:full basic}

\begin{figure}
    \centering
    \includegraphics[width=10cm, trim={2cm 2cm 5cm 18cm}, clip]{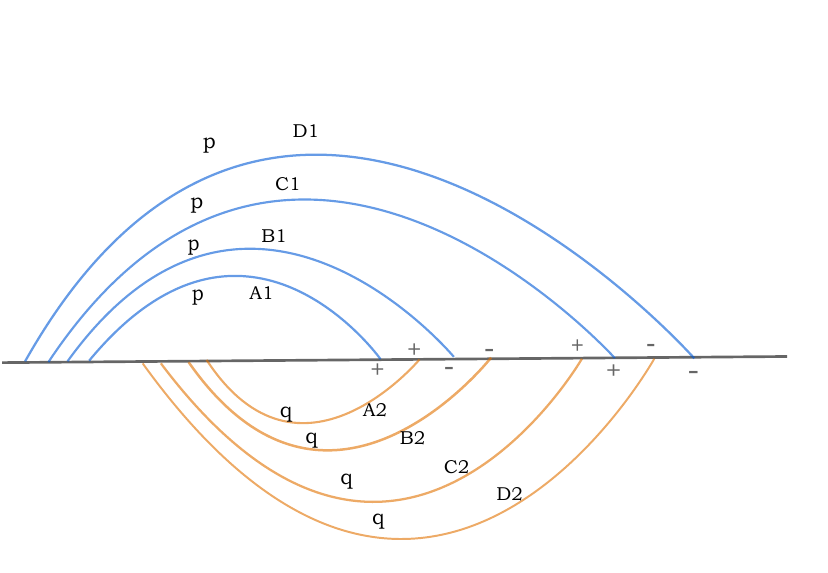}
    \caption{Chord diagram for $F_0$}
    \label{fig: F0 chord diagram}
\end{figure}

We show here that the element of $\pi_2\EmbIM$ constructed by the chord diagram in \cref{fig: F0 chord diagram} is trivial in 
$\pi_2\EmbIM$, and we see this map $I^2\to \EmbIM$ appearing in horizontal $I^2$ slices of $\Gpqr$ in parallel with the green lassos. We describe a specific null homotopy we call the capping null homotopy $H_F$ of this in $\EmbIM$, which is defined by ``capping off" the lassos in this sum with copies of the full null homotopy $F$ in the definition of $\Gpqr$. 

The capping null homotopy is a composition of several homotopies. The main idea is that $L_{A_1B_1C_1D_1} \prl  L_{A_2B_2C_2D_2}$ has a null homotopy where we first apply $U_{B_2C_2}$ followed by $U_{B_1C_1}$, then $U_{A_2D_2}$ and finally $U_{A_1D_1}$. We then need to extend this to the edge squares of the concatenation diagram of $F_0$ which involves transitioning between the undo null homotopies at the edges to the full null homotopy first before folding. 

\begin{figure}
\centering
\includegraphics[width=9cm, trim={2cm 12cm 10cm 9cm}, clip]{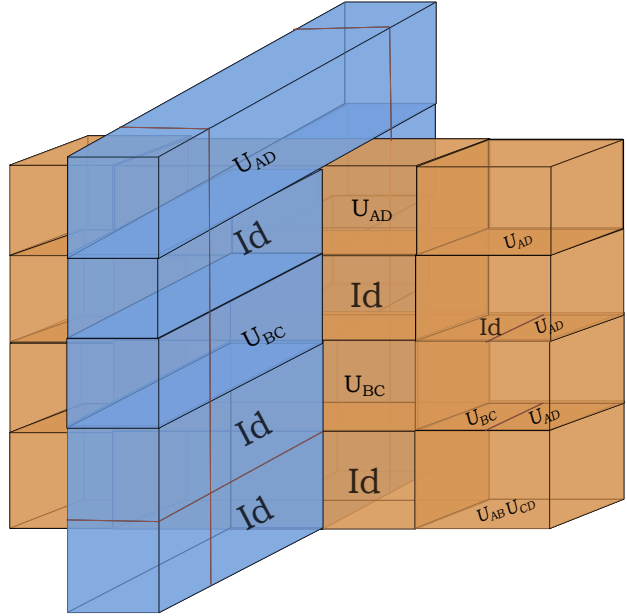}
\caption{Geometric picture of $H_F$}
\label{fig:capping}
\end{figure}

Each stage can be represented by the colored picture in \cref{fig:capping}. The null homotopy can be formally described using five components that, when concatenated, go from $F_0$ to $s$. These five components correspond to the five vertical levels of the picture in \cref{fig:capping} and the five 3-dimensional concatenation diagrams in the sequence represented by \eqref{eqn.capping1} followed by \eqref{eqn.capping2} (namely those depicted as arrows pointing in the 3rd direction).

\begin{equation}\label{eqn.capping1}
\hspace{-0.6cm}
\def\arraystretch{1.1}
\begin{tikzcd}[ampersand replacement=\&, row sep=25ex]
{
\flip{\begin{array}{| c | c |}
    \hline  U_{A_1B_1} \prl U_{C_1D_1} & \id\gamma \\ & \\
    \hline \substack{\id_2 L_{A_1B_1C_1D_1} \; \prl \\ \id_1 L_{A_2B_2C_2D_2}} & \rot_{1,2}(U_{A_2B_2}\prl U_{C_2D_2})  \\
    \hline
\end{array}}}
\dar[shift right=10,outer sep=10]{
\flip{\begin{array}{| c | c |}
    \hline  T_{UF} & \id\gamma \\ 
    & \\
    \hline \id_3\left(\substack{\id_2 L_{A_1B_1C_1D_1} \; \prl \\ \id_1 L_{A_2B_2C_2D_2}}\right) & \rot_{1,2}(T_{UF})  \\
    \hline
\end{array}}}
\dar[shift right=10,swap]{3}
\\
\flip{\begin{array}{| c | c  |}
    \hline  U_{A_1D_1} & \id\gamma  \\
    \cline{1-1} U_{B_1C_1} \prl \id_2 L_{A_1D_1}& \\
    \hline \substack{\id_2 L_{A_1B_1C_1D_1} \; \prl \\ \id_1 L_{A_2B_2C_2D_2}} & \begin{array}{c |c} \substack{\rot_{1,2}U_{B_2C_2} \; \prl \\ \id_1 L_{A_2D_2}} & \rot_{1,2}U_{A_2D_2} \end{array}  \\
    \hline
\end{array}}
\dar[shift right=10,outer sep=10]{
\flip{\begin{array}{| c | c |}
    \hline \id_3 U_{A_1D_1} & \id\gamma  \\
    \cline{1-1}  \id_3(U_{B_1C_1} \prl \id_2 L_{A_1D_1})& \\
    \hline \substack{\id_1 U_{B_2C_2} \; \prl \\ \id_3\left(\substack{\id_2 L_{A_1D_1} \; \prl \\ \id_1 L_{A_2B_2C_2D_2}}\right)} & \begin{array}{c|c} \substack{\rot_{1,2}\fold^{0,0}_{2,3} U_{B_2C_2} \; \prl \\ \id_{3,1} L_{A_2D_2}} & \id_3 \rot_{1,2}U_{A_2D_2}\end{array}\\
    \hline
\end{array}}}
\dar[shift right=10,swap]{3}
\\
\flip{\begin{array}{| c | c |}
    \hline   U_{A_1D_1} & \id\gamma  \\
    \cline{1-1} U_{B_1C_1} \prl \id_2 L_{A_1D_1} & \\
    \hline \substack{\id_2 L_{A_1B_1C_1D_1} \; \prl \\ \id_1 L_{A_2D_2}} & \begin{array}{c|c} \id_1 L_{A_2D_2} & \rot_{1,2}U_{A_2D_2}  \end{array}\\
    \hline
\end{array}}
\end{tikzcd} \hspace{-3cm}\compassss{1}{2}
\end{equation}

\begin{equation}\label{eqn.capping2}
\hspace{-2cm}
\begin{tikzcd}[ampersand replacement=\&, row sep=15ex]
\flip{
\begin{array}{| c | c |}
    \hline   U_{A_1D_1} & \id\gamma  \\
    \cline{1-1} U_{B_1C_1} \prl \id_2 L_{A_1D_1} & \\
    \hline \id_2 L_{A_1B_1C_1D_1} \prl \id_1 L_{A_2D_2} & \begin{array}{c|c} \id_1 L_{A_2D_2} & \rot_{1,2}U_{A_2D_2}  \end{array}\\
    \hline
\end{array}}
\dar[outer sep=10]{
\flip{\begin{array}{| c | c |}
    \hline  \id  & \id\gamma \\
    \cline{1-1} \fold^{0,0}_{2,3} U_{B_1C_1} \prl \;\id_{3,2} L_{A_1D_1} & \\
    \hline \id_2 U_{B_1C_1} \prl \id_3\left(\id_2 L_{A_1D_1} \prl \; \id_1 L_{A_2D_2}\right) & \textbf{unitor}  \\
    \hline
\end{array}}}
\dar[swap]{3}
\\[5ex]
\flip{\begin{array}{| c | c |}
    \hline  U_{A_1D_1} & \id\gamma  \\
    \cline{1-1} \id_2 L_{A_1D_1} & \\
    \hline \id_2 L_{A_1D_1} \prl \id_1 L_{A_2D_2} & \rot_{1,2} U_{A_2D_2}  \\
    \hline
\end{array}}
\dar[outer sep=10]{
\flip{\begin{array}{| c | c |}
    \hline \textbf{unitor} & \id\gamma \\
    \hline \id_1 U_{A_2D_2} \prl \;\id_{3,2} L_{A_1D_1} & \rot_{1,2}\fold^{0,0}_{2,3} U_{A_2D_2}  \\
    \hline
\end{array}}}
\dar[swap]{3}
\\
\flip{\begin{array}{| c | c |}
    \hline  U_{A_2D_2} & \id\gamma  \\
    \hline \id_1L_{A_1D_1} & \id\gamma  \\
    \hline
\end{array}}
\dar[outer sep=10]{
\flip{\begin{array}{| c | c |}
    \hline \fold^{0,0}_{2,3}U_{A_1D_1} & \id\gamma \\
    \hline \id_2 U_{A_1D_1}  & \id\gamma  \\
    \hline
\end{array}}}
\dar[swap]{3}
\\
\flip{\begin{array}{| c | c |}
    \hline  \id\gamma & \id\gamma  \\
    \hline \id\gamma & \id\gamma  \\
    \hline
\end{array}}
\end{tikzcd}\hspace{-3cm}\compasss{1}{2}
\end{equation}

Note that the full null homotopy on the blue part of $F_0$, namely the second column of the form
\[
\def\arraystretch{1.3}
\begin{array}{| c | c | c |}
    \hline \rev_2(U_{A_1B_1} \prl U_{C_1D_1})  & \id_2 L_{A_1B_1C_1D_1} & U_{A_1B_1} \prl U_{C_1D_1}  \\
    \hline
\end{array}\qquad\compass{2}{1}
\]
in \eqref{eqn.capping1} and \eqref{eqn.capping2}, is homotopic to the concatenation diagram in \eqref{eqn.cappingzoom}. 
\begin{equation}\label{eqn.cappingzoom}
\def\arraystretch{1.8}
\begin{tikzcd}[ampersand replacement=\&, row sep=20ex]\begin{array}{| c | c | c |}
    \hline \rev_2 \fold^{0,0}_{2,3} F & \id_2 F & \fold^{0,0}_{2,3} F \\
    \hline \rev_2 T_{UF} & \id L_{A_1B_1C_1D_1} & T_{UF} \\
    \hline
\end{array}\end{tikzcd} \qquad\compass{2}{3}
\end{equation}

The capping null homotopy on the orange part of $F_0$ (the second row) is a rotated version of this (but as we described it as a null homotopy of $F_0$, aspects of the orange and the blue full capping null homotopy need to alternate).

We can use this shorthand to describe the ``blue part" of $\Gpqr$ (which is denoted in blue in \cref{fig: gpqr}) by the concatenation diagram below.

\[
\flip{\begin{array}{| c | c |}
    \hline \id_2 F & \fold^{0,0}_{2,3} F \\
    \hline \id L_{A_1B_1C_1D_1} & T_{UF} \\
    \cline{2-2} & \id_3 U \\ 
    \hline
\end{array}}
\]

\subsection{Other transition homotopies} \label{subsec: transitions}

For Chapter 4, we will require transitions between every pair of $B, U$, and $F$ similar to $T_{UF}$, as well as a transition between these transitions, which we now describe.

\subsubsection{The ``backtrack-undo" transition homotopy $T_{BU}$}\label{subsubsec: TBU}

Suppose $A_1, B_1$ are parallel chords (traversing the same element of $\pi_1(M)$) with opposite signs. The transition homotopy from $B_{A_1B_1}$ to $U_{A_1B_1}$ is denoted by $T_{BU_{A_1B_1}}$ (or $T_{BU}$ in short) is indexed by $\sigma$. We will use $\tau$ to denote the parameter of each of the null homotopies $T_{BU}(\sigma)$.

The backtrack homotopy of $L_{A_1B_1}$ involves pulling back the chords more and more as we go from $\tau =0$ to $\tau =1$. 

The undo homotopy of $L_{A_1B_1}$ involves zipping the bands from time $\tau = 0$ to $\tau = 0.25$, then zipping the lasso disk from $\tau = 0.25$ to $\tau=0.5$ followed by lifting/pulling out the portion of $I$ that pierces the lasso disk out of it from $\tau = 0.5$ to $\tau = 0.75$ and finally retracting the zipped up band and lasso disk from $\tau=0.75$ to $\tau = 1$. 

\begin{itemize}

    \item The backtrack-undo transition involves delaying the backtrack part by adding in the zipping up of the bands before the chords are retracted. So, starting at $\sigma = 0$, no bands are zipped before the chords are backtracked. At $\sigma = 0.25$, the bands are fully zipped before the chords are backtracked.

    \item For $\sigma = 0.25$ to $\sigma = 0.5$, we zip the bands fully and then start zipping up the lasso bit by bit and then retract the chords. So, the lasso disks are unzipped at $\sigma =0.25$ and fully zipped by $\sigma  = 0.5$

    \item For $\sigma = 0.5$ to $\sigma = 0.75$, we begin by zipping the bands and the lassos fully and then pulling out $I$ and then retracting the now-zipped-up chords which is exactly the undo homotopy. 

    \item For $\sigma = 0.75$ to $\sigma = 1$, we begin by zipping the bands and then start zipping up the lasso, pull the piercing chord out, and then retracting the chords a little and then retracting the chords fully.  
\end{itemize}

\begin{defn}
    $T_{BU_{A_1B_1}}$ (or $T_{BU}$) is the transition homotopy from the backtrack homotopy $B_{A_1B_1}$ to the undo null homotopy $U_{A_1B_1}$. We depict $T_{BU}$ as a concatenation diagram below.
\[
\begin{tikzcd}[row sep=large, column sep=large]
    L_{A_1B_1} \rar{U_{A_1B_1}} & \id\gamma \\
    L_{A_1B_1} \uar[equals] \rar[swap]{B_{A_1B_1}} \ar[phantom]{ur}{T_{BU_{A_1B_1}}} & \id\gamma \uar[equals]
\end{tikzcd}\qquad\compass{2}{3}
\]

    We may denote $\rev_3T_{BU}$ as $T_{UB}$ because it is a homotopy from $U$ to $B$.
\end{defn}

\subsubsection{The ``backtrack-full" transition homotopy $T_{BF}$}\label{subsubsec: TBF}

The Full null homotopy of $L_{ABCD}$ where $A,B,C,D$ are parallel chords of alternating signs (like in \cref{fig: full chords}) consists of the undo homotopy on chords $B$ and $C$ followed by the undo homotopy on chords $A$ and $D$. 

In order to transition to the backtrack homotopy $B$ we just concatenate the backtrack-undo transitions on each of the 2 undo homotopies involved in $F$.

\begin{defn}
$T_{BF}$ is the transition homotopy from $B_{ABCD}$ to $F_{ABCD}$, defined as a concatenation diagram below.

\[
\begin{tikzcd}[row sep=huge, column sep=70]
    L_{ABCD} \rar{U_{BC}\;\prl\; \id_2L_{AD}} & L_{AD} \rar{U_{AD}} &\id\gamma \\
    L_{ABCD} \uar[equals] \rar[swap]{B_{BC}\;\prl \;\id_2L_{AD}} \ar[phantom]{ur}{T_{BU_{BC}}\prl \id_{3,2}L_{AD}} & L_{AD} \uar[equals] \rar[swap]{B_{AD}} \ar[phantom]{ur}{T_{BU_{AD}}} & \id\gamma \uar[equals]
\end{tikzcd}\qquad\compass{2}{3}
\]
\end{defn}

\subsubsection{The triple transition homotopies}\label{subsubsec: TBBUF}

We now define a transition homotopy $T_{BBUF}$ mediating between $T_{BU}$, $T_{BF}$, and $T_{UF}$ as in \eqref{TBBUF}.

\begin{equation}\label{TBBUF}
    \begin{tikzcd}[row sep=large, column sep=large]
    U \rar{T_{UF}} & F \\
    B \uar{T_{BU}} \rar[equals]  \ar[phantom]{ur}{T_{BBUF}} & B \uar[swap]{T_{BF}}
\end{tikzcd}\qquad\compass{3}{4}
\end{equation}

We describe $T_{BBUF}$ as a family of homotopies $T_{BBUF}(\sigma)$ from $B$ to $T_{UF}(\sigma)$, where $\sigma \in I$ ranges over the 3-direction. 

\begin{defn}
  We define $T_{BBUF}(\sigma)$ as follows.
  \begin{itemize}
    \item When $\sigma=0$, we have $T_{UF}(0) = U$ and $T_{BBUF}(0) = T_{BU}$. 
    \item When $\sigma=1$, we have $T_{UF}(1) = F$ and $T_{BBUF}(1) = T_{BF}$.
    \item $T_{UF}(\sigma)$ in the beginning involves deforming $U$ to a null homotopy with all bands and lasso disks zipped into 1 band and 1 lasso disk. In each of these stages, we imitate $T_{BU}$ where we increasingly zip the bands before we backtrack (and then increasingy zip the lassos before backtracking). 
    \item The same can be done at the end where $T_{UF}$ involves deforming $F$ to a null homotopy where all bands and lasso disk(s) are fully zipped before pulling out the arc that pierces the lasso disk(s). 
    \item Apart from the fully zipping portion, $T_{UF}(\sigma)$ involves lowering vertex 1 in the arc piercing the lasso disc by $\sigma$ units. To make a transition from $B$. we increase how much of the band gets fully zipped, and then increase how much the lasso gets zipped and increase (as a fraction of $\sigma$) how much vertex $4$ gets pulled up, and finally increase how much of vertices $2,3,4$ get pulled down before backtracking the bands and lassos. 
  \end{itemize}
\end{defn}

Once we define $T_{BBUF}$ we can define $T_{XYZW}$ for any combination of $X, Y, Z, W \in \{B, U, F\}$ by concatenating with appropriate folds and/or composing with rotations, for instance as in \eqref{TUBUF}. 
\begin{equation}\label{TUBUF}
    \begin{tikzcd}[row sep=large, column sep=large]
    U \rar{T_{UF}} & F \\
    U \uar[equals] \rar[swap]{T_{UB}}  \ar[phantom]{ur}{T_{UBUF}} & B \uar[swap]{T_{BF}}
\end{tikzcd}
=
\begin{tikzcd}[row sep=large, column sep=70]
    U \rar[equals] & U \rar{T_{UF}} & F \\
    U \rar[swap]{T_{UB}} \uar[equals] \ar[phantom]{ur}{\fold^{1,0}_{3,4}T_{UB}} & B \uar[swap]{T_{BU}} \rar[equals]  \ar[phantom]{ur}{T_{BBUF}} & B \uar[swap]{T_{BF}}
\end{tikzcd}
\qquad\compass{3}{4}
\end{equation}

%
%

\section{$\Gpqr$ is Null Homotopic in $\pi_3(T_3(\EmbIM))$}\label{sec.nullhomotopy}

\subsection{Null homotopy of $G(p,q)$ in $T_2(\EmbIM)$}\label{subsec: Gpqnull}

This section describes the null homotopy from \cite{budney} in language developed in this thesis. These ideas will be extended in the subsequesnt sections for $\Gpqr$

We have $G(p,q)\colon I^2 \to \EmbIM$, which induces a map $T_2G(p,q) \colon I^2 \times C_2\lr{I} \to C_2\lr{M}$. To show that $G(p,q)$ is trivial in $\pi_2T_2\EmbIM$, we need to construct a null homotopy $N^\ast\colon I \times I^2 \times C_2\lr{I}  \to C_2\lr{M}$ of $T_2G(p,q)$.

$N_B$ and $N_U$ are null homotopies of $G(p,q)$ in $\ImmIM$ which apply the backtrack (and respectively undo) homotopies on the lasso portions of $G(p,q)$. 

To be precise,
\[
\hspace{-1cm}
N_B \coloneq\qquad
\begin{array}{c}
\flip{
\begin{array}{|c|c|}
    \hline T_{UB_{A_1B_1}} & \id\gamma \\
    \hline \id_3(\id_2 L_{A_1B_1} \prl \id_1 L_{A_2B_2}) & \rot_{1,2}T_{UB_{A_2B_2}} \\
    \hline
\end{array}
}
\\\ast_3 
\\\flip{
\begin{array}{|c|c|}
    \hline \fold^{0,0}_{2,3}B_{A_1B_1} & \id\gamma \\
    \hline \id_2 B_{A_1B_1} \prl \id_1 B_{A_2B_2} & \rot_{1,2}\fold^{0,0}_{2,3}B_{A_2B_2} \\
    \hline
\end{array}
}
\end{array}
\compass{1}{2}
\]

\[
\hspace{-1cm}
N_U \coloneq \qquad 
\begin{array}{c}
\flip{
\begin{array}{|c|c|}
    \hline \id_3U_{A_1B_1} & \id\gamma \\
    \hline \id_3(\id_2 L_{A_1B_1} \prl \id_1 L_{A_2B_2}) & \rot_{1,2}\id_3U_{A_2B_2} \\
    \hline
\end{array}
}
\\\ast_3
\\\flip{
\begin{array}{|c|c|}
    \hline \fold^{0,0}_{2,3}U_{A_1B_1} & \id\gamma \\
    \hline \id_2 U_{A_1B_1} \prl \id_1 U_{A_2B_2} & \rot_{1,2}\fold^{0,0}_{2,3}U_{A_2B_2} \\
    \hline
\end{array}
}
\end{array}
\compass{1}{2}
\]

and 

\[
\hspace{-1cm}
N_{BU} \coloneq\quad 
\begin{array}{c}
\flip{
\begin{array}{|c|c|}
    \hline \fold^{1,0}_{3,4}T_{UB_{A_1B_1}} & \id\gamma \\
    \hline \id_{3,4}(\id_2 L_{A_1B_1} \prl \id_1 L_{A_2B_2}) & \rot_{1,2}\fold^{1,0}_{3,4}T_{UB_{A_2B_2}} \\
    \hline
\end{array}
}
\\\ast_3
\\\flip{
\begin{array}{|c|c|}
    \hline \fold^{0,0}_{2,3}T_{BU_{A_1B_1}} & \id\gamma \\
    \hline \id_2 T_{BU_{A_1B_1}} \prl \id_1 T_{BU_{A_2B_2}} & \rot_{1,2}\fold^{0,0}_{2,3}T_{BU_{A_2B_2}} \\
    \hline
\end{array}
}
\end{array}
\compass{1}{2}
\]

We will consider the interval as partitioned (up to overlapping endpoints) into the sub-intervals $I_1,I_{1'},I_2,I_{2'},I_3$, listed in order. The chords $A_1, B_1$ originate at $I_1$, chords $A_2, B_2$ originate at $I_2$, and all chords lasso around points on $I_3$.  

We first define $N\colon C_2\lr{I} \to \text{ Map }(I\times I^2, \ImmIM)$ which takes $(p_1, p_2)$ to $N_{p_1, p_2}$ which will be among $N_B, N_U$ or $N_{BU}(t)$ (which is an intermediate stage of the transition homotopy from $N_B$ to $N_U$). The main feature of this is that if $(p_1, p_2) \in I_a\times I_b$ (where $1\leq a \leq b \leq 3$), then $T_2N_{p_1, p_2}$ is well defined when restricted to that specific $I_a\times I_b$.  For example, when $p_1\in I_1$ and $p_2\in I_2$, $N_U$ doesn't map points in $I_1$, $I_2$ to distinct points in $M$ because the undo homotopies collide $I_1$ and $I_2$ when done simultaneously. However, we don't see intersections between $I_1$ and $I_3$ (and $I_2$ and $I_3$) because neither of $p_1, p_2$ is in $I_3$, so $N_B$ maps points in $I_1, I_2$ to distinct points in $M$. Thus in \eqref{eqn: NstarGpq}, we see $N_B$ in the square $I_1\times I_2$. The complete $N$ is defined in the concatenation diagram \eqref{eqn: NstarGpq}. The directions for $C_2\lr I $ are $4$ and $5$ because $N_X$ (for $X=U, B$ use up directions $1,2,3$

This allows us to define the null homotopy  $N^\ast$ that we want as $$N^\ast(t, a, b, p_1, p_2) = N_{p_1, p_2}(t)(a,b)(p_1, p_2).$$

\begin{equation}\label{eqn: NstarGpq}
    \begin{array}{|c||*{5}{c|}}
        \hline  
                & I_1     & I_{1'}     & I_2      & I_{2'}     & I_3  \\\hline
	I_3:    & \id_{4,5} N_U & \id_{4,5} N_U    & \id_{4,5} N_U  & \id_{4,5} N_U    & \id_{4,5} N_U  \\\cline{1-6}
	I_{2'}: &  \id_4 N_{BU} &  \id_4N_{BU}    &  \id_4 N_{BU}  &  \id_4 N_{BU}           \\\cline{1-5}
	I_2:    & \id_{4,5} N_B & \id_{4,5} N_B    &\id_{4,5} N_B                  \\\cline{1-4}
	I_{1'}  & \id_{4,5} N_B & \id_{4,5} N_B                           \\\cline{1-3}
	I_1     & \id_{4,5} N_B                                 \\\cline{1-2}
    \end{array}\qquad\compass{4}{5}
\end{equation}

\subsection{Null homotopy of $G(p,q,r)$ in $T_3\EmbIM$}\label{subsec: define N}

In this section, we use similar ideas to \cref{subsec: Gpqnull} to define the null homotopy of $G(p,q,r)$ in $\pi_3T_3\EmbIM$. 
The element $\Gpqr\colon I^3 \to \EmbIM$ is null homotopic in $\ImmIM$. We will use three such null homotopies: ``Back-track" $N_B$, ``Undo" $N_U$, and ``Full" $N_F$. 
\[
\Gpqr \Tto{N_B}{4} \id\gamma \qquad \Gpqr \Tto{N_U}{4} \id\gamma \qquad \Gpqr \Tto{N_F}{4} \id\gamma 
\]
that we define in \cref{subsec: NBNUNF}. We also define homotopies between each pair 
\[
N_B \Tto{N_{BU}}{5} N_U \qquad N_U \Tto{N_{UF}}{5} N_F \qquad N_B \Tto{N_{BF}}{5} N_F
\]
as well as homotopies between these homotopies such as
\[
\begin{tikzcd}
    N_F \rar{N_{FU}} & N_U \\
    N_B \rar[equals,swap]{\id N_B} \uar{N_{BF}} \ar[phantom]{ur}{N_{BBFU}} & N_B \uar[swap]{N_{BF}}
\end{tikzcd}
\qquad\qquad\compass{5}{6}
\]
in \cref{subsec: major transitions}, where we also write $N_{YX}$ for $\rev_5 N_{XY}$.

\begin{figure}
    \centering
    \includegraphics[width=8cm, trim={4cm 19cm 1.8cm 2cm}, clip]{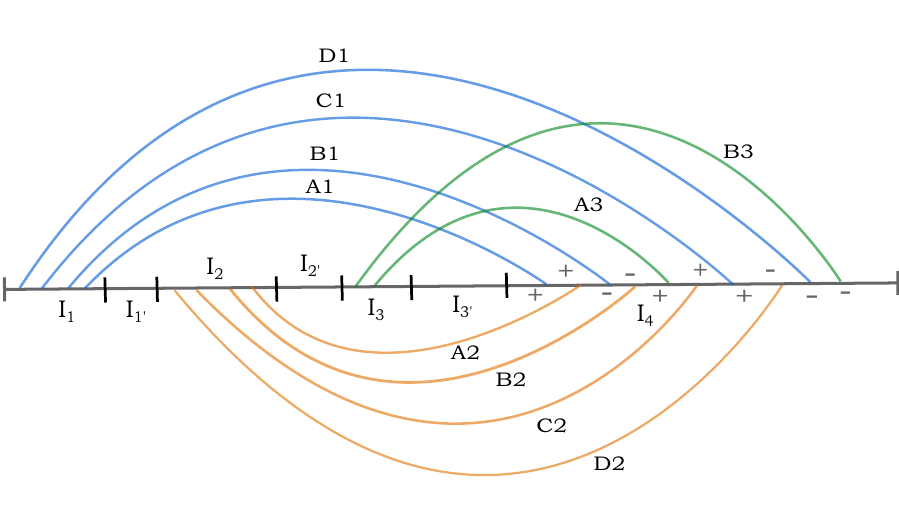}
    \caption{Division of $I$ into subintervals $I_1$ through $I_4$}
    \label{fig: subintervals}
\end{figure}

We will consider the interval as partitioned (up to overlapping endpoints) into the sub-intervals $I_1,I_{1'},I_2,I_{2'},I_3,I_{3'},I_4$, listed in order. The chords $A_1, B_1, C_1, D_1$ originate at $I_1$, chords $A_2, B_2, C_2, D_2$ originate at $I_2$, and chords $A_3, B_3$ originate at $I_3$. All 10 chords lasso around points on $I_4$.  See \cref{fig: subintervals}.

First, we define $N\colon C_3\lr{I} \to \text{ Map }(I\times I^3, \EmbIM)$ where $(p_1, p_2, p_3) \mapsto N_{p_1, p_2, p_3}$ such that if $p_1, p_2, p_3$ are in specified intervals as in table \eqref{tab: intervals for null hom}, $N_{p_1, p_2, p_3}$ is as specified in the rightmost column.

\begin{equation}\label{tab: intervals for null hom}   
\begin{array}{ |c|c|c|c| } 
 \hline 
 p_1 & p_2 & p_3 & \text{ Null homotopy } N_{p_1 p_2 p_3}\\ 
 \hline
 I_1/I_2/I_3 & I_1/I_2/I_3 & I_1/I_2/I_3 & N_B\\ 
 \hline
 I_2/I_3/I_4 & I_2/I_3/I_4 & I_4 & N_U\\ 
 \hline
 I_1 & I_1/I_3/I_4 & I_4 & N_U\\ 
 \hline
 I_1 & I_2 & I_4 & N_F\\ 
 \hline
\end{array}
\end{equation}

\begin{figure}
    \centering
    \includegraphics[width=15cm, trim ={0.5cm 4cm 0 11.5cm}, clip]{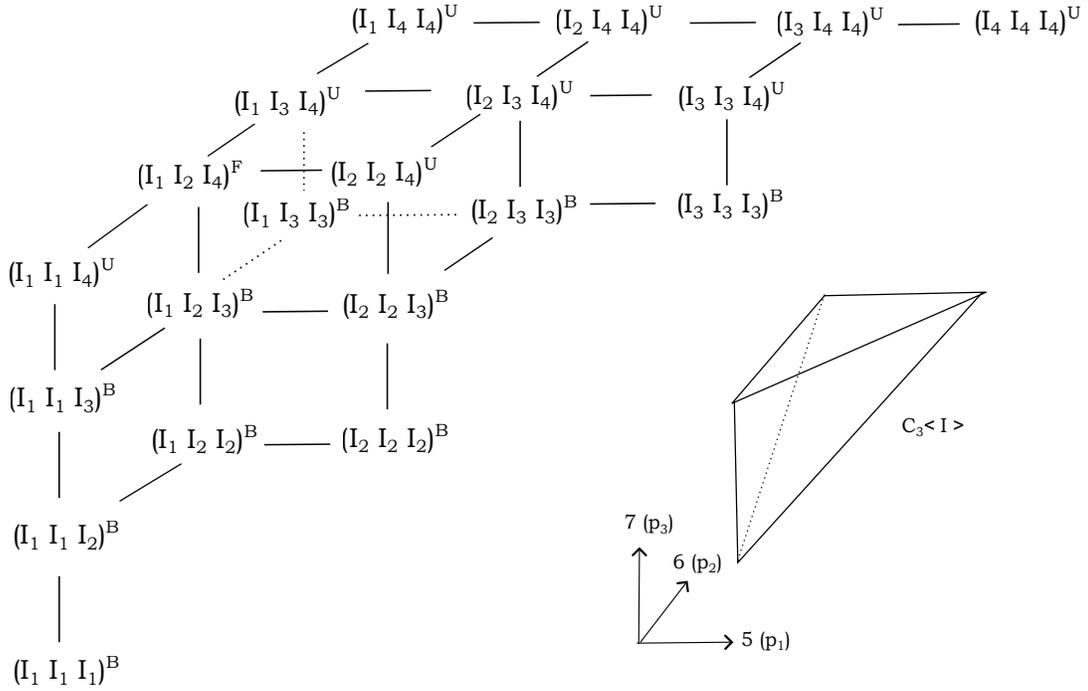}
    \caption{Adjacency graph of products $I_a\times I_b\times I_c$ with their null homotopies}
    \label{fig: adjacency graph}
\end{figure}

\cref{fig: adjacency graph} is a visualization of each of the products of intervals $I_1, I_2, I_3, I_4$ in $C_3\lr{I}$ and the superscript is $N_{p_1, p_2, p_3}$ from Table \eqref{tab: intervals for null hom}. This allows us to see which transition homotopies we will need to extend $N$ to the entirety of $C_3\lr{I}$.

Now we will extend $N$ to the rest of $C_3\lr{I}$ as a composition of concatenation diagrams below. We will assume the directions of $C_3\lr{I}$ are $5,6,7$ respectively because $N_X$ (for $X\in \{B, U, F \}$) uses directions $1$ through $4$. Furthermore, $N_{XY}$ is a homotopy from $N_X$ to $N_Y$ in the $5$ direction (which will be used to `fill edges' in \cref{fig: adjacency graph}), and $N_{XYZW}$ is a homotopy in the $6$ direction between two of these homotopies (where $X,Y,Z,W \in \{B,U,F\}$), which will be used to `fill faces'. 

The concatenation diagram \eqref{N} for $N$, is written from bottom up to form $C_3\lr{I}$. Recall that $\id N_X$ means $\id_{7,6,5}N_X$ for $X=\{U, B, F\}$ and $\id N_{XY}$ means $\id_{7,6}N_{XY}$.

\clearpage
\begin{equation}\label{N}
\def\arraystretch{.8}
\begin{array}{cc}
    \hspace{-2.5cm}p_3\in I_{1} \hspace{3cm} p_3\in I_{1'} \hspace{4cm} p_3\in I_{2} & \\
    \begin{array}{|c||c|}
        \hline
	    & I_1    \\\hline
	I_1 & \id N_B    \\\hline
    \end{array} \quad \star_7 \quad

    \begin{array}{|c||*{2}{c|}}
        \hline  
                & I_1 & I_{1'}   \\\hline
        I_{1'}: & \id N_B & \id N_B      \\\cline{1-3}
	I_{1}:  & \id N_B           \\\cline{1-2}
    \end{array} \quad \star_7 \quad

    \begin{array}{|c||*{3}{c|}}
        \hline  
                & I_1 & I_{1'} & I_2  \\\hline
        I_2:    & \id N_B & \id N_B    & \id N_B  \\\cline{1-4}
	I_{1'}  & \id N_B & \id N_B           \\\cline{1-3}
	I_1     & \id N_B                 \\\cline{1-2}
    \end{array} & \\
    &\\
    \hspace{-1cm}p_3\in I_{2'} \hspace{4cm} p_3\in I_{3} & \\
    \star_7 \quad
    \begin{array}{|c||*{4}{c|}}
        \hline  
                & I_1 & I_{1'} & I_2 & I_{2'} 
        \\\hline
        I_{2'}: & \multicolumn{4}{c|}{} 
        \\\cline{1-1}\cline{5-5}
	I_2:  & \multicolumn{3}{c|}{\id N_B}
        \\\cline{1-1}\cline{4-4}
	I_{1'}: & \multicolumn{2}{c|}{}
        \\\cline{1-1}\cline{3-3}
	I_1:    &  
        \\\cline{1-2}
    \end{array} \quad \star_7 \quad 

    \begin{array}{|c||*{5}{c|}}
        \hline  
                & I_1 & I_{1'} & I_2  & I_{2'} & I_3  \\\hline
	I_3:    & \multicolumn{5}{c|}{} 
        \\\cline{1-1}\cline{6-6}
	I_{2'}: & \multicolumn{4}{c|}{\id N_B} 
        \\\cline{1-1}\cline{5-5}
	I_2:  & \multicolumn{3}{c|}{}
        \\\cline{1-1}\cline{4-4}
	I_{1'}: & \multicolumn{2}{c|}{}
        \\\cline{1-1}\cline{3-3}
	I_1:    &  
        \\\cline{1-2}
    \end{array}  &  \\
    \hspace{-3cm} p_3\in I_{3'} & \\
    \hspace{-5cm}\star_7\quad
     \begin{array}{|c||*{6}{c|}}
        \hline  
                & I_1                & I_{1'}       & I_2       & I_{2'}     & I_3  & I_{3'}
        \\\hline
        I_{3'}: & \multicolumn{6}{c|}{}     
        \\\cline{7-7}\cline{1-1}
	I_3:    & \multicolumn{5}{c|}{\id N_{BU}}            \\\cline{6-6}\cline{1-3}
	I_{2'}: & \multicolumn{2}{c|}{} & \multicolumn{2}{c|}{}
        \\\cline{1-2}\cline{5-5}
	I_2:  & \id N_{BF} &  N_{BBFU} & {}
        \\\cline{4-4}\cline{1-2}
	I_{1'}: & \multicolumn{2}{c|}{}
        \\\cline{1-3}
	I_1:    & \id N_{BU}                                                        \\\cline{1-2}
    \end{array} & \hspace{-5cm}\compass{5}{6} \\
    p_3\in I_{4} & \\
    
    \star_7 \quad \begin{array}{|c||*{7}{c|}}
        \hline  
                & I_1 & I_{1'}                      & I_2             & I_{2'}              & I_3 & I_{3'} & I_4 \\\hline
        I_4:    & \multicolumn{7}{c|}{} 
        \\\cline{8-8}\cline{1-1}
        I_{3'}: & \multicolumn{6}{c|}{\id N_U}     
        \\\cline{7-7}\cline{1-1}
	I_3:    & \multicolumn{5}{c|}{}                
        \\\cline{6-6}\cline{1-3}
	I_{2'}: & \multicolumn{2}{c|}{} & \multicolumn{2}{c|}{}
        \\\cline{1-2}\cline{5-5}
	I_2:  & \id N_F & \id N_{FU} & {}
        \\\cline{4-4}\cline{1-2}
	I_{1'}: & \multicolumn{2}{c|}{}
        \\\cline{1-3}
	I_1:    & \id N_U                                                        \\\cline{1-2}
    \end{array} & 

\end{array}
\end{equation}

\begin{defn}
Given $N_X\colon I \to \text{ Map }(I^3, \ImmIM)$, define $N_X^\ast\colon I \times I^3 \times I^3 \to M^3$ as the induced map 
\[
N_X^\ast(t,a,b,c,p_1,p_2,p_3) = \bigg(
N_X(t)(a,b,c)(p_1),
N_X(t)(a,b,c)(p_2),
N_X(t)(a,b,c)(p_3)
\bigg)
\]
that applies the immersion to tuples $(p_1, p_2, p_3)$ in $I^3$. 
\end{defn}

 The null homotopies $N_B^\ast,N_U^\ast,N_F^
 \ast\colon I \times I^3 \times I^3\to M^3$  land in $C_3\lr{M}$ when we restrict to certain products $I \times I^3 \times I_a\times I_b\times I_c$ (see Table \eqref{tab: intervals for null hom} of subintervals $I_1,I_2,I_3,I_4$). For example, on the block $I_1\times I_2\times I_2 \subset I^3$, $N_B$ sends distinct triples $(p_1\times p_2\times p_3)$ to distinct triples in $M$: even though the backtrack homotopy contains non-embedded intervals in general, $p_3$ is in $I_2$ and never in $I_4$. So the image of $I_1\times I_2\times I_2$ in $M^3$ doesn't detect the self intersection of the immersion (which is only seen in the product $I_a\times I_b\times I_4$ where either $a$ or $b$ is $1, 2,$ or $3$).  

We use the map $N\colon C_3\lr{I} \to \text{ Map }(I\times I^3, \EmbIM)$ to define a map $N^\ast \colon I \times I^3 \times C_3\lr I \to C_3\lr M$ which will be our null homotopy of $T_3G(p,q,r)$. 

\begin{defn}
We now define $N^\ast\colon I\times I^3 \times C_3\lr{I} \to C_3 \lr{M}$ as $$N^\ast(t,a, b,c, p_1, p_2, p_3) = N_{p_1,p_2,p_3}(t)(a, b,c) (p_1, p_2, p_3)$$ 
\end{defn}

\subsection{The homotopies $N_B$, $N_U$, $N_F$}\label{subsec: NBNUNF}

In \cref{defn: undo,defn: backtrack,defn: full} we defined the undo, backtrack, and full null homotopies for lassos. Using them, we show in this section that the map $\Gpqr\colon I^3 \to \EmbIM$ is null homotopic in $\ImmIM$. We construct three such null homotopies: ``Back-track" $N_B$, ``Undo" $N_U$, and ``Full" $N_F$ using the homotopies in the earlier subsections mentioned. 

\subsubsection{The undo homotopy $N_U$ of $\Gpqr$}

On the green portion of $\Gpqr$, we define the undo null homotopy of $\Gpqr$ by applying the undo null homotopy of $L_{A_3B_3}$ to the center of the concatenation and folding the undo homotopy around the outside as in \eqref{eqn.NUgreen}.

\begin{equation}\label{eqn.NUgreen}
\def\arraystretch{1.0}
\begin{tikzcd}[ampersand replacement=\&, row sep=20ex]{\revolve{\id_{2,1}L_{A_3B_3}}{U_{A_3B_3}}}
\dar[outer sep=10]{\revolve{\id_{2,1}U_{A_3B_3}}{\fold^{0,0}_{2,3}U_{A_3B_3}}}
\dar[swap]{4}
\\
{\revolve{\id\gamma}{\id\gamma}}
\end{tikzcd} \qquad\compassss{1}{2}
\end{equation}

We define the undo homotopy on the entire blue part of $\Gpqr$ in \eqref{eqn.NUblue}. The undo homotopy for the orange is the same and we would put them together perpendicular to each other using the $\rot$ operation as we did when defining $G(p,q)$ and $F_0$. 

\begin{equation}\label{eqn.NUblue}
\hspace{-1.0cm}
\begin{tikzcd}[ampersand replacement=\&, row sep=35ex]
\flip{\begin{array}{| c | c |}
    \hline \id_2 F & \fold^{0,0}_{2,3} F \\
    \hline \id L_{A_1B_1C_1D_1} & T_{UF} \\
    \cline{2-2} & \id_3 U \\ 
    \hline
\end{array}}
\dar[outer sep=10]{
\flip{\begin{array}{| c | c |}
    \hline \id_2 T_{FU} & \fold^{0,0}_{2,3} T_{FU} \\
    \hline \id L_{A_1B_1C_1D_1} & \fold^{1,0}_{3,4}T_{UF} \\
    \cline{2-2} & \id_{4,3} U\\
    \hline
\end{array}}
}
\dar[swap]{4}
\\
\flip{\begin{array}{| c | c |}
    \hline \id_2 U & \fold^{0,0}_{2,3} U \\
    \hline \id L_{A_1B_1C_1D_1} & \id_3 U \\
    \cline{2-2} & \id_3 U \\
    \hline
\end{array}}
\rar[equals]{\sim}
\&
{\revolve{\id L_{A_1B_1C_1D_1}}{U}}
\dar[swap,outer sep=10]{{\revolve{\id U}{\fold^{0,0}_{2,3}U}}}
\dar{4}
\\[-10ex]
\& {\revolve{\id\gamma}{\id\gamma}}
\end{tikzcd}
\hspace{-1cm}\compassss{2}{3}
\end{equation}

\subsubsection{The full homotopy $N_F$ of $\Gpqr$}

The Full Homotopy $N_F$ for the green portion can be either $N_U$ or $N_B$ (because we only need the full end homotopy when the green chords are not visible in the induced map $C_3\lr{I} \to C_3\lr{M}$. We shall choose $N_U$.  

For the blue part of $\Gpqr$ we define $N_F$ as follows. (the orange part is done similarly but with a $1-2$ rotation with the chords $A_1, B_1, C_1, D_1$ replaced by $A_2, B_2, C_2, D_2$),

\[
\def\arraystretch{1.0}
\begin{tikzcd}[ampersand replacement=\&, row sep=22ex]
\flip{\begin{array}{| c | c |}
    \hline  \id_2 F & \fold^{0,0}_{2,3} F \\
    \hline  \id L_{A_1B_1C_1D_1} & T_{UF} \\
    \cline{2-2} & \id_3 U \\ 
    \hline
\end{array}}
\dar[outer sep=10]{
\flip{\begin{array}{| c | c |}
    \hline \id_{4,2} F & \id_4 \fold^{0,0}_{2,3} F \\
    \hline \id L_{A_1B_1C_1D_1} & \fold^{0,0}_{3,4} T_{UF} \\
    \cline{2-2} & \id_3 T_{UF}\\
    \hline
\end{array}}}
\dar[swap]{4}\\
\flip{\begin{array}{| c | c |}
    \hline \id_2 F & \fold^{0,0}_{2,3} F \\
    \hline  \id L_{A_1B_1C_1D_1} & \id_3 F \\
    \cline{2-2} & \id_3 F \\
    \hline
\end{array}}
\rar[equals]{\sim}
\&
{\revolve{\id L_{A_1B_1C_1D_1}}{F}}
\dar[outer sep=10,swap]{\revolve{\id F}{\fold^{0,0}_{2,3}F}}
\dar{4}
\\
\& {\revolve{\id\gamma}{\id\gamma}}
\end{tikzcd} \hspace{-1cm}\compassss{2}{3}
\]

\subsubsection{The backtrack homotopy $N_B$ of $\Gpqr$}

On the green part of $\Gpqr$, we first transition the $U$ border of $L_{A_3B_3}$ to $B$ and then apply $B$ to the center and we fold $B$ on the border.

\[
\begin{tikzcd}[ampersand replacement=\&, row sep=17ex]
{\revolve{\id_{2,1}L_{A_3B_3}}{U}}
\dar[outer sep=10]{\revolve{\id_{4,2,1}L_{A_3B_3}}{T_{UB}}}
\dar[swap]{4}
\\
{\revolve{\id_{2,1}L_{A_3B_3}}{B}}
\dar[outer sep=10]{\revolve{\id_{2,1}B}{\fold^{0,0}_{2,3}B}}
\dar[swap]{4}
\\
{\revolve{\id\gamma}{\id\gamma}}
\end{tikzcd} \qquad\compass{1}{2}
\]

On the blue part (and orange by rotation), we define $N_B$ as follows.

\[
\def\arraystretch{1.0}
\begin{tikzcd}[ampersand replacement=\&, row sep=22ex]
\flip{\begin{array}{| c | c |}
    \hline  \id_2 F & \fold^{0,0}_{2,3} F \\
    \hline  \id L_{A_1B_1C_1D_1} & T_{UF} \\
    \cline{2-2} & \id_3 U \\ 
    \hline
\end{array}}
\dar[outer sep=10]{
\flip{\begin{array}{| c | c |}
    \hline \id_2 T_{FB} & \fold^{0,0}_{2,3} T_{FB} \\
    \hline \id L_{A_1B_1C_1D_1} & T_{UFBB} \\
    \cline{2-2} & \id_3 T_{UB}\\
    \hline
\end{array}}}
\dar[swap]{4}\\
\flip{\begin{array}{| c | c |}
    \hline \id_2 B & \fold^{0,0}_{2,3} B \\
    \hline \id L_{A_1B_1C_1D_1} & \id_3 B \\
    \cline{2-2} & \id_3 B \\
    \hline
\end{array}}
\dar[outer sep=10]{\revolve{\id B}{\fold^{0,0}_{2,3}B}}
\dar[swap]{4}
\\
{\revolve{\id\gamma}{\id\gamma}}
\end{tikzcd} \qquad\compass{2}{3}
\]

\subsection{Transition homotopies between $N_U, N_B$, and $N_F$}\label{subsec: major transitions}

\subsubsection{The Undo-Full transition homotopy $N_{UF}$ from $N_U$ to $N_F$}\label{subsec: NUF}

The full-undo transition is needed when either the orange chords are not seen ($p_1$ is in the intervals $I_1'$ or $p_2$ is in $I_2'$) or the blue chords aren't seen ($p_1$ is in the interval $I_2'$ or $p_2$ is in the interval $I_3'$). 

The essential piece in this transition is the homotopy presented in \cref{subsec:TUF} which is what we shall do on the lasso portion. For the borders, we show how to transition from $N_U$ to $N_F$ using piece-wise transitions and concatenations.

\[
\hspace{-.2cm}
\begin{tikzcd}[ampersand replacement=\&, row sep=27ex]
\;\;\;\;\;
\flip{\begin{array}{| c | c |}
    \hline \id_2 T_{FU} & \fold^{0,0}_{2,3} T_{FU} \\
    \hline \id L_{A_1B_1C_1D_1} & \fold^{1,0}_{3,4}T_{UF} \\
    \cline{2-2} & \id_{4,3} U\\
    \hline
\end{array}} 
 \quad\star_4\quad  {\revolve{\id_{3,2}U}{\fold^{0,0}_{2,3}U}}
    \dar[shift right=30,outer sep=10]
    {
    \flip{\begin{array}{| c | c |}
        \hline \id_2\fold^{1,0}_{3,4}T_{FU} & \fold^{0,0}_{2,3}\fold^{1,0}_{3,4}T_{FU} \\
        \hline \id L_{A_1B_1C_1D_1} & \rot_{4,5}\twist_3 T_{UF} \\
        \cline{2-2} & \id_3\fold^{1,1}_{3,4} T_{UF}\\
        \hline
    \end{array}}
    \;\star_4\;   
    \revolve{\id_{3,2} T_{UF}}{\fold^{0,0}_{2,3}T_{UF}}
    }
    \dar[shift right=30,swap]{5}
\\
\flip{\begin{array}{| c | c |}
    \hline \id_{4,2} F & \fold^{0,0}_{2,3} \id_3 F \\
    \hline \id L_{A_1B_1C_1D_1} & \fold^{0,0}_{3,4} T_{UF} \\
    \cline{2-2} & \id_3 T_{UF}\\
    \hline
\end{array}}
\quad\star_4\quad
\revolve{\id_{3,2}F}{\fold^{0,0}_{2,3}F}
\end{tikzcd}
\qquad\hspace{-3cm}\compassss{2}{3}
\]
The map $\twist_3 T_{UF}$ is defined in \cref{lem:twist}.

The top right squares of the first piece of the source, the target, and the morphism can be described as $\fold^{0,0}_{2,3}$ of the square 
\[
\begin{tikzcd}[row sep=huge, column sep=70]
    F \rar[equals]{\id F} & F \\
    F \uar[equals]{\id F} \rar[swap]{T_{FU}}  \ar[phantom]{ur}{ \fold^{1,0}_{3,4}T_{FU}} & U \uar[swap]{T_{UF}}
\end{tikzcd} \qquad \compass{3}{4} 
\]
which becomes 
\[
\begin{tikzcd}[row sep=huge, column sep=70]
    \fold^{0,0}_{2,3}F \rar[equals]{\id\fold^{0,0}_{2,3} F} & \fold^{0,0}_{2,3}F \\
    \fold^{0,0}_{2,3} F \uar[equals]{\id\fold^{0,0}_{2,3} F} \rar[swap]{\fold^{0,0}_{2,3} T_{FU}}  \ar[phantom]{ur}{\fold^{0,0}_{2,3}\fold^{1,0}_{3,4}T_{FU}} & \fold^{0,0}_{2,3}U \uar[swap]{\fold^{0,0}_{2,3}T_{UF}}
\end{tikzcd} \qquad\compass{4}{5}
\]

\subsubsection{The transition homotopy $N_{BU}$ from $N_B$ to $N_U$}\label{subsec: NBU}

We now show how to transition from $N_B$ to $N_U$.

\[
\hspace{-.2cm}
\begin{tikzcd}[ampersand replacement=\&, row sep=27ex]
\;\;\;\;\;{
\flip{\begin{array}{| c | c |}
    \hline \id_2 T_{FB} & \fold^{0,0}_{2,3} T_{FB} \\
    \hline \id L_{A_1B_1C_1D_1} & T_{UFBB} \\
    \cline{2-2} & \id_3 T_{UB}\\
    \hline
\end{array}} 
\quad\star_4\quad  
\revolve{\id_{3,2}B}{\fold^{0,0}_{2,3}B}}
    \dar[shift right=30,outer sep=10]
    {    
    \flip{
    \begin{array}{| c | c |}
        \hline \id_2 T_{FBFU} & \fold^{0,0}_{2,3} T_{FBFU} \\
        \hline \id L_{A_1B_1C_1D_1} & 
        \substack{
            \fold^{1,0}_{3,4}T_{BBUF} 
            \\ \star_4 \\
            \fold^{1,0}_{4,5}T_{UFBB}
        } \\
        \cline{2-2} & \id_3\fold^{1,0}_{3,4}T_{UB}\\
        \hline
    \end{array}}
    \;\star_4\;   
    \revolve{\id_{4,3,2} T_{BU}}{\fold^{0,0}_{2,3}T_{BU}}
    }
\dar[shift right=30,swap]{5}
\\
\flip{\begin{array}{| c | c |}
    \hline \id_2 T_{FU} & \fold^{0,0}_{2,3} T_{FU} \\
    \hline \id L_{A_1B_1C_1D_1} & \fold^{1,0}_{3,4}T_{UF} \\
    \cline{2-2} & \id_{4,3} U\\
    \hline
\end{array}}
\quad\star_4\quad
\revolve{\id_{3,2}U}{\fold^{0,0}_{2,3}U} 
\end{tikzcd}
\qquad\hspace{-3cm}\compassss{2}{3}
\]

\subsubsection{The transition homotopy $N_{BF}$ from $N_B$ to $N_F$}\label{subsec: NBF}

We now show how to transition from $N_B$ to $N_F$.

\[
\hspace{-.2cm}
\begin{tikzcd}[ampersand replacement=\&, row sep=27ex]
\;\;\;\;\;
{
\flip{\begin{array}{| c | c |}
    \hline \id_2 T_{FB} & \fold^{0,0}_{2,3} T_{FB} \\
    \hline \id L_{A_1B_1C_1D_1} & T_{UFBB} \\
    \cline{2-2} & \id_3 T_{UB}\\
    \hline
\end{array}}
\quad\star_4\quad  
\revolve{\id_{3,2}B}{\fold^{0,0}_{2,3}B}
}
\dar[shift right=30,outer sep=10]
    {
    \flip{\begin{array}{| c | c |}
    \hline \id_2 \fold^{1,0}_{3,4}T_{FB} & \fold^{0,0}_{2,3}\fold^{1,0}_{3,4}T_{FB} \\
    \hline \id L_{A_1B_1C_1D_1} & 
    \substack{
        \fold^{0,0}_{3,4}T_{BBUF}
        \\ \star_4 \\
        \fold^{1,0}_{4,5}T_{UFBB}
    } \\
    \cline{2-2} & \id_3 T_{UBUF}\\
    \hline
    \end{array}}
    \;\star_4\; 
    \revolve{\id_{4,3,2} T_{BF}}{\fold^{0,0}_{2,3}T_{BF}}
    }
\dar[shift right=30,swap]{5}
\\
\flip{\begin{array}{| c | c |}
    \hline \id_{4,2} F & \fold^{0,0}_{2,3}\id_3 F \\
    \hline \id L_{A_1B_1C_1D_1} & \fold^{0,0}_{3,4} T_{UF} \\
    \cline{2-2} & \id_3 T_{UF}\\
    \hline
\end{array}}
\quad\star_4\quad
\revolve{\id_{3,2}F}{\fold^{0,0}_{2,3}F}
\end{tikzcd}
\qquad\hspace{-3cm}\compassss{2}{3}
\]

\subsubsection{Transition of transitions}\label{subsec: NUFBB}

We now describe a transition homotopy $N_{BUBF}$ from $N_{BU}$ to $N_{BF}$ such that on the sides it has $\id N_B $ and $N_{UF}$

\[
\begin{tikzcd}[ampersand replacement=\&, row sep=27ex]
\!\!\!\!\!\!\!\!
{
\flip{
    \begin{array}{| c | c |}
        \hline \id_2 T_{FBFU} & \fold^{0,0}_{2,3} T_{FBFU} \\
        \hline \id L_{A_1B_1C_1D_1} & 
        \substack{
            \fold^{1,0}_{3,4}T_{BBUF} 
            \\ \star_4 \\
            \fold^{1,0}_{4,5}T_{UFBB}
        } \\
        \cline{2-2} & \id_3\fold^{1,0}_{3,4}T_{UB}\\
        \hline
    \end{array}}
    \;\star_4\;   
    \revolve{\id_{4,3,2} T_{BU}}{\fold^{0,0}_{2,3}T_{BU}}
}
\dar[shift right=62,outer sep=10]
    {
    \flip{\begin{array}{| c | c |}
    \hline \id_2 \left(\substack{
        \id_5\fold^{1,0}_{3,4}T_{FB}
        \\ \star_3 \\
        \fold^{01}_{3,5}T_{BBFU}
    }\right) 
    & \fold^{0,0}_{2,3} \left(\substack{
        \id_5\fold^{1,0}_{3,4}T_{FB}
        \\ \star_3 \\
        \fold^{0,1}_{3,5}T_{BBFU}
    }\right) \\
    \hline \id L_{A_1B_1C_1D_1} & 
    \substack{
        \rot_{4,6}\twist_3 T_{BBUF}
        \\ \star_4 \\
        \id_6\fold^{1,0}_{4,5}T_{UFBB}
    } \\
    \cline{2-2} & \id_3 \left(\substack{
        \id_5\fold^{1,0}_{3,4}T_{UB}
        \\ \star_3 \\
        \fold^{1,0}_{3,5}T_{BBUF}
    }\right) \\
    \hline
    \end{array}}
    \;\star_4\; 
    \revolve{\id_{4,3,2} T_{BUBF}}{\fold^{0,0}_{2,3}T_{BUBF}}
    }
\dar[shift right=62,swap]{6}
\\
\flip{\begin{array}{| c | c |}
    \hline \id_2 \fold^{1,0}_{3,4}T_{FB} & \fold^{0,0}_{2,3}\fold^{1,0}_{3,4}T_{FB} \\
    \hline \id L_{A_1B_1C_1D_1} & 
    \substack{
        \fold^{0,0}_{3,4}T_{BBUF}
        \\ \star_4 \\
        \fold^{1,0}_{4,5}T_{UFBB}
    } \\
    \cline{2-2} & \id_3 T_{UBUF}\\
    \hline
    \end{array}}
    \;\star_4\; 
    \revolve{\id_{4,3,2} T_{BF}}{\!\!\!\!\!\!\!\!\fold^{0,0}_{2,3}T_{BF}}
\end{tikzcd}
\hspace{-2.4cm}\compassss{2}{3}
\]
The map $\twist_3 T_{BBUF}$ is defined in \cref{lem:twist} and has as its target in the 5-direction the map $\twist_3 T_{UF}$ from \cref{subsec: NUF}.


In \cref{subsec: define N}, we use $N_{BBFU}$ which is $\rev_5 \rot_{5,6} N_{BUBF}$.

\section{$\pi_3\EmbIS$ via Generators and Relations in $\pi_7(\CSB{4})$}\label{sec.piCk}

In this chapter, $\pi_nC_k\lr{M}$ denotes the rational homotopy groups unless specified. 

\subsection{Background on $\pi_m(\CSB{k})$}\label{subsec: proppinck}

We recall the relations satisfied by the rational generators $\twij{i}{p}{i}{j}$ of $\pi_{k}\CSB{k}$ from \cite{budney}.
\begin{itemize}
    \item $w_{ij}=(-1)^{k+1}w_{ji}$, which becomes $w_{ij}=w_{ji}$ when $k=3$.
    \item $[w_{ij}, w_{jk}] = [w_{jk}, w_{ki}]=[w_{ki}, w_{ij}]= - [w_{jk}, w_{ij}]$
    \item Jacobi identity: $[A, [B, C]]+[B, [C, A]] + [C, [A, B]]=0$ (when $A, B, C$ have the same degree)
    \item $t_1^{a_1} t_2^{a_2}\dots  t_m^{a_m}w_{ij} = t^{a_i - a_j}_iw_{ij} = t_j^{a_j - a_i } w_{ij}$
    \item $[w_{ij}, w_{kl}]=0$ where $\{i,j\} \cap \{ k, l \} = \phi$.
\end{itemize}

We will say that the ``cyclic shifts" of $[A, [B, C]]$ are $[B, [C, A]]$ and $[C,[A, B]]$.

\subsection{Generators of $\pi_7(\CSB{3})$}

\noindent The linearly independent generators of $\pi_7(\CSB{3})$ are given by:
\begin{enumerate}[(A)]
    \item $[\twij{1}{p}{1}{2}, [\twij{1}{q}{1}{2}, \twij{1}{r}{1}{2}]]$
    \item $[\twij{1}{p}{1}{2}, [\twij{1}{q}{1}{2}, \twij{2}{r}{2}{3}]]$ and a cyclic shift $[\twij{2}{p}{2}{3}, [\twij{1}{q}{1}{2}, \twij{1}{r}{1}{2}]]$
    \item $[\twij{1}{p}{1}{2}, [\twij{2}{q}{2}{3}, \twij{2}{r}{2}{3}]]$ and a cyclic shift $[\twij{2}{p}{2}{3}, [\twij{1}{q}{1}{2}, \twij{2}{r}{2}{3}]]$ 
    \item $[\twij{2}{p}{2}{3}, [\twij{2}{q}{2}{3}, \twij{2}{r}{2}{3}]]$
    \item $[\twij{1}{p}{1}{3}, [\twij{1}{q}{1}{3}, \twij{1}{r}{1}{3}]]$
\end{enumerate}
All other generators can be shown to be a linear combination of the above by using Jacobi and other relations repeatedly.

\subsection{Generators and relations of $\pi_7\CSB{4}$}\label{subsec: relations}

R is the subgroup of $\pi_7(\CSB{4})$ generated by torsion and the images of the maps
$$\pi_7(\CSB{3})\hookrightarrow \pi_7(\CSB{4})$$ 
induced by the 5 boundary faces of $\CSB{4}$: 
$$p_1=\ast, \; p_1=p_2, \; p_2=p_3, \; p_3=p_4, \sand p_4=\ast.$$
We now describe the relations on elements of $\pi_7(\CSB{4})$ that arise from quotienting by $R$.
\begin{itemize}
    \item The face $p_1=\ast$ gives us that 
$$[\twij{i}{p}{i}{j}, [\twij{k}{q}{k}{l}, \twij{m}{r}{m}{n}]] = 0$$ 
for $i, j, k, l, m, n \in \{2, 3, 4\}$. 
    \item The face $p_4=\ast$ gives the same when $i, j, k, l, m, n \in \{1, 2, 3\}$.
    \item From the face $p_1=p_2$, $t_1 \mapsto t_1t_2$, $t_2 \mapsto t_3$, and $t_3 \mapsto t_4 $. It hence maps 
$$\twij{1}{p}{1}{2} \mapsto \twij{1}{p}{1}{3} + \twij{2}{p}{2}{3},$$
$$\twij{2}{q}{2}{3} \mapsto \twij{3}{q}{3}{4},$$ 
$$\twij{1}{p}{1}{3} \mapsto \twij{1}{p}{1}{4}+ \twij{2}{p}{2}{4}.$$ 
Note that generators (A) and (D) map to relations already obtained from faces $p_4=\ast$ and $p_1= \ast$ respectively. 

We start with the first generator of (B), where we have
\begin{equation}\label{eqn.Bgen}
[\twij{1}{p}{1}{2}, [\twij{1}{q}{1}{2}, \twij{2}{r}{2}{3}]] \mapsto \left(\begin{array}{cl}
& [\twij{1}{p}{1}{3}, [\twij{1}{q}{1}{3}, \twij{3}{r}{3}{4}]] \\
+ & [\twij{2}{p}{2}{3}, [\twij{1}{q}{1}{3}, \twij{3}{r}{3}{4}]] \\
+ & [\twij{1}{p}{1}{3}, [\twij{2}{q}{2}{3}, \twij{3}{r}{3}{4}]] \\
+ & [\twij{2}{p}{2}{3}, [\twij{2}{q}{2}{3}, \twij{3}{r}{3}{4}]] 
\end{array}\right)
\end{equation}
\[
=  
\left(\begin{array}{cl}
& [\twij{1}{p}{1}{3}, [\twij{1}{q}{1}{3}, \twij{3}{r}{3}{4}]] \\
+ & [\twij{2}{p}{2}{3}, [\twij{1}{q}{1}{3}, \twij{3}{r}{3}{4}]] \\
+ & [\twij{1}{p}{1}{3}, [\twij{2}{q}{2}{3}, \twij{3}{r}{3}{4}]] 
\end{array}\right)
\]
because the bottom term comes from $p_1=\ast$.

We also obtain a relation similar to a relation in \cite{budney} by cyclically shifting the relation in \eqref{eqn.Bgen} to get 
\[
\left(\begin{array}{cl}
& [\twij{3}{r}{3}{4}, [\twij{1}{p}{1}{3}, \twij{1}{q}{1}{3}]] \\
+ & [\twij{3}{r}{3}{4}, [\twij{2}{p}{2}{3}, \twij{1}{q}{1}{3}]] \\
+ & [\twij{3}{r}{3}{4}, [\twij{1}{p}{1}{3}, \twij{2}{q}{2}{3}]] 
\end{array}\right)
\]
\[
=
\left[\twij{3}{r}{3}{4}, \left([\twij{1}{p}{1}{3}, \twij{1}{q}{1}{3}] + (t_2^p t_1^q - t_1^p t_2^q) [w_{23}, w_{13}]\right)\right] 
=0
\]
 
Now we do the second generator of (C), where we have
\[
\begin{array}{c}
[\twij{2}{p}{2}{3}, [\twij{1}{q}{1}{2}, \twij{2}{r}{2}{3}]] \mapsto \\ 
{[\twij{3}{p}{3}{4}, [\twij{1}{q}{1}{3}, \twij{3}{r}{3}{4}]] + [\twij{3}{p}{3}{4}, [\twij{2}{q}{2}{3}, \twij{3}{r}{3}{4}]]} \\
= [\twij{3}{p}{3}{4}, [\twij{1}{q}{1}{3}, \twij{3}{r}{3}{4}]]
\end{array}
\]
because the latter is 0 from $p_1=\ast$. 

So, this gives us 
\begin{equation}\label{eqn.134is0}[\twij{3}{p}{3}{4}, [\twij{1}{q}{1}{3}, \twij{3}{r}{3}{4}]]=0
\end{equation}
as well as  
$$[\twij{1}{p}{1}{3}, [\twij{3}{q}{3}{4}, \twij{3}{r}{3}{4}]]=0$$ 
(by using the Jacobi relation). 

The generator (E) will be dealt with later. 
    \item The face $p_3=p_4$ works almost analogously to the $p_1=p_2$ face. This face inclusion maps $t_1 \mapsto t_1$, $t_2 \mapsto t_2$, $t_3 \mapsto t_3t_4 $. It hence maps
$$\twij{2}{p}{2}{3} \mapsto \twij{2}{p}{2}{3} + \twij{2}{p}{2}{4},$$
$$\twij{1}{q}{1}{2} \mapsto \twij{1}{q}{1}{2},$$ 
$$\twij{1}{p}{1}{3} \mapsto \twij{1}{p}{1}{3}+ \twij{1}{p}{1}{4}.$$ 
Note that generators (A) and (D) map to relations already obtained from faces $p_4=\ast$ and $p_1= \ast$ respectively. 

We start with the first generator of (B). 
\begin{equation}\label{eqn: 3Term0b}
\begin{array}{c}
[\twij{1}{p}{1}{2}, [\twij{1}{q}{1}{2}, \twij{2}{r}{2}{3}]] \mapsto \\
{[\twij{1}{p}{1}{2}, [\twij{1}{q}{1}{2}, \twij{2}{r}{2}{3}]] + [\twij{1}{p}{1}{2}, [\twij{1}{q}{1}{2}, \twij{2}{r}{2}{4}]]} \\
= [\twij{1}{p}{1}{2}, [\twij{1}{q}{1}{2}, \twij{2}{r}{2}{4}]]
\end{array}
\end{equation}

This gives $ [\twij{1}{p}{1}{2}, [\twij{1}{q}{1}{2}, \twij{2}{r}{2}{4}]]=0$ as well as  $ [\twij{2}{p}{2}{4}, [\twij{1}{q}{1}{2}, \twij{1}{r}{1}{2}]]=0$ by the Jacobi identity.

Now the first generator of (C). 
\begin{equation*}
[\twij{1}{p}{1}{2}, [\twij{2}{q}{2}{3}, \twij{2}{r}{2}{3}]] \mapsto 
\left(\begin{array}{cl}
& [\twij{1}{p}{1}{2}, [\twij{2}{q}{2}{3}, \twij{2}{r}{2}{3}]] 
\\ + & [\twij{1}{p}{1}{2}, [\twij{2}{q}{2}{3}, \twij{2}{r}{2}{4}]] 
\\ + & [\twij{1}{p}{1}{2}, [\twij{2}{q}{2}{4}, \twij{2}{r}{2}{3}]] 
\\ + & [\twij{1}{p}{1}{2}, [\twij{2}{q}{2}{4}, \twij{2}{r}{2}{4}]] 
\end{array} \right)
\end{equation*}
\[
 =
 \left(\begin{array}{cl}
 &[\twij{1}{p}{1}{2}, [\twij{2}{q}{2}{3}, \twij{2}{r}{2}{4}]] 
 \\ + & [\twij{1}{p}{1}{2}, [\twij{2}{q}{2}{4}, \twij{2}{r}{2}{3}]] 
 \\ + & [\twij{1}{p}{1}{2}, [\twij{2}{q}{2}{4}, \twij{2}{r}{2}{4}]]
\end{array} \right)
\]

(The first term is 0 from the $p_4 = \ast $ face)


Using the Jacobi identity to cyclically shift that relation we get $$[\twij{2}{p}{2}{4}, [\twij{1}{q}{1}{2}, \twij{2}{r}{2}{4}]] + [\twij{2}{p}{2}{3}, [\twij{1}{q}{1}{2}, \twij{2}{r}{2}{4}]] + [\twij{2}{p}{2}{4}, [\twij{1}{q}{1}{2}, \twij{2}{r}{2}{3}]]=0 .$$

We will need further algebraic manipulation to this relation for later. 
\begin{equation}\label{eqn: Break3}
\begin{array}{rrl}
 &&[\twij{2}{p}{2}{4}, [\twij{1}{q}{1}{2}, \twij{2}{r}{2}{4}]]  
\\= & - &
[\twij{2}{p}{2}{3}, [\twij{1}{q}{1}{2}, \twij{2}{r}{2}{4}]] 
\\ & - &[\twij{2}{p}{2}{4}, [\twij{1}{q}{1}{2}, \twij{2}{r}{2}{3}]]  
\\= & - & [t_2^{-q}\twij{3}{-p-q}{2}{3},[\twij{2}{-q}{1}{2}, t_2^{-q}\twij{4}{-r-q}{2}{4}]] 
\\ & - & [t_2^{-q}\twij{4}{-p-q}{2}{4}, [\twij{2}{-q}{1}{2}, t_2^{-q}\twij{3}{-r-q}{2}{3}]]  
\\= & -& (t_2^{-q} t_3^ {-p-q} t_4^{-q-r}) [w_{23}, [w_{12}, w_{24}]] 
\\ & - & (t_2^{-q} t_3^ {-q-r} t_4^{-p-q}) [w_{24}, [w_{12}, w_{23}]] 
\\= & - & (t_2^{-q} t_3^ {-p-q} t_4^{-q-r})[w_{23}, [w_{14}, w_{12}]] 
\\ & - & (t_2^{-q} t_3^ {-q-r} t_4^{-p-q}) [w_{24}, [w_{13}, w_{12}]] 
\\= & + & t_2^{-q} t_3^ {-p-q} t_4^{-q-r}\left( [w_{14}, [w_{12}, w_{23}]] + [w_{12}, [w_{23}, w_{14}]]\right)
\\ & + & t_2^{-q} t_3^ {-q-r} t_4^{-p-q}\left( [w_{13}, [w_{12}, w_{24}]] + [w_{12}, [w_{24}, w_{13}]]  \right)
\\= & + & (t_2^{-q} t_3^ {-p-q} t_4^{-q-r}) [w_{14}, [w_{13}, w_{12}]] 
\\ & + & (t_2^{-q} t_3^ {-q-r} t_4^{-p-q}) [w_{13}, [w_{14}, w_{12}]] 
\end{array}
\end{equation}
We use the Jacobi identity at the 5th equal sign which changes the sign of the whole expression. For the last equality, we can delete the 2nd and 4th term because they have a whitehead product of the form $[w_{ij}, w_{kl}]$ where $\{i,j\} \cap \{ k, l \} = \phi$.

The generator (E) will be done later. 
\end{itemize}

\begin{rem}
So far, we have shown that $[\twij{2}{p}{2}{4}, [\twij{1}{q}{1}{2}, \twij{2}{r}{2}{4}]] $ and  $[\twij{1}{p}{1}{3}, [\twij{1}{q}{1}{3}, \twij{3}{r}{3}{4}]]$ (and their cyclic shifts) can be written as a sum of terms with all four indices. 
\\[0.2cm]On the other hand, $[\twij{1}{p}{1}{2}[\twij{1}{q}{1}{2}, \twij{2}{r}{2}{4}]] $ and  $[\twij{3}{p}{3}{4}, [\twij{1}{q}{1}{3}, \twij{3}{r}{3}{4}]]$ (and their cyclic shifts) are 0. 
\\[0.2cm]Furthermore, any term which has only three indices (say $\{1, 3, 4\}$ or $\{1, 2, 4\}$ are generated by the above terms with three indices.  Thus we can conclude so far that, $\pi_7(\CSB{4})/R$ is generated by just the terms with all four indices included.
\end{rem}

\begin{itemize}
    \item The face $p_2=p_3$ maps $t_1 \mapsto t_1$, $t_2 \mapsto t_2t_3$, and  $t_3 \mapsto t_4 $. It hence maps $$\twij{2}{p}{2}{3} \mapsto \twij{2}{p}{2}{4} + \twij{3}{p}{3}{4}, $$ 
$$\twij{1}{q}{1}{2} \mapsto \twij{1}{q}{1}{2}+ \twij{1}{p}{1}{3},$$
$$\twij{1}{p}{1}{3} \mapsto  \twij{1}{p}{1}{4}.$$ 

Note that generators (A) and (D) map to relations already obtained from faces $p_4=\ast$ and $p_1= \ast$ respectively. 
\\Generator (E) maps to $[\twij{1}{p}{1}{4}, [\twij{1}{q}{1}{4}, \twij{1}{r}{1}{4}]] $ making that zero. 

Let 
$$\twij{2}{p}{2}{3} \mapsto  \twij{2}{p}{2}{4} + \twij{3}{p}{3}{4} = A_2+A_3, \; \; \;  \; \; 
 \twij{1}{q}{1}{2} \mapsto  \twij{1}{q}{1}{2}+ \twij{1}{p}{1}{3} =  B_2+B_3  $$
 $$\qand \twij{2}{r}{2}{3} \mapsto  \twij{2}{p}{2}{4} + \twij{3}{p}{3}{4} = C_2+C_3.$$
Note that $[A_i, [B_2, C_3]] = [A_i, [B_3, C_2]] = 0$ because $[w_{ij}, w_{kl}]=0$ when $\{i,j\} \cap \{ k, l \} = \phi$
 \\[0.2cm]The second generator of (C) included into $p_2=p_3$ maps to the element in \eqref{eqn.generatorC} after setting the above terms 0.
 \begin{equation}\label{eqn.generatorC}
 [\twij{2}{p}{2}{3}, [\twij{1}{q}{1}{2}, \twij{2}{r}{2}{3}]]  \mapsto  
 \end{equation}
\[
\left(\begin{array}{cl}
 &[A_2, [B_2, C_2]] 
 \\ + & [A_3, [B_2, C_2]]
 \\ + & [A_2, [B_3, C_3]]
 \\ + & [A_3, [B_3, C_3]]
 \end{array}\right)
=
\left(\begin{array}{cl}
&[\twij{2}{p}{2}{4}, [\twij{1}{q}{1}{2}, \twij{2}{r}{2}{4}]] 
\\ + &  [\twij{3}{p}{3}{4}, [\twij{1}{q}{1}{2}, \twij{2}{r}{2}{4}]] 
\\ + &  [\twij{2}{p}{2}{4}, [\twij{1}{q}{1}{3}, \twij{3}{r}{3}{4}]] 
\\ + &  [\twij{3}{p}{3}{4}, [\twij{1}{q}{1}{3}, \twij{3}{r}{3}{4}]] 
\end{array} \right)
\]

The last term is 0 from \eqref{eqn.134is0} in face $p_1=p_2$. The first term will be rewritten by \eqref{eqn: Break3}. So we get the 4-term equation in \eqref{eqn: 4Term}.
\begin{equation}\label{eqn: 4Term}
\begin{array}{l}
\;\;\;\;\left(\begin{array}{l}
+(t_2^{-q} t_3^ {-p-q} t_4^{-q-r}) [w_{14}, [w_{13}, w_{12}]]
\\+(t_2^{-q} t_3^ {-q-r} t_4^{-p-q}) [w_{13}, [w_{14}, w_{12}]]
\\+[\twij{3}{p}{3}{4}, [\twij{1}{q}{1}{2}, \twij{2}{r}{2}{4}]] 
\\+ [\twij{2}{p}{2}{4}, [\twij{1}{q}{1}{3}, \twij{3}{r}{3}{4}]] 
\end{array}\right)

 \\ \\ = \left(\begin{array}{l}
 +(t_2^{-q} t_3^ {-p-q} t_4^{-q-r}) [w_{14}, [w_{13}, w_{12}]] 
 \\+ (t_2^{-q} t_3^ {-q-r} t_4^{-p-q}) [w_{13}, [w_{14}, w_{12}]] 
 \\+[t_4^{-q-r}\twij{3}{p-q-r}{3}{4},t_2^{-q} t_4^{-q-r}[w_{12}, w_{24}]] 
 \\+ [t_4^{-q-r}\twij{2}{p-q-r}{2}{4},t_3^{-q} t_4^{-q-r} [w_{13}, w_{34}]] 
 \end{array}\right)
 \\
\\= 
\left(\begin{array}{r}
+(t_2^{-q} t_3^ {-p-q} t_4^{-q-r}) [w_{14}, [w_{13}, w_{12}]] 
\\+ (t_2^{-q} t_3^ {-q-r} t_4^{-p-q}) [w_{13}, [w_{14}, w_{12}]] 
\\+(t_2^{-q} t_3^{p-q-r} t_4^{-q-r})[w_{34}, [w_{12}, w_{24}]] 
\\+(t_2^{p-q-r} t_3^{-q} t_4^{-q-r}) [w_{24} ,[w_{13}, w_{34}]]
\end{array}\right)
\\ \\= \left(\begin{array}{l}
+(t_2^{-q} t_3^ {-p-q} t_4^{-q-r}) [w_{14}, [w_{13}, w_{12}]] 
\\+ (t_2^{-q} t_3^ {-q-r} t_4^{-p-q}) [w_{13}, [w_{14}, w_{12}]] 
\\+(t_2^{-q} t_3^{p-q-r} t_4^{-q-r})[w_{34}, [w_{14}, w_{12}]] 
\\+(t_2^{p-q-r} t_3^{-q} t_4^{-q-r}) [w_{24} ,[w_{14}, w_{13}]]
 \end{array}\right)
  \\
\\= 
\left(\begin{array}{r}
+(t_2^{-q} t_3^ {-p-q} t_4^{-q-r}) [w_{14}, [w_{13}, w_{12}]] 
\\+ (t_2^{-q} t_3^ {-q-r} t_4^{-p-q}) [w_{13}, [w_{14}, w_{12}]] 
\\-(t_2^{-q} t_3^{p-q-r} t_4^{-q-r})[w_{12}, [w_{34}, w_{14}]] 
\\-(t_2^{p-q-r} t_3^{-q} t_4^{-q-r}) [w_{13} ,[w_{24}, w_{14}]]
\end{array}\right)
\end{array}
\end{equation}
\[
= \left(\begin{array}{l}
+(t_2^{-q} t_3^ {-p-q} t_4^{-q-r}) [w_{14}, [w_{13}, w_{12}]] 
\\+ (t_2^{-q} t_3^ {-q-r} t_4^{-p-q}) [w_{13}, [w_{14}, w_{12}]] 
\\-(t_2^{-q} t_3^{p-q-r} t_4^{-q-r})[w_{12}, [w_{14}, w_{13}]] 
\\-(t_2^{p-q-r} t_3^{-q} t_4^{-q-r}) [w_{13} ,[w_{14}, w_{12}]] 
 \end{array}\right)
\]

Setting $p=q=r=0$ in \eqref{eqn: 4Term} we get that 
\begin{equation}\label{eqn: 2of3equal}
[w_{14}, [w_{13}, w_{12}]] - [w_{12}, [w_{14}, w_{13}]]  = 0
\end{equation}
Using the Jacobi identity, we further get 
$$[w_{13}, [w_{12}, w_{14}]]= -2  [w_{14}, [w_{13}, w_{12}]] = -2  [w_{12}, [w_{14}, w_{13}]]  = -[w_{13}, [w_{14}, w_{12}]] $$
\end{itemize}

Our four term relation then becomes 
\[
\left(
t_2^{-q} t_3^ {-p-q} t_4^{-q-r} + 2 \; t_2^{-q} t_3^ {-q-r} t_4^{-p-q} -t_2^{-q} t_3^{p-q-r} t_4^{-q-r} - 2\; t_2^{p-q-r} t_3^{-q} t_4^{-q-r}
 \right)
 [w_{14}, [w_{13}, w_{12}]]
\]
\vspace{-1.5cm}
\begin{equation}\label{eqn: 4Term2}
=0
\end{equation}
which says that  $\pi_7(\CSB{4})/R$ maps surjectively to $$G\coloneq \sfrac{\mathbb{Q}[t_2^{\pm 1}, t_3^{\pm 1}, t_4^{\pm 1}]}{\left((t_2^{-q} t_3^ {-p-q} t_4^{-q-r} + 2 \; t_2^{-q} t_3^ {-q-r} t_4^{-p-q} -t_2^{-q} t_3^{p-q-r} t_4^{-q-r} - 2\; t_2^{p-q-r} t_3^{-q} t_4^{-q-r})=0\right)}.$$

When we set $p=q=r$ in the 4-term relation \eqref{eqn: 4Term2}, it becomes 
\[
(3 \; t_2^{-q} t_3^ {-2q} t_4^{-2q} - 3\;t_2^{-q} t_3^ {-q} t_4^{-2q} )[w_{14}, [w_{13}, w_{12}]]=0,
\]
which gives us $t_3^{-q}=1$ in G (as we are working with rational homotopy groups). We don't consider $t_i=0$ because each of these are invertible. Using $t_3=1$ and setting $q=r=0$ we get $1+2t_4^{-p} = 1 + 2 t_2^{p}$ which gives $t_2 = t^{-1}_4$ in G.

So, at this point we have $$\pi_7\CSB{4}/R \twoheadrightarrow \mathbb{Q}[t_2^{\pm 1}]$$
where $1 \in \mathbb{Q} =  [w_{14}, [w_{13}, w_{12}]]$

\subsubsection{Relations with three indices}

We will now see how there is one more relation we get when including generator $(E)$ into face $p_1=p_2$. This gives us a relation with 3-indices containing $[\twij{1}{p}{1}{2}, [\twij{1}{q}{1}{4}, \twij{2}{r}{2}{4}]]$ and its two cyclic shifts which we then rewrite in terms of $[w_{14}, [w_{13}, w_{12}]]$ (and similarly for the relations containing $[\twij{1}{p}{1}{3}, [\twij{1}{q}{1}{4}, \twij{3}{r}{3}{4}]]$ and its two cyclic shifts).

We will use frequently that
$$[w_{14}, w_{24}] = [w_{12}, w_{14}] = [w_{24}, w_{12}].$$

First note that from \cref{eqn: 3Term0b}, we get that $$[\twij{1}{p}{1}{2}, [\twij{1}{q}{1}{4}, \twij{2}{r}{2}{4}]] = - [\twij{1}{p}{1}{2}, [t_1^q\twij{2}{r}{1}{2}, \twij{2}{r}{2}{4}]] = 0$$
This means $$[\twij{2}{r}{2}{4}, [\twij{1}{p}{1}{2}, \twij{1}{q}{1}{4}]]  = -  [\twij{1}{q}{1}{4}, [ \twij{2}{r}{2}{4}, \twij{1}{p}{1}{2}]] =  [\twij{1}{q}{1}{4}, [ \twij{1}{p}{1}{2}, \twij{2}{r}{2}{4}]] $$

In \eqref{eqn: simple 24 12 14} we use $[w_{ij}, w_{jk}] = [w_{jk}, w_{ki}]$ in the first equality and \eqref{eqn: Break3} for the second equality and use results from the end of the previous section to further simplify.

\begin{equation}\label{eqn: simple 24 12 14}
\begin{array}{ll}
        &  [\twij{2}{p}{2}{4}, [\twij{1}{q}{1}{2}, \twij{1}{q+r}{1}{4}]]  
    \\ =& -[\twij{2}{p}{2}{4}, [\twij{1}{q}{1}{2}, \twij{2}{r}{2}{4}]]  
    \\ = & -(t_2^{-q} t_3^ {-p-q} t_4^{-q-r}) [w_{14}, [w_{13}, w_{12}]] -(t_2^{-q} t_3^ {-q-r} t_4^{-p-q}) [w_{13}, [w_{14}, w_{12}]]
    \\ = & -t_2^r [w_{14}, [w_{13}, w_{12}]] -2 t_2^p[w_{14}, [w_{13}, w_{12}]]
    \\ = & -(t_2^r + 2t_2^p) [w_{14}, [w_{13}, w_{12}]] 
\end{array}
\end{equation}

We also have 
\begin{equation}\label{eqn: simple 14 12 24}
\begin{array}{ll}
        &  [\twij{1}{q+r}{1}{4}, [\twij{1}{q}{1}{2}, \twij{2}{p}{2}{4}]]   
    \\ =& [\twij{2}{p}{2}{4}, [\twij{1}{q}{1}{2}, \twij{1}{q+r}{1}{4}]]
    \\ = & -(t_2^r + 2t_2^p) [w_{14}, [w_{13}, w_{12}]].
\end{array}
\end{equation}

Generator (E) included into face $p_1=p_2$: 
\[
[\twij{1}{p}{1}{3}, [\twij{1}{q}{1}{3}, \twij{1}{r}{1}{3}]]  \mapsto 
\]
\[
\left(\begin{array}{cll}
&[\twij{1}{p}{1}{4}, [\twij{1}{q}{1}{4}, \twij{1}{r}{1}{4}]]  &+ [\twij{1}{p}{1}{4}, [\twij{1}{q}{1}{4}, \twij{2}{r}{2}{4}]]  
\\ + &  [\twij{1}{p}{1}{4}, [\twij{2}{q}{2}{4}, \twij{1}{r}{1}{4}]]  &+ [\twij{1}{p}{1}{4}, [\twij{2}{q}{2}{4}, \twij{2}{r}{2}{4}]] 
\\+ &[\twij{2}{p}{2}{4}, [\twij{1}{q}{1}{4}, \twij{1}{r}{1}{4}]]  &+ [\twij{2}{p}{2}{4}, [\twij{1}{q}{1}{4}, \twij{2}{r}{2}{4}]]  
\\ + &  [\twij{2}{p}{2}{4}, [\twij{2}{q}{2}{4}, \twij{1}{r}{1}{4}]]  &+ [\twij{2}{p}{2}{4}, [\twij{2}{q}{2}{4}, \twij{2}{r}{2}{4}]] 
\end{array} \right)
\]

The first term is 0 from generator (E) included into $p_2=p_3$ and the last term is 0 from $p_1=\ast$. We break up the 4th and 5th terms using the Jacobi relation and use relation $[w_{ij}, w_{jk}] = [w_{jk}, w_{ki}]$ repeatedly on all remaining terms to get
\[
[\twij{1}{p}{1}{3}, [\twij{1}{q}{1}{3}, \twij{1}{r}{1}{3}]]  \mapsto 
\]
\[
\left(\begin{array}{cll}
- &[\twij{1}{p}{1}{4}, [\twij{1}{q-r}{1}{2}, \twij{2}{r}{2}{4}]]  &+ [\twij{1}{p}{1}{4}, [\twij{1}{r-q}{1}{2}, \twij{2}{q}{2}{4}]]
\\ + &  [\twij{2}{q}{2}{4}, [\twij{1}{p-r}{1}{2}, \twij{1}{p}{1}{4}]]  & - [\twij{2}{r}{2}{4}, [\twij{1}{p-q}{1}{2}, \twij{1}{p}{1}{4}]] 
\\+ &[\twij{1}{q}{1}{4}, [\twij{1}{r-p}{1}{2}, \twij{2}{p}{2}{4}]]  & - [\twij{1}{r}{1}{4}, [\twij{1}{q-p}{1}{2}, \twij{2}{p}{2}{4}]]  
\\ + &  [\twij{2}{p}{2}{4}, [\twij{1}{q-r}{1}{2}, \twij{1}{q}{1}{4}]]  & - [\twij{2}{p}{2}{4}, [\twij{1}{r-q}{1}{2}, \twij{1}{r}{1}{4}]] 
\end{array} \right).
\]

Using \eqref{eqn: simple 24 12 14} and \eqref{eqn: simple 14 12 24} will greatly simplify calculations to make the last relation.

\[
\left(\begin{array}{ll}
    (t_2^{p-q+r} + 2t_2^{r}) &  -(t_2^{p+q-r} + 2t_2^q)  \\
    -(t_2^{r} + 2 t_2^{q}) & +(t_2^{q} + 2 t_2^{r})\\
    -(t_2^{p+q-r} + 2 t_2^{p}) & + (t_2^{r+p-q} + 2 t_2^{p}) \\
    -  (t_2^{r} + 2 t_2^{p}) &+  (t_2^{q} + 2 t_2^{p})
\end{array}\right)[w_{14}, [w_{13}, w_{12}]]
\]

Thus in $G$ this becomes
$$ 2t_2^{p-q+r} - 2t_2^{p+q-r} +2t_2^r - 2t_2^q = 0$$
which when we set $p=q=0$, we get $4t_2^r - 2 t_2^{-r} - 2 = 0 $ which setting $r=1$ is $$2t_2^{2} - 1 - t_2 = 0$$

If we set $r=-1$, we get
$$4t_2^{-1} - 2 t_2^{1} -2 = 0 \implies 2 -  t_2^{2} -t_2 = 0 $$
which if we subtract from the equation we got when setting $r=1$, this gives us $3t_2^2 -3 = 0$ which can be plugged back in to $2 -  t_2^{2} -t_2 = 0$ to get $t_2=1$.
 
and we finally get \cref{thm: pi7c4} that holds rationally: $$\pi_7\CSB{4}/R\coloneq \mathbb{Q} \qquad \text{generated by } [w_{12}, [w_{13}, w_{14}]]$$

\subsection{Computations in the Bousfield Kan spectral sequence for $\pi_\ast\EmbIM$}\label{subsec: relationsspectral}
We now turn our attention to the Bousfield Kan spectral sequence for $\pi_k(\EmbIM)$ which Sinha constructs in \cite{sinha} and using which Scannell and Sinha \cite{sinha-scannell} compute various differentials in the case of $M=B^4$. We recall that $$E^{-p, q}_1 = \bigcap ker(s_i) \subset \pi_qC'_{p}\lr{M} \text{ and  } d_1 = \sum (-1)^i\partial^i $$

We compute $d_1$ in our case of $M = S^1 \times B^3$ similar to \cite{sinha-scannell}. Here, the $E_1$ page has infinite dimensional cells as shown in \cref{fig: spectralE1} (by virtue of the $\pi_1$ action on homotopy groups of $C_k\lr{S^1\times B^3}$). 

\begin{figure}
\[\begin{array}{c|c|c|c|c|c|c|c|}
           &            &            &           &        &      &      & 10 \\
\hline
           &            & \Q^\infty     & \Q^\infty    & \Q^\infty&   \Q^\infty   &      & 9  \\
\hline
           &            &            &           &        &      &      & 8  \\
\hline
           &            &            & \Q^\infty    & \Q^\infty &  \Q^\infty    &      & 7  \\
\hline
           &            &            &           &        &      &      & 6  \\
\hline
           &            &            &           & \Q^\infty  & \Q^\infty     &      & 5  \\
\hline
           &            &            &           &        &      &      & 4  \\
\hline
           &            &            &           &        & \Q[t^{\pm1}]  & \Q& 3  \\
\hline
           &            &            &           &        &      &      & 2  \\
\hline
           &            &            &           &        & &  & 1  \\
\hline
  -7 & -6  & -5 & -4 \hspace{1ex} & -3 \hspace{1.5ex}  & -2 \hspace{1ex} &  -1 \hspace{1ex} & \\
\hline
\end{array}\]
\caption{$E_1$ page for spectral sequence computing $\EmbIS$}
\label{fig: spectralE1}
\end{figure}

For instance, $E^{-2, 5}_1$ is generated by $[\twij{1}{p}{1}{2}, \twij{1}{q}{1}{2}]$ where $p > q$ because this is the only non zero whitehead product and is also, trivially, in the kernals of all $s_i\colon \pi_5C_2\lr{M} \to \pi_5C_1\lr{M}$. 

We also have $E^{-3, 5}_1$ is generated by $[\twij{1}{p}{1}{3}, \twij{2}{q}{2}{3}]$ because any element in the intersection $\bigcap_i \text{ker }s_i\colon \pi_5C_3\lr{M} \to \pi_5C_2\lr{M}$ has to have all three indices involved, so that forgetting any point sends at least one $w_{ij}$ to 0.  

We will now compute the $\partial^i$ differentials from $E^{-2,5}_1 \to E^{-3,5}_1$.

\begin{align*}
\partial^0([\twij{1}{p}{1}{2}, \twij{1}{q}{1}{2}]) &= [\twij{2}{p}{2}{3}, \twij{2}{q}{2}{3}]
\\\\\partial^1([\twij{1}{p}{1}{2}, \twij{1}{q}{1}{2}]) &= [\twij{1}{p}{1}{3}+ \twij{2}{p}{2}{3}, \twij{1}{q}{1}{3}+ \twij{2}{q}{2}{3}] 
\\&= +[\twij{1}{p}{1}{3}, \twij{1}{q}{1}{3}] + [\twij{2}{p}{2}{3}, \twij{1}{q}{1}{3}]
\\& \quad+ [\twij{1}{p}{1}{3}, \twij{2}{q}{2}{3}] +[\twij{2}{p}{2}{3}, \twij{2}{q}{2}{3}]
\\\\\partial^2([\twij{1}{p}{1}{2}, \twij{1}{q}{1}{2}]) &= [\twij{1}{p}{1}{2}+ \twij{1}{p}{1}{3}, \twij{1}{q}{1}{2}+ \twij{1}{q}{1}{3}] \\&= +[\twij{1}{p}{1}{2}, \twij{1}{q}{1}{2}] + [\twij{1}{p}{1}{2}, \twij{1}{q}{1}{3}]
\\&  \quad + [\twij{1}{p}{1}{3}, \twij{1}{q}{1}{2}] +[\twij{1}{p}{1}{3}, \twij{1}{q}{1}{3}]
\\\\\partial^3([\twij{1}{p}{1}{2}, \twij{1}{q}{1}{2}]) &= [\twij{1}{p}{1}{2}, \twij{1}{q}{1}{2}]
\end{align*}

When we put these together into $d_1 = \sum\limits_{i=0}^3 (-1)^i\partial^i $, we get 

\begin{align*}
d_1([\twij{1}{p}{1}{2}, \twij{1}{q}{1}{2}]) &= -[\twij{2}{p}{2}{3}, \twij{1}{q}{1}{3}] - [\twij{1}{p}{1}{3}, \twij{2}{q}{2}{3}] 
\\ & \;\;+ [\twij{1}{p}{1}{2}, \twij{1}{q}{1}{3}]+ [\twij{1}{p}{1}{3}, \twij{1}{q}{1}{2}]
\\&= t_1^qt_2^p [w_{13}, w_{23}] - t_1^pt_2^q[w_{13}, w_{23}] 
\\ & \;\; + t_2^{-p}t_3^{-q}[w_{13}, w_{23}]+ t_3^{-p}t_2^{-q}[w_{23}, w_{13}]
\\& = (t_1^qt_2^p - t_1^pt_2^q  +  t_1^qt_2^{q-p}- t_1^p t_2^{p-q})[w_{13}, w_{23}]
\end{align*}
 which is precisely the hexagonal relation in Remark 3.5 in \cite{budney}. This makes $E_2^{-3, 5} = \sfrac{\Q[t_1^{\pm 1}, t_2^{\pm 1}]}{\lr{t_1^pt_2^q + t_1^p t_2^{p-q} = t_1^qt_2^p +  t_1^qt_2^{q-p} }}$.

We also note that terms we get in the images of $\partial^i$ that are not in $\bigcap ker\; s_i$ (in $\pi_5(C_3\lr{S^1\times B^3}$), like $[\twij{2}{p}{2}{3}, \twij{2}{q}{2}{3}]$, cancel out in the alternating sum to make a well defined $d_1$ to $E^{-3, 5}_1$.  


To show $E_2^{-2, 5} = 0$, we see the kernal of $d_1: E_1^{-2,5} \to E_1^{-3,5}$ is trivial (Rather, $[\twij{1}{p}{1}{2}, \twij{1}{q}{1}{2}] - [\twij{1}{q}{1}{2}, \twij{1}{p}{1}{2}]$ is in the kernal but is already 0). Hence the $E_2$ page in cells row 5 and lower looks like \cref{fig: spectralE2part}.

\begin{figure}
\[\begin{array}{c|c|c|c|c|c|c|c|}
           &            &            &           &        &      &      & 6  \\
\hline
           &            &            &           & \sfrac{\Q[t_1^{\pm 1}, t_2^{\pm 1}]}{\lr{t_1^pt_2^q + t_1^p t_2^{p-q} = t_1^qt_2^p +  t_1^qt_2^{q-p} }}  &      &      & 5  \\
\hline
           &            &            &           &        &      &      & 4  \\
\hline
           &            &            &           &        & \sfrac{\Q[t^{\pm1}]}{\lr{t^0}}& & 3  \\
\hline
  -7 & -6  & -5 & -4 \hspace{1ex} & -3 \hspace{1.5ex}  & -2 \hspace{1ex} &  -1 \hspace{1ex} & \\
\hline
\end{array}\]
\caption{Part of $E_2$ page}
\label{fig: spectralE2part}
\end{figure}

Because there is no other $d_r$ that hits $E^{-3, 5}_r$ and also there is nothing else in the $-p+q=2$ diagonal of $E^{-p,q}$, we can say that rationally, $$\pi_2\EmbIS \cong \Q[t_1^{\pm 1}, t_2^{\pm 1}]/\lr{t_1^pt_2^q + t_1^p t_2^{p-q} = t_1^qt_2^p +  t_1^qt_2^{q-p} }.$$ In other words, the $W_3$ map in \cite{budney} is an isomorphism.  

The computation of $d_1\colon E^{-3, 7}_1 \to E^{-4, 7}_1$ is more computationally challenging, but most of the work has been done in \cref{subsec: relations} where we computed the images of various generators under the face inclusions (in the context of cosimplicial spaces here, they will be called coface maps).

First we must determine $E^{-3, 7}_1$ and $E^{-4, 7}_1$.  $E^{-3, 7}_1$ contains iterated whitehead products of $\twij{i}{\alpha}{i}{j}$ where all three indices are present, so $[\twij{1}{p}{1}{2}, [\twij{1}{q}{1}{2}, \twij{1}{r}{1}{3}]]$ and $[\twij{1}{p}{1}{3},[\twij{1}{q}{1}{2}, \twij{1}{r}{1}{3}] $ (and a cyclic shift of each) will be generators.  Similarly $E^{-4, 7}_1$ contains iterated whitehead products of $\twij{i}{\alpha}{i}{j}$ where all four indices are present, so $[\twij{1}{p}{1}{2}, [\twij{1}{q}{1}{3}, \twij{1}{r}{1}{4}]]$ (and a cyclic shift) will be generators. 

\[
\begin{array}{l}
\partial^0([\twij{1}{p}{1}{2}, [\twij{1}{q}{1}{2}, \twij{1}{r}{1}{3}]]) =[\twij{2}{p}{2}{3}, [\twij{2}{q}{2}{3}, \twij{2}{r}{2}{4}]]
\\\\
\;\;\;\;\partial^1([\twij{1}{p}{1}{2}, [\twij{1}{q}{1}{2}, \twij{1}{r}{1}{3}]]) 
\\= [\twij{1}{p}{1}{3}+ \twij{2}{p}{2}{3}, [\twij{1}{q}{1}{3}+ \twij{2}{q}{2}{3}, \twij{1}{r}{1}{4}+ \twij{2}{r}{2}{4}] ]
\\=[\twij{1}{p}{1}{3}, [\twij{1}{q}{1}{3}, \twij{1}{r}{1}{4}] ]+ [\twij{2}{p}{2}{3}, [\twij{1}{q}{1}{3}, \twij{1}{r}{1}{4}] ]
\\+\;[\twij{1}{p}{1}{3}, [\twij{2}{q}{2}{3}, \twij{1}{r}{1}{4}] ]+[ \twij{2}{p}{2}{3}, [\twij{2}{q}{2}{3}, \twij{1}{r}{1}{4}] ]
\\+\;[\twij{1}{p}{1}{3}, [\twij{1}{q}{1}{3}, \twij{2}{r}{2}{4}] ]+[\twij{2}{p}{2}{3}, [\twij{1}{q}{1}{3},  \twij{2}{r}{2}{4}] ]
\\+\;[\twij{1}{p}{1}{3}, [\twij{2}{q}{2}{3},\twij{2}{r}{2}{4}] ]+[ \twij{2}{p}{2}{3}, [\twij{2}{q}{2}{3},  \twij{2}{r}{2}{4}] ]
\\\\
\;\;\;\;\partial^2([\twij{1}{p}{1}{2}, [\twij{1}{q}{1}{2}, \twij{1}{r}{1}{3}]]) 
\\=[\twij{1}{p}{1}{2}+\twij{1}{p}{1}{3}, [\twij{1}{q}{1}{2}+\twij{1}{q}{1}{3}, \twij{1}{r}{1}{4}]]
\\=[\twij{1}{p}{1}{2}, [\twij{1}{q}{1}{2}, \twij{1}{r}{1}{4}]]+[\twij{1}{p}{1}{2}, [\twij{1}{q}{1}{3}, \twij{1}{r}{1}{4}]]
\\+\;[\twij{1}{p}{1}{3}, [\twij{1}{q}{1}{2}, \twij{1}{r}{1}{4}]]+[\twij{1}{p}{1}{3}, [\twij{1}{q}{1}{3}, \twij{1}{r}{1}{4}]]
\\\\
\;\;\;\;\partial^3([\twij{1}{p}{1}{2}, [\twij{1}{q}{1}{2}, \twij{1}{r}{1}{3}]]) [\twij{1}{p}{1}{2}, [\twij{1}{q}{1}{2}, \twij{1}{r}{1}{3}+ \twij{1}{r}{1}{4}]]
\\=[\twij{1}{p}{1}{2}, [\twij{1}{q}{1}{2}, \twij{1}{r}{1}{3}]]+[\twij{1}{p}{1}{2}, [\twij{1}{q}{1}{2}, \twij{1}{r}{1}{4}]]
\\\\\partial^4([\twij{1}{p}{1}{2}, [\twij{1}{q}{1}{2}, \twij{1}{r}{1}{3}]]) =  [\twij{1}{p}{1}{2}, [\twij{1}{q}{1}{2}, \twij{1}{r}{1}{3}]]
\end{array}
\]
When we put these together into $d_1 = \sum\limits_{i=0}^4 (-1)^i\partial^i $, we get

\begin{equation}\label{eqn: 4termspectral}
\begin{array}{rl}
& d_1([\twij{1}{p}{1}{2}, [\twij{1}{q}{1}{2}, \twij{1}{r}{1}{3}]]) \\
= & - [\twij{2}{p}{2}{3}, [\twij{1}{q}{1}{3}, \twij{1}{r}{1}{4}] ] - [\twij{1}{p}{1}{3}, [\twij{2}{q}{2}{3},\twij{2}{r}{2}{4}] ] \\
& + [\twij{1}{p}{1}{2}, [\twij{1}{q}{1}{3}, \twij{1}{r}{1}{4}]]+[\twij{1}{p}{1}{3}, [\twij{1}{q}{1}{2}, \twij{1}{r}{1}{4}]] \\\\
= & + [\twij{1}{r}{1}{4}, [\twij{2}{p}{2}{3}, \twij{1}{q}{1}{3}]] + [\twij{2}{r}{2}{4},[\twij{1}{p}{1}{3}, \twij{2}{q}{2}{3}]] \\
& + (t_2^{-p}t_3^{-q}t_4^{-r}) [w_{12}[w_{13}, w_{14}]] + t_2^{-q}t_3^{-p}t_4^{-r}[w_{13},[w_{12}, w_{14}]] \\\\
= & + [\twij{1}{r}{1}{4}, t_2^pt_1^q[w_{13}, w_{12}]]+ [\twij{2}{r}{2}{4}, t_1^pt_2^q[w_{12}, w_{13}]] \\
& + (t_2^{-p}t_3^{-q}t_4^{-r}) [w_{12}[w_{13}, w_{14}]] + t_2^{-q}t_3^{-p}t_4^{-r}[w_{13},[w_{12}, w_{14}]] \\\\
= & + t_2^{p-q}t_3^{-q}t_4^{-r}[w_{14}, [w_{13}, w_{12}]] + [\twij{2}{r}{2}{4}, [\twij{1}{p-q}{1}{2}, \twij{1}{p}{1}{3}]] \\
& + (t_2^{-p}t_3^{-q}t_4^{-r}) [w_{12}[w_{13}, w_{14}]] + t_2^{-q}t_3^{-p}t_4^{-r}[w_{13},[w_{12}, w_{14}]] \\\\
= & + t_2^{p-q}t_3^{-q}t_4^{-r}[w_{14}, [w_{13}, w_{12}]] - [\twij{1}{p}{1}{3} , [\twij{4}{-r}{2}{4}, \twij{1}{p-q}{1}{2}]] \\
& +(t_2^{-p}t_3^{-q}t_4^{-r}) [w_{12}[w_{13}, w_{14}]] + t_2^{-q}t_3^{-p}t_4^{-r}[w_{13},[w_{12}, w_{14}]] \\\\
= & + (t_2^{p-q}t_3^{-q}t_4^{-r})[w_{14}, [w_{13}, w_{12}]] - (t_2^{q-p}t_3^{-p}t_4^{q-p-r})[w_{13}, [w_{12}, w_{14}]]] \\
& + (t_2^{-p}t_3^{-q}t_4^{-r}) [w_{12}[w_{13}, w_{14}]] + t_2^{-q}t_3^{-p}t_4^{-r}[w_{13},[w_{12}, w_{14}]]
\end{array}
\end{equation}

When we set $p=q=r=0$ in \eqref{eqn: 4termspectral} and set that expression to 0 in $E^{-4, 7}_1$ we get 
$$[w_{12}[w_{13}, w_{14}]] = -[w_{14}, [w_{13}, w_{12}]]= [w_{14}[w_{12}, w_{13}]]$$
which is the same relation we obtained in \eqref{eqn: 2of3equal}.  So we also have 
$$ [w_{13}, [w_{12}, w_{14}]] = 2 [w_{12}, [w_{13}, w_{14}]]$$
This also appears in $d_1
[[w_{13}, w_{23}], w_{23}]$ in  \cite{sinha-scannell} for $M=B^4$. Hence \eqref{eqn: 4termspectral} becomes
$$\left(-t_2^{p-q}t_3^{-q}t_4^{-r} - 2t_2^{q-p}t_3^{-p}t_4^{q-p-r}+t_2^{-p}t_3^{-q}t_4^{-r} +2t_2^{-q}t_3^{-p}t_4^{-r}\right)[w_{12}, [w_{13}, w_{14}]] = 0$$

Setting $p=q=r$, we get 
$$0=-t_3^{-p}t_4^{-p} - 2t_3^{-p}t_4^{-p} + t_2^{-p}t_3^{-p}t_4^{-p} + 2t_2^{-p}t_3^{-p}t_4^{-p}=-3(t_3t_4)^{-p}(1-t_2^{-p})$$

This gives us $t_2=1$ and setting $r=0$ in \eqref{eqn: 4termspectral} gives us 
$$0 = -t_3^{-q}-2t_3^{-p}t_4^{q-p}+t_3^{-q}+2t_3^{-p} = -2t_3^{-p}(t_4^{q-p}-1) = 0$$

which gives us $t_4=1$, in addition to $t_2=1$ we got previously. 

We  compute $d_1([\twij{1}{p}{1}{3}, [\twij{1}{q}{1}{2}, \twij{1}{r}{1}{2}]])$ (the cyclic shift of $[\twij{1}{p}{1}{2}, [\twij{1}{q}{1}{2}, \twij{1}{r}{1}{3}]]$). A similar calculation gives us 

\begin{equation}\label{eqn: spectralsecgen}
\begin{array}{rl}
& d_1([\twij{1}{p}{1}{3}, [\twij{1}{q}{1}{2}, \twij{1}{r}{1}{2}]]) \\
= & - [\twij{1}{p}{1}{4}, [\twij{1}{r}{1}{3}, \twij{2}{r}{2}{3}] ] - [\twij{1}{p}{1}{4}, [\twij{2}{q}{2}{3},\twij{1}{r}{1}{3}] ] 
\\ & - [\twij{2}{p}{2}{4}, [\twij{1}{q}{1}{3}, \twij{2}{r}{2}{3}] ] - [\twij{2}{p}{2}{4}, [\twij{2}{q}{2}{3},\twij{1}{r}{1}{3}] ]
\\ & + [\twij{1}{p}{1}{4}, [\twij{1}{q}{1}{2}, \twij{1}{r}{1}{3}] ] + [\twij{1}{p}{1}{4}, [\twij{1}{q}{1}{3},\twij{1}{r}{1}{2}] ]
\\\\ = & -[\twij{1}{p}{1}{4}, t_1^{q}t_2^{r}[w_{12}, w_{13}]]-[\twij{1}{p}{1}{4}, t_2^{q}t_1^{r}[w_{13}, w_{12}]] 
\\ & - [\twij{2}{p}{2}{4}, t_1^{q}t_2^{r}[w_{12}, w_{13}]]-[\twij{2}{p}{2}{4}, t_2^{q}t_1^{r}[w_{13}, w_{12}]] 
\\ & + t_2^{-q}t_3^{-r}t_4^{-p}[w_{14}, [w_{12}, w_{13}]] + t_2^{-r}t_3^{-q}t_4^{-p}[w_{14},[w_{13},w_{12}]]
\\\\ = & - t_2^{r-q}t_3^{-q}t_4^{-p}[w_{14}, [w_{12}, w_{13}]] - t_2^{q-r}t_3^{-r}t_4^{-p}[w_{14}, [w_{13},w_{12}]]
\\ & - t_2^{r-q}t_3^{-q}t_4^{r-q-p}[w_{24}, [w_{12}, w_{13}]] - t_2^{q-r}t_3^{-r}t_4^{q-r-p}[w_{24}, [w_{13},w_{12}]]
\\ & + t_2^{-q}t_3^{-r}t_4^{-p}[w_{14}, [w_{12}, w_{13}]] + t_2^{-r}t_3^{-q}t_4^{-p}[w_{14},[w_{13},w_{12}]]
\\\\ = & - t_2^{r-q}t_3^{-q}t_4^{-p}[w_{14}, [w_{12}, w_{13}]] - t_2^{q-r}t_3^{-r}t_4^{-p}[w_{14}, [w_{13},w_{12}]]
\\ & + t_2^{r-q}t_3^{-q}t_4^{r-q-p}[w_{13}, [w_{12}, w_{14}]] + t_2^{q-r}t_3^{-r}t_4^{q-r-p}[w_{13}, [w_{14},w_{12}]]
\\ & + t_2^{-q}t_3^{-r}t_4^{-p}[w_{14}, [w_{12}, w_{13}]] + t_2^{-r}t_3^{-q}t_4^{-p}[w_{14},[w_{13},w_{12}]]
\\\\ = & (-t_3^{-q}+t_3^{-r}+2t_3^{-q}-2t_3^{-r}+t_3^{-r}-t_3^{-q})[w_{12},[w_{13}, w_{14}]] 
\\\\ = & 0
\end{array}
\end{equation}

We  now have to compute and $d_1([\twij{1}{p}{1}{3}, [\twij{1}{q}{1}{2}, \twij{1}{r}{1}{3}]])$, and for its cyclic shift. We get the following by a similar alternating sum.

\[
\begin{array}{rl}
& d_1([\twij{1}{p}{1}{3}, [\twij{1}{q}{1}{2}, \twij{1}{r}{1}{3}]])
\\ = & - [\twij{1}{p}{1}{4}, [\twij{2}{q}{2}{3}, \twij{2}{r}{2}{4}]] - [\twij{2}{p}{2}{4}, [\twij{1}{q}{1}{3}, \twij{1}{r}{1}{4}]]
\\ & - [\twij{1}{p}{1}{3}, [\twij{1}{q}{1}{2}, \twij{1}{r}{1}{4}]] - [\twij{1}{p}{1}{4}, [\twij{1}{q}{1}{2}, \twij{1}{r}{1}{3}]]
\\\\ = & - t_1^{p-r}t_3^{-q}t_4^{-r}[w_{14}, [w_{23}, w_{24}]] - t_2^{p-r}t_3^{-q}t_4^{-r}[w_{24}, [w_{13}, w_{14}]
\\ & - t_2^{-q}t_3^{-p}t_4^{-r}[w_{13},[w_{12}, w_{14}]] - t_2^{-q} t_3^{-r}t_4^{-p}[w_{14}, [w_{12}, w_{13}]]
\\\\ = & + t_1^{p-r}t_3^{-q}t_4^{-r}[w_{23}, [w_{24}, w_{14}]] + t_2^{p-r}t_3^{-q}t_4^{-r}[w_{13},[ w_{14}, w_{24}]
\\ & - ( 2t_3^{-p} + t_3^{-r})[w_{14}, [w_{12}, w_{13}]]
\\\\ = & + t_1^{p-r}t_3^{-q}t_4^{-r}[w_{23}, [w_{14}, w_{12}]] + t_2^{p-r}t_3^{-q}t_4^{-r}[w_{13},[ w_{12}, w_{14}]
\\ & - ( 2t_3^{-p} + t_3^{-r})[w_{12}, [w_{13}, w_{14}]]
\\\\ = & -  t_1^{p-r}t_3^{-q}t_4^{-r}[w_{14}, [w_{12}, w_{23}]]+(2t_3^{-q}- 2t_3^{-p} - t_3^{-r})[w_{12}, [w_{13}, w_{14}]]
\\\\ = & - [t_4^{-p}w_{14}, [t_3^{r-q-p}w_{13}, t_2^{r-p}w_{12}]]+(2t_3^{-q}- 2t_3^{-p} - t_3^{-r})[w_{12}, [w_{13}, w_{14}]]
\\\\ = & (t_3^{r-q-p} + 2t_3^{-q}- 2t_3^{-p} - t_3^{-r})[w_{12}, [w_{13}, w_{14}]]
\end{array}
\]

When we set $q=-1$ and $p=r=0$, we get $3t_3^{1}-3 = 0$ giving us $t_3=1$. 

Similar to the cyclic shift of $[\twij{1}{p}{1}{2}, [\twij{1}{q}{1}{2}, \twij{1}{r}{1}{3}]]$, we can show $d_1([\twij{1}{p}{1}{2}, [\twij{1}{q}{1}{3}, \twij{1}{r}{1}{3}]])$ equals 0. Thus we have reduced $E^{-4,7}/\text{im }d_1$ to $\Q $, where $1\in \Q$ corresponds to $[w_{12}, [w_{13}, w_{14}]]$. Thus the $E_2$ page only has $E^{-4, 7}_2 = \Q$ on the $-p+q = 3$ diagonal which proves \cref{thm: pi3EmbS1B3}. 

\section{Strategies for Showing $\Gpqr$ is Nontrivial}\label{sec.progress}

\subsection{Detecting whitehead products: Linking numbers and Hopf invariants}

Whitehead products $[f,g]: S^{n+m-1} \to X$ factor as $S^{n+m-1} \to S^n\vee S^m\xrightarrow{f\;\vee\; g} X$ where the first map is the whitehead product of inclusions of $S^n$ and $S^m$ into their wedge. We denote that map $\phi: S^{n+m-1} \to S^n \vee S^m$. If $a\in S^n , b\in S^m$ are non wedge points, then $\phi^{-1}(a)$ and $\phi^{-1}(b)$ are homeomorphic to $S^{m-1}$ and $S^{n-1}$ that are linked as a generalized Hopf link. We can use this idea to define an invariant of homotopy classes of maps $f: S^{n+m-1} \to S^{n}\vee S^{m}$ as the linking number between $f^{-1}(a)$ and $f^{-1}(b)$. 

Sinha and Walter in \cite{sinha-walter} describe a theory of Hopf invariants to detect homotopy groups. We describe here that theory applied to the special case above (See \cite[Example 1.9]{sinha-walter}). Suppose $A$ and $B$ are disjoint submanifolds of a manifold $X$ with co-dimensions $d_A$ and $d_B$ respectively, and we want to create a homotopy invariant of a map $f: S^{d_A + d_B - 1} \to X$. Let $\omega_A$ (and similarly $\omega_B)$ denote a representative of a $d_A$ dimensional \emph{Thom cochain} dual to $A$ (dual in the sense of the cap product). Then the linking number of $f^{-1}(A)$ and $f^{-1}(B)$ in $S^{(d_A+d_B-1)}$ is the same as an evaluation of a certain top dimensional cohomology class on $[S^{(d_A+d_B-1)}]$. The invariant turns out to be: 

$$\int\limits_{S^{(d_A+d_B-1)}} (d^{-1}f^\ast \omega_A \wedge f^\ast\omega_B)$$
where $d^{-1}\omega$ picks out some representative of $H^{dim(\omega)-1}( \cdot)$ that the coboundary $d$ maps to the cochain $\omega$. This would also be equal to $\int\limits_{S^{(d_A+d_B-1)}} (f^\ast \omega_A \wedge d^{-1}f^\ast\omega_B)$. The analogy is that linking number between $X_0, X_1\in X$ is the intersection number of $X_i$ with a manifold that bounds $X_{(1-i)}$. 

\begin{figure}
    \centering
    \includegraphics[width=8cm]{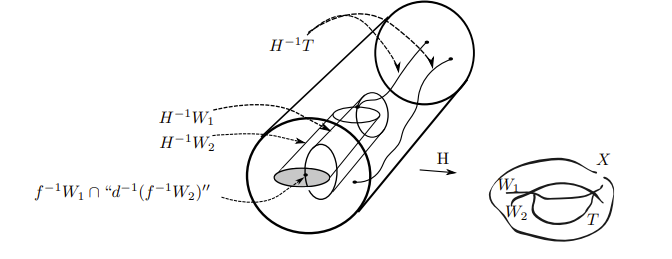}
    \caption{Depiction of linking with correction from \cite{sinha-walter}}
    \label{fig: linkwcorrect}
\end{figure}

Now if $A$ and $B$ intersect, then they describe a generalized linking invariant with correction to be 

$$\int\limits_{S^n} (d^{-1}f^\ast \omega_A \wedge f^\ast\omega_B) \pm f^\ast\omega_{(A\cap B)}$$

where $\omega_{(A\cap B)}$ is the Thom co-chain for $A\cap B$. This allows us to detect a homotopy class of $f$ that may have representatives that have intersections between $f^{-1}(A)$ and $f^{-1}(B)$, but we we keep track of those intersection points with sign. This is best described in \cite[Figure 2]{sinha-walter} which is copied here as \cref{fig: linkwcorrect} for convenience.

Now we will describe Budney and Gabai's linking invariant to detect elements of \\$\pi_5\CSB{3}$ induced by $G(p,q)$. Let $\tilde{C}_{k}(S^1\times B^3)$ denote the universal covering space of $\CSB{k}$ that is seen as a subset of $C_k(\mathbb{R}\times B^3)$ where each point has $\mathbb{Z}$ orbits. \cite{budney} define $t^\alpha Co_i^j \subset \CSB{k}$ to be the subspace of points $(p_1, p_2, \cdots p_k)$ such that $t^\alpha p_j - p_i$ in $\tilde{C}_{k}(S^1\times B^3)$ is parallel to a chosen vector $\zeta$. (Here $t^\alpha p_j$ denotes the endpoint of the lift of the loop $\alpha$ based at $p_j$). These $t^\alpha Co_i^j$ detect $t_j^\alpha w_{ij}$.

Let $A$ be $t^p Co^3_1$ and $B$ be $t^q Co^3_2$. $A$ detects $t_1^p w_{13}$ and $B$ detects $t_2^q w_{23}$ (see \cite[Figure 9]{budney}). Let $lk(A,B)$ denote the linking number between $(G(p,q))^{-1}(A)$ and $(G(p,q))^{-1}(B)$. If $A$ and $B$ didn't intersect, $lk(A,B)$ would be the coefficient of $t_1^pt_2^q[w_{13}, w_{23}]$ that $G(p,q)$ maps to, however $A$ and $B$ do intersect. They make an appropriate correction to account for this intersection. Let $C = t^{p-q}Co^2_1$ and $D = t^{q-p}Co^1_2$. They show that $lk(A-D,B-C)$ is an invariant that detects the coefficient of $t_1^pt_2^q[w_{13}, w_{23}]$. In the next section (\cref{subsec: cohorinvt}), we show using ideas from \cite{sinha-walter} how the sum/difference of linking numbers of preimages of cohorizontal manifolds that Budney and Gabai use to detect $[t_1^pw_{13}, t_2^qw_{23}]$ is an invariant. This is essentially proved in \cite[Section 3.3]{sinha-walter}. 

\begin{figure}
    \centering
    \includegraphics[width=6cm, trim = {1cm 21cm 13cm 1cm}, clip]{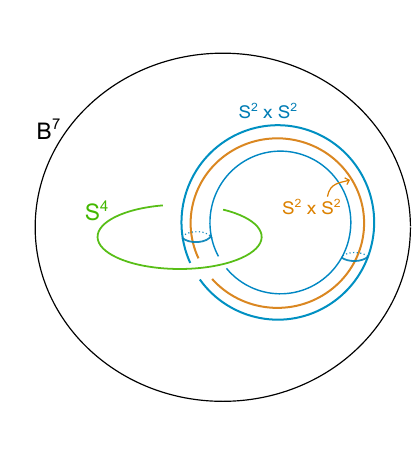}
    \caption{Schematic picture for linking in $S^7$ that detects $[[w_{14}, w_{24}], w_{34}]$}
    \label{fig: iterwhite}
\end{figure}

In the case of iterated Whitehead products (a focus of this thesis): $\phi: S^7 \to (S_a^3\vee S_b^3)\vee S_c^3$, factors as $S^7 \xrightarrow{\phi_1} S_p^5 \vee S_c^{3} \xrightarrow{\phi_2} (S_a^3 \vee S_b^3)\vee S_c^3$. For points $a\in S_a^3, b\in S_b^3, c\in S_c^3$, we will have $\phi^{-1}(c) = S^4$, and $\phi_1^{-1}(p)=S^2$. Inside $S_p^5$, we have $\phi_2^{-1}(a) = \phi_2^{-1}(b) = S^2$. So, $\phi^{-1}(a) = \phi^{-1}(b) = \phi_1^{-1}(p)\times \phi_2^{-1}(b) = S^2\times S^2$. So, to detect $[[w_{14}, w_{24}], w_{34}]$, we should expect the submanifolds that detect $w_{14}$ and $w_{24}$ would have preimages $S^2\times S^2$ under $G(0,0,0)$ and the submanifold that detects $w_{34}$ has preimage $S^4$ arranged in a triple linked configuration depicted schematically in \cref{fig: iterwhite}.

$Co^4_i$ detects $w_{i4}$ for $i=1, 2, 3$. However, these submanifolds intersect each other, furthermore, in its current state, $G(0,0,0)$ intersects $Co^4_i$ at these mutual intersections. One way to show directly that $G(0,0,0)\mapsto\pm [[w_{14}, w_{24}], w_{34}]$ would be to deform it to a map similar to the one in \cite[Definition 12.16]{budney} where the blue chords are only seen when $p_1\in I_1$ (and similarly for the green and orange chords). Another approach could be creating a well defined linking invariant using ideas that we describe in \cref{subsec: cohorinvt}. 

\subsection{Showing the cohorizontal intersection number is an invariant}\label{subsec: cohorinvt}

Recall that $A = t^p Co^3_1$, $B = t^q Co^3_2$, $C = t^{p-q}Co^2_1$ and $D = t^{q-p}Co^1_2$. We see $A\cap B = (p_1, p_2, p_3)$ such that $(t^p p_3 - p_1)$ is parallel to $(t^q p_3 - p_2)$ and both are parallel to $\zeta$. So either $(p_3, t^{-p} p_1,  t^{-q} p_2)$ are collinear along $\zeta$ (in that order) or $(p_3, t^{-q} p_2, t^{-p} p_1)$ are collinear along $\zeta$. We have $C\cap D = \phi$, $A\cap D$ consists of points $(p_3, t^{-p} p_1,  t^{-q} p_2)$ along $\zeta$ and $B\cap C$ consists of points $(p_3,  t^{-q} p_2, t^{-p} p_1)$ along $\zeta$ in those orders.

Hence $(A\cap B) \cup (C\cap D) = (A\cap C) \cup (B\cap D)$. As long as the map we are detecting does not intersect $A, B, C, D$ in any of their mutual intersections, this allows us to define an invariant as follows. We calculate the invariant as the linking number between the preimages of $A - C$ and $B - D$. In other words, we add $lk(A, B)$ and $lk(C,D)$ and subtract $lk(A,D)$ and $lk(B,C)$. One intuition for why this is an invariant is because $lk(A,B)$ detects the coefficient of $t_1^pt_2^q[w_{13}, w_{23}]$, while $-lk(A,D)$ detects $-[t_1^{p}w_{13}, t_2^{q-p}w_{12}]$ and $-lk(C,B)$ detects  $-[t_1^{p-q}w_{12}, t_2^{q}w_{23}]$ which are all homotopic Whitehead products. (The pair $(C, D)$ also detects a whitehead product but that one is 0). 

We use ideas from \cite{sinha-walter} to show that any submanifolds $A,B,C,D$ such that $(A\cap B) \cup (C\cap D) = (A\cap C) \cup (B\cap D)$ creates such an invariant. Recall that $$\int\limits_{S^n} (d^{-1}f^\ast \omega_A \wedge f^\ast\omega_B) \pm f^\ast\omega_{(A\cap B)}$$ is an invariant. However, if we add/subtract the linking numbers in the Budney-Gabai invariant, the second term cancels out and we are left with 
\begin{align*}
&+ \int\limits_{S^n} d^{-1}f^\ast (\omega_A ) \wedge f^\ast(\omega_B) -\int\limits_{S^n} d^{-1}f^\ast (\omega_A) \wedge f^\ast(\omega_D) 
\\& +\int\limits_{S^n} d^{-1}f^\ast (\omega_C) \wedge f^\ast(\omega_D) -\int\limits_{S^n} d^{-1}f^\ast (\omega_C) \wedge f^\ast(\omega_B) 
\end{align*}

which would simply find the linking number (by geometry) of the preimages of the pair $(A-C, B-D)$. This would then be a homotopy invariant even though $A,B,C,D$ have pairwise intersections. 

In \cite{budney}, they proved this is an invariant by constructing a cobordism between certain collinear manifolds to the above mentioned sum/difference of cohorizontal manifolds. The alternative proof (nearly identical to \cite[Section 3.3]{sinha-walter}) presented here can hopefully be generalized to detect elements of $\pi_3\EmbIM$ like $\Gpqr$. We would need to find what combination of $Co_i^j$ manifolds would cancel out intersections to make a well defined generalized linking invariant. Furthermore, we would have to understand better how to compute the correction terms given that we would have 3 intersecting 4-dimensional manifolds (with possibly many components) in $S^7$.

\section{Future Goals}\label{sec.goals}

We describe here some future ambitions of this research project.

\subsection{Further computations in the spectral sequence for $\pi_\ast\EmbIM$}

We have computed differentials in the spectral sequence from \cref{subsec: cosimplicial}. We would like to be able to compute higher homotopy groups and potentially have a general result for the homotopy groups of $\EmbIS$.

Another curious fact is that $E^{-4,7}_1/ker(d_1(E^{-3,7}_1))$ is isomorphic to $\pi_7(\CSB{4})/R$ where $R$ is the images of the 5 possible face inclusions. The generators and relations of the former are a strict subset of the generators and relations of the latter. However in both situations of $\pi_7\CSB{4}/R$ and $\pi_5\CSB{3}$, they have been isomorphic to the corresponding groups from the spectral sequence. It is plausible this holds for higher dimensional groups $\pi_{2n-1}\CSB{n}/R$. 

\subsection{Showing $G(0,0,0)$ generates $\pi_3\EmbIS$}

The immediate next step of this thesis would be to show that $G(0,0,0)$ is the (rational) generator of $\EmbIS$. Some strategies and challenges to this were described in \cref{sec.progress}. 

\subsection{Develop linking/intersection invariants to detect elements of $\pi_{2n+1}C_n$}

We would like to create well defined linking invariants to detect higher degree iterated whitehead products like $[[A_1, [A_{2}, [\cdots A_{m_1}]]\cdots ]] , [A_{m_1+1}, \cdots A_{m_2}] \cdots ] $.

\subsection{Generalizing $G(p,q)$ to higher dimensions}
We constructed $G(p,q,r)$ by ``smashing" a null homotopic map (orange \& blue chords) $S^2\to \EmbIM$ and a null homotopic map (green chords) $S^1\to \EmbIM$. 

We can generalize this to construct maps $S^4\to \EmbIM$ by using 2 different null homotopic maps $S^2\to \EmbIM$ or a map $S^1\to\EmbIM$ smashed with a map $S^3\to \EmbIM$. (We could expect some relations between these two constructions given that $[A, [B, [C, D]]] + [B, [[C,D], A]] + [[C,D], [A, B]] = 0$ on the Whitehead product side). We can further generalize these to higher homotopy groups of $\EmbIM$ and the question to ask would be if these are the generators of those groups.

\begin{figure}
    \centering
    \includegraphics[width=12cm]{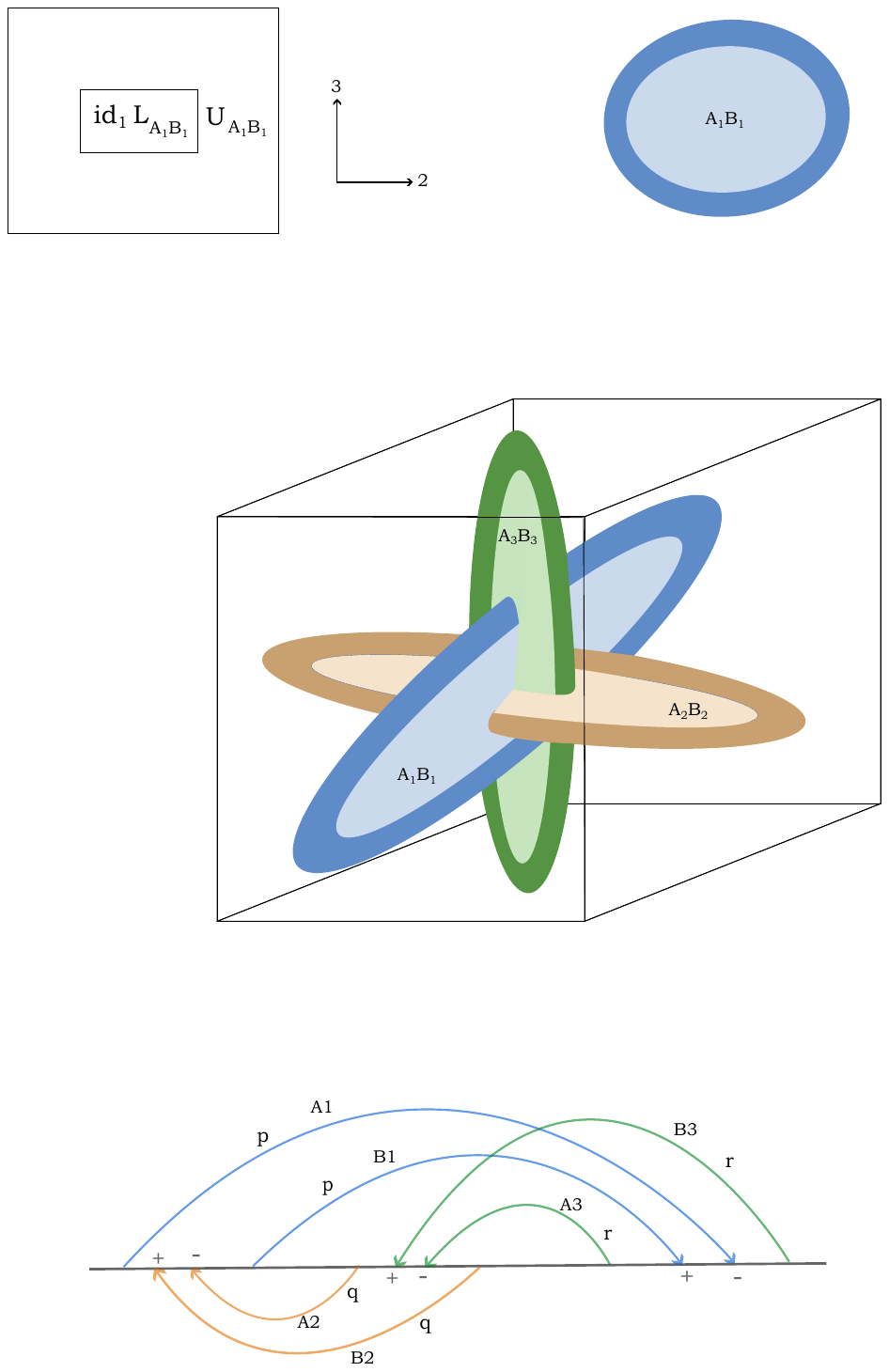}
    \caption{$E(p,q,r)$ construction}
    \label{fig: epqr}
\end{figure}

\subsection{$E(p,q,r)$}

We can somewhat generalize the construction of $E(p,q)$ from \cite{budney} to a map $(I^3, \partial I^3) \to (\EmbIM, \gamma)$ as shown in \cref{fig: epqr}. This does not appear to be a linear combination of $\Gpqr$ unlike $E(p,q)$ which equals $-G(-q, p) + G(p, -q)$. However, this element is also null homotopic in $T_3\EmbIM$ by a similar but simpler argument as we did for $\Gpqr$ in \cref{sec.nullhomotopy} (because it only requires transitions between undo and backtrack homotopies), but it remains to be seen if this is non trivial in $\pi_3\EmbIM$.

Other properties of the equivalence classes (for analogues, see Lemma 2.24, Prop 2.28 in \cite{budney}) of $E(p,q,r)$ like independence of end homotopies, multilinearity (up to certain restrictions) also hold. 

\bibliography{bibliography}
\bibliographystyle{alpha}

\end{document}